\documentclass[11pt,twoside]{article}
\usepackage[hmargin=0.8in,vmargin=0.9in]{geometry}
\usepackage{fancyhdr}
\usepackage{graphicx}
\usepackage{subfigure}
\usepackage{amssymb}
\usepackage{amsmath}
\usepackage{amsfonts}
\usepackage{theorem}
\usepackage{mathrsfs}
\usepackage{mathtools}
\usepackage{bm}
\usepackage{color}
\usepackage{setspace}
\usepackage{exscale}
\usepackage{relsize}
\usepackage{epstopdf}
\usepackage{float}
\usepackage{cite}
\usepackage{hyperref}

\usepackage{diagbox}
\usepackage{array}

\usepackage{color}

\usepackage{booktabs,multirow} 
\usepackage{array} 
\usepackage{paralist} 
\usepackage{verbatim} 
\usepackage{subfigure} 

\theoremstyle{plain}			
\newtheorem{thm}{Theorem}[section]

{\theorembodyfont{\rmfamily}}

\newenvironment{DA}{{\flushleft \bf Declarations:}}{}

\allowdisplaybreaks[1]

\numberwithin{equation}{section}
\numberwithin{figure}{section}
\numberwithin{table}{section}

\newcommand\eref[1]{(\ref{#1})}

\newcommand*\xbar[1]{%
  \hbox{%
    \vbox{%
      \hrule height 0.5pt 
      \kern0.4ex
      \hbox{%
        \kern-0.05em
        \ensuremath{#1}%
        \kern-0.00em
      }%
    }%
  }%
}

\newcommand{\mf}{\bm{f}}
\newcommand{\mF}{\bm{F}}

\newcommand{\bmu}{\bm{u}}

\newcommand{\mo}{\bm{0}}

\newcommand{\dx}{\Delta x}

\newcommand{\hf}{{\frac{1}{2}}}

\newcommand{\jph}{{j+\frac{1}{2}}}
\newcommand{\jmh}{{j-\frac{1}{2}}}

\newtheorem{rmk}[thm]{Remark}

\def\softd{{\leavevmode\setbox1=\hbox{d}%
          \hbox to 1.05\wd1{d\kern-0.4ex{\char039}\hss}}}

\title{High-Order Schemes for Hyperbolic Conservation Laws Using Young Measures}
\author{Shaoshuai Chu\thanks{Department of Mathematics, RWTH Aachen University, 52056, Aachen, Germany;
{\tt chu@igpm.rwth-aachen.de}}~ and ~Michael Herty\thanks{Department of Mathematics, RWTH Aachen
University, 52056, Aachen, Germany;
Department of Mathematics and Applied Mathematics, University of Pretoria, Private Bag X20, Hatfield 0028, South Africa; {\tt herty@igpm.rwth-aachen.de}
}}
\date{}
\begin{document}
\date{}
\maketitle
\begin{abstract}
We develop high-order numerical schemes to solve random hyperbolic conservation laws using linear programming. The proposed schemes are high-order extensions of the existing first-order scheme introduced in [{\sc S. Chu, M. Herty,  M. Luk\'a\v{c}ov\'a-Medvi{\softd}ov\'a, and Y. Zhou}, SIAM J. Sci. Comput., 48 (2026)], where a novel structure-preserving numerical method using a concept of generalized, measure-valued solutions to solve random hyperbolic systems of conservation laws is proposed, yielding a linear partial differential equation concerning the Young measure and allowing the computation of approximations based on
linear programming problems. The second-order extension is obtained using piecewise linear reconstructions of the one-sided point values of the unknowns. The fifth-order scheme is developed using the finite-difference alternative weighted essentially non-oscillatory (A-WENO) framework. These extensions significantly improve the resolution of discontinuities, as demonstrated by a series of numerical experiments on both random (Burgers equation, isentropic Euler equations) and deterministic (discontinuous flux, pressureless gas dynamics, Burgers equation with non-atomic support) hyperbolic conservation laws.
\end{abstract}



\section{Introduction}
This paper focuses on numerical solutions of random hyperbolic conservation laws, which arise in various fields including engineering, atmospheric, and geophysical
applications. It is well known that these systems often develop complex wave patterns, including discontinuities such as shocks and rarefaction waves, even when the initial data are very smooth, which introduces significant challenges in designing suitable numerical methods for such systems.

In the last two decades, a series of numerical methods have been proposed for solving random conservation laws, including the non--intrusive (Quasi-) Monte Carlo simulations, stochastic collocation, and (intrusive) generalized polynomial chaos expansion; see, e.g., \cite{MR3674794,MR3793404,MR3471175,MR2199923,MR2723020,MR3202521,MR4045233,MR3328389,MR2501693,MR4066974,MR1507356} and the references therein. Monte Carlo methods are widely valued for their robustness, though they are not always computationally efficient. The stochastic collocation method benefits from straightforward applicability, as it allows existing deterministic numerical solvers to be reused with minimal adaptation. In contrast, intrusive approaches such as stochastic Galerkin methods, which rely on generalized polynomial chaos expansions (gPC), impose greater constraints. These methods require projecting the governing equations onto a gPC basis, generating an expanded deterministic system for the expansion coefficients \cite{MR2660316,MR3129545}. Although this intrusive approach can achieve high accuracy with low polynomial degrees for smooth solutions \cite{MR2823001,MR2855645}, it faces two key challenges. First, the resulting coupled system for gPC coefficients may lose hyperbolicity, complicating its numerical solution \cite{MR3202523,MR4228317,MR3968112,MR4066974,MR4772617,MR4699586}. Second, the method demands significant re--implementation of the deterministic solver, as the governing equations must be reformulated and projected onto the gPC basis, creating substantial implementation overhead. At the same time, both stochastic collocation and Galerkin methods are susceptible to Gibbs phenomena near discontinuities due to their spectral approximation frameworks. Mitigating these oscillations often requires additional techniques such as spectral filtering \cite{MR1874071,MR4045233,MR4228317,MR2333926,MR3874527}, spectral viscosity \cite{MR977947,MR1770237}, or adaptive diffusion \cite{chertockchallenges,MR2979844,CKK_24}. The advantages and disadvantages of different approaches have been extensively discussed in the literature, and we refer to \cite{chertockchallenges} for recent discussions.

In \cite{CHLZ25}, a novel numerical approach based on a measure-valued framework is proposed using parameterized Young measures  rather than traditional weak (distributional) solutions \cite{MR3409135,MR257325} (see also \cite{DiPerna1985}, where the concept of measure-valued solutions for conservation laws was introduced in a rigorous analytical framework and has inspired numerous developments).  This reformulation transforms the problem into a linear optimization problem, circumventing several challenges associated with existing methods. Inspired by recent advances in dissipative measure-valued solutions for compressible fluid flows, the approach directly computes the Young measures within the numerical scheme, unlike previous methods where it emerges from weakly convergent approximations. Additionally, while previous works (see, e.g.,\cite{MR4777936}) approximate Young-measure moments using multiple entropy functions—often unavailable for general hyperbolic systems—this method considers only moments in the random dimension, enabling standard finite-volume schemes for moment propagation. Compared with existing moment-based approaches to random conservation laws; see, e.g.,\cite{MR3828934, MR4045233}, it introduces a different closure formulation that eliminates the need for filters and extends naturally to systems of equations. The proposed method has been applied so far only to one-dimensional (1-D) hyperbolic conservation laws, including the isentropic Euler system and Burgers equation. In this paper, we extend the discussion to a series of further problems as well as high-order finite-volume (FV) and finite-difference (FD) methods. 

The approach introduced in \cite{CHLZ25}, which uses a FV Lax--Friedrichs fluxes, offers a robust and accurate tool for solving  random nonlinear conservation laws. However, it is only first-order accurate and requires very fine meshes to achieve the desired resolution for shock waves and other complex wave structures. In order to address this issue, in this paper, we adopt the semi-discretization form and extend the first-order approach to higher orders (second and fifth) of accuracy by employing high-order local Lax–Friedrichs (LLF) fluxes in the spatial dimension. The second-order extension is achieved by applying piecewise linear reconstruction to obtain the one-sided point values of the unknowns, and a fifth-order scheme is developed within the FD alternative weighted essentially non-oscillatory (A-WENO) framework. The latter has been proven to be a powerful tool for generalizing low-order FV schemes to higher-order FD ones; see, e.g., \cite{Liu17,CKX23,Jiang13,Jiang96}. In order to reduce  non--physical oscillations, we use the local characteristic decomposition (LCD) in the high-order interpolation as done in \cite{Jiang13,Liu17,Qiu02,Shu20,CK2023,CCHKL_22}. For evolving the numerical solution in time, we solve the system of ordinary differential equations (ODEs) using a standard three-stage third-order strong stability preserving Runge-Kutta (SSP RK3) solver; see, e.g., \cite{Gottlieb12,Gottlieb11}.

The studied higher-order schemes are applied to a series of numerical examples, including both random and deterministic hyperbolic conservation laws. In the random case, the obtained results of Burgers equation and isentropic Euler equations of gas dynamics clearly demonstrate the robustness and higher efficiency of the studied high-order schemes. The schemes have also been tested on several challenging examples including discontinuous flux, pressureless gas dynamics, and Burgers equation with non-atomic support. Conservation laws with discontinuous flux arise in numerous applications, including multi-phase flows in heterogeneous porous media and traffic flows on highways with varying surface conditions. The mathematical and numerical analysis of such problems is highly challenging, as many difficulties already appear in the simplest case of 1-D scalar conservation laws with discontinuous flux. As demonstrated in \cite{AMG2005,Mishra2017}, systems of conservation laws with discontinuous flux may admit multiple entropy solutions and exhibit small-scale dependent shock waves. The pressureless gas dynamics system is only weakly hyperbolic, and its solutions develop delta-shocks and vacuum states in finite time, which makes the design of stable and accurate numerical methods particularly challenging; see, e.g.\cite{Sen2025,BBT2006}. Finally, the solutions of the Burgers equation are genuinely measure-valued and cannot be represented by atomic states. Standard numerical methods tend to collapse such distributions into atomic solutions, making it challenging to design schemes that accurately capture the non-atomic structure; see, e.g., \cite{DiPerna1985,FMT2016,Tadmor_1Dburgers}. The obtained  numerical results demonstrate that the designed schemes can be applied to these examples and the resolution of the computed solution improves significantly when higher-order schemes are used.

The paper is organized as follows. In \S\ref{sec2}, we briefly review the existing first-order Lax-Friedrichs Young-measure scheme. In \S \ref{sec3}, we extend the first-order scheme to second order. In \S \ref{sec4}, we present the fifth-order extension in the framework of the A-WENO scheme. Finally, in \S\ref{sec5}, we apply the designed schemes to a number of numerical examples to demonstrate the performance of the studied schemes for both random and deterministic hyperbolic conservation laws.

\section{First-Order Lax--Friedrichs Scheme: A Brief Review}\label{sec2}
In this section, we recall the first-order Lax--Friedrichs Young-measure scheme introduced in \cite{CHLZ25}. Since the high-order schemes developed in Sections~\ref{sec3} and~\ref{sec4} are built upon this first-order method, we first introduce the notation used throughout the paper and then briefly recall the measure-valued closure leading to the linear programming problem.

\subsection{Random Hyperbolic System and Notations}
We consider the one-dimensional random hyperbolic system
\begin{equation}
\bmu_t+\mf(\bmu)_x=\mo,\qquad\bmu(x,\xi,0)=\bm\psi(x,\xi),
\label{1.1a}
\end{equation}
where \(t\) denotes time, \(x\) is the spatial variable, and $\xi\in\Omega$ is the random parameter. The unknown
$$
\bmu=\bmu(x,\xi,t)=\big(u^{(1)},\ldots,u^{(d)}\big)^\top\in\mathbb{R}^d
$$
is the vector of conservative variables, where $d$ is the number of equations in the system. The nonlinear flux is denoted by $\mf:\mathbb{R}^d\to\mathbb{R}^d$. We assume that \eqref{1.1a} is hyperbolic, that is, the Jacobian $A(\bmu):=\frac{\partial \mf}{\partial \bmu}$ has real eigenvalues and is diagonalizable; see, e.g., \cite{Bressan2000,Dafermos2016}. We also assume that \eqref{1.1a} is equipped with a convex entropy $\eta$ and the corresponding entropy flux $q$, satisfying
\begin{equation*}
\partial_t\eta(\bmu(x,\xi,t))+\partial_x q(\bmu(x,\xi,t))\leq 0.
\end{equation*}

We use the following indices. The index $i=1,\ldots,N_\xi$ refers to the discretization in the random variable $\xi$, the index $j=1,\ldots,N_x$ refers to the spatial cell, and $n=0,1,2,\ldots$ denotes the time level $t^n$. The index $\ell=1,\ldots,N_{\bmu}$ is used for the discretization of the phase space, namely, the space of admissible values of the conservative vector $\bmu$. We use $\bm z=(z^{(1)},\ldots,z^{(d)})^\top\in\mathbb{R}^d$ to denote a generic phase-space variable, namely, a possible value of the conservative vector $\bmu$. The corresponding phase-space grid points are denoted by $\bm z_\ell\in\mathbb{R}^d$, $\ell=1,\ldots,N_{\bmu}$. For a system, $\bm z_\ell$ is a vector-valued phase-space point. If the phase space is discretized by a tensor-product grid, then $N_{\bmu}$ denotes the total number of phase-space grid points. We denote by $\Delta V$ the corresponding phase-space cell volume. More precisely, if the grid is uniform in each phase-space direction, then $\Delta V=\prod_{r=1}^d \Delta z^{(r)}$. 
In the scalar case, $d=1$, and thus $\Delta V=\Delta u$. In this case, we also write $N_u$ instead of $N_{\bmu}$.

We assume that $\xi\sim \mathcal{U}(\Omega)$, namely, $\xi$ is uniformly distributed on $\Omega$. Hence, $p(\xi)=p_0:=\frac{1}{|\Omega|}$. Let $\{\xi_i\}_{i=1}^{N_\xi}$ be equidistant points in $\Omega$ with mesh size $\Delta\xi$. Following \cite{CHLZ25}, we use the piecewise constant orthonormal
basis functions
\begin{equation*}
\phi_i(\xi)=\frac{1}{\sqrt{p_0\Delta\xi}}\,\chi_{[\xi_i-\Delta\xi/2,\xi_i+\Delta\xi/2]}(\xi),\quad i=1,\ldots,N_\xi.
\end{equation*}
This choice mimics the stochastic collocation approach in the random variable.

We denote by $\bmu_{i,j}^n\in\mathbb{R}^d$ the discrete moment of the corresponding discrete Young measure, at $(\xi_i,x_j,t^n)$. For a fixed spatial cell $j$ and time level $n$, we collect all random-space moments into the vector $\boldsymbol{\bmu}_j^n:=\big(\bmu_{1,j}^n,\ldots,\bmu_{N_\xi,j}^n\big)^\top$. The dependence of the discrete Young measure on $\boldsymbol{\bmu}_j^n$ will be emphasized below.

\subsection{Moment Formulation and Young-Measure Closure}
We now recall the idea behind the linear programming formulation used in \cite{CHLZ25}. Testing \eqref{1.1a} against the basis function $\phi_i(\xi)$ and integrating with respect to $\xi$ formally gives the moment equation
\begin{equation*}
\partial_t\bm M_i(t,x)+\partial_x\bm G_i(t,x)=\mo,\qquad i=1,\ldots,N_\xi,
\end{equation*}
where
$$
\bm M_i(t,x):=\int_\Omega \phi_i(\xi)\bmu(t,x,\xi)p(\xi)\,d\xi \quad {\rm and }\quad \bm G_i(t,x):=\int_\Omega \phi_i(\xi)\mf(\bmu(t,x,\xi))p(\xi)\,d\xi.
$$

To obtain a closure, the method in \cite{CHLZ25} uses a measure-valued formulation. Let $\nu_{t,x,\xi}\in\mathcal P(\mathbb{R}^d)$ be a Young measure parameterized by
$(t,x,\xi)$. The first moment and the corresponding averaged flux with respect to this Young measure are $\int_{\mathbb{R}^d}\bm z\,d\nu_{t,x,\xi}(\bm z)$ and 
$\int_{\mathbb{R}^d}\mf(\bm z)\,d\nu_{t,x,\xi}(\bm z)$, respectively.  If $\nu_{t,x,\xi}$ is a Dirac measure concentrated at the classical state $\bmu(t,x,\xi)$, namely,
$\nu_{t,x,\xi}=\delta_{\bmu(t,x,\xi)}$, then the usual weak formulation of \eqref{1.1a} is recovered.

For given moments, one seeks an admissible family of probability measures $\{\nu_{\xi;t,x}\}_{\xi\in\Omega}$ by minimizing the entropy subject to the moment constraints. More precisely, for a fixed \((t,x)\), the closure measure is obtained from
\begin{align}
\nu_{\xi;t,x}^*&=\mbox{\rm argmin}_{\nu_{\xi;t,x}}\int_\Omega\int_{\mathbb{R}^d}\eta(\bm z)\,d\nu_{\xi;t,x}(\bm z)\,p_0\,d\xi,\label{eq:2.continuousLP-a}\\
\mbox{\rm subject to}\qquad& \int_\Omega\int_{\mathbb{R}^d}\phi_i(\xi)\bm z\,d\nu_{\xi;t,x}(\bm z)\,p_0\,d\xi=\bm M_i(t,x),\qquad i=1,\ldots,N_\xi .
\label{eq:2.continuousLP-b}
\end{align}
The corresponding closed flux is then given by
\begin{equation}
\bm G_i(t,x)\approx\int_\Omega\int_{\mathbb{R}^d}\phi_i(\xi)\mf(\bm z)\,d\nu_{\xi;t,x}^*(\bm z)\,p_0\,d\xi .
\label{eq:2.closedflux}
\end{equation}
This is the basic motivation for the linear programming problem used in the numerical scheme. Indeed, both the objective functional in \eqref{eq:2.continuousLP-a} and the constraints in \eqref{eq:2.continuousLP-b} are linear with respect to the measure. Therefore, after discretizing the phase space, the closure problem becomes a finite-dimensional linear programming problem.

\subsection{Fully Discrete Linear-Programming Problem}

We now describe the fully discrete closure used in this paper. Let $\{\bm z_\ell\}_{\ell=1}^{N_{\bmu}}$ be a uniform discretization of the phase space.
For given moments \(\bmu_{i,j}^n\) at each fixed spatial cell \(j\) and time level \(n\), we compute a discrete approximation
$$
\mu_{i,j,\ell}^{*,n}\approx\nu_{\xi_i;t^n,x_j}^*(\bm z_\ell),\qquad i=1,\ldots,N_\xi,\quad \ell=1,\ldots,N_u,
$$
by solving the following linear programming problem:
\begin{align}
\{\mu_{i,j,\ell}^{*,n}\}&=\mbox{\rm argmin}_{\{\mu_{i,\ell}\}}\Delta\xi\,\Delta V\sum_{i=1}^{N_\xi}\sum_{\ell=1}^{N_{\bmu}}\eta(\bm z_\ell)p_0\mu_{i,\ell},\label{linprog-a}\\
  \mbox{\rm subject to}\qquad &\mu_{i,\ell}\geq 0,\quad \forall (i,\ell),\label{linprog-b}\\
  &\mu_{i,\ell}\leq \lambda_F /\Delta V,\quad \forall (i,\ell),\label{linprog-c}\\
  &\Delta V\sum_{\ell=1}^{N_{\bmu}}\mu_{i,\ell}=1,\quad \forall i,\label{linprog-d}\\
  &\Delta V\sum_{\ell=1}^{N_{\bmu}}\bm z_\ell\mu_{i,\ell}=\bmu_{i,j}^n,\forall i.\label{linprog-e}
\end{align}

The constraints \eqref{linprog-b}--\eqref{linprog-d} ensure that, for each random index $i$, the coefficients $\{\mu_{i,\ell}\}_{\ell=1}^{N_{\bmu}}$ define a discrete probability measure. The constraint \eqref{linprog-e} enforces that the first moment of this measure is equal to the given vector moment $\bmu_{i,j}^n$. Thus,
\eqref{linprog-a}--\eqref{linprog-e} selects, among all admissible discrete probability measures with the same first moment, one that minimizes the entropy.

\begin{rmk}
Here, $0<\lambda_F\leq 1$ controls the admissible support of the discrete Young measure. The case $\lambda_F=1$ allows atomic solutions, while choosing $\lambda_F<1$ prevents the entire mass from concentrating at a single phase-space grid point and thus permits non-atomic discrete Young measures.
\end{rmk}

\begin{rmk}
We note that the linear programming problem \eqref{linprog-a}--\eqref{linprog-e} is written in a coupled form for all random indices $i$. However, since the constraints are independent for different $i$, it can be solved as $N_\xi$ smaller linear programming problems. This observation is useful in the implementation and reduces the computational cost. Notice that the accuracy of the approximated values of $\mu^*$ depends on the tolerance  of the black–box solver $\texttt{linprog}$. In all of the examples reported in \S \ref{sec5}, we use the default settings of this function. In order to obtain more accurately approximated values of $\mu^*$, we also propose a new approach for the case $\lambda_F=1$ to obtain the values of $\mu^*$ by using the KKT-based solver; see Appendix \ref{appd}. 
\end{rmk}

\subsection{First-order Lax--Friedrichs Scheme}
After obtaining $\mu_{i,j,\ell}^{*,n}$, we define the closed numerical flux
\begin{equation}\label{2.4}
\bm F_{i,j}^n=\Delta V \sum_{\ell=1}^{N_{\bmu}}\mf(\bm z_\ell)\mu_{i,j,\ell}^{*,n}.
\end{equation}
This is the discrete counterpart of the measure-valued flux in \eqref{eq:2.closedflux}. The first-order Lax--Friedrichs scheme for the discrete moments is then given by
\begin{equation*}
\bmu_{i,j}^{n+1} =\frac{1}{2}\left(\bmu_{i,j+1}^{n}+\bmu_{i,j-1}^{n}\right)-\frac{\Delta t}{2\Delta x}\left(\bm F_{i,j+1}^{n}-\bm F_{i,j-1}^{n}\right),
\end{equation*}
where $\Delta x$ is the spatial mesh size and $\Delta t$ is the time step.

The time step is determined by the CFL condition
\begin{equation}
\max_{1\leq i\leq N_\xi,\;1\leq j\leq N_x}\left\{\Delta V\sum_{\ell=1}^{N_{\bmu}}\sigma(\bm z_\ell)\mu_{i,j,\ell}^{*,n}\right\}\Delta t=\mathrm{CFL}\Delta x.
\label{eq:2.amax}
\end{equation}
Here, $\sigma(\bm z_\ell)$ denotes the spectral radius of the Jacobian $\partial\mf/\partial\bmu$ evaluated at the phase-space point $\bm z_\ell$.

The first-order scheme \eqref{linprog-a}--\eqref{eq:2.amax} is the starting point for the high-order extensions developed below. In Sections~\ref{sec3} and \ref{sec4}, we replace the first-order Lax--Friedrichs spatial discretization by second- and fifth-order local Lax--Friedrichs discretizations, respectively, while the Young-measure closure through the linear programming problem remains the central reconstruction step.

\section{Second-Order Local Lax–Friedrichs (LLF) Scheme}\label{sec3}
We now extend the first-order Lax-Friedrichs scheme introduced in \S \ref{sec2} to the second-order of accuracy. The moments $\bmu^n_{i,j}$ (in the rest of
the paper, we will suppress the time-dependence of all of the indexed quantities for the sake of brevity) of the 1-D random hyperbolic conservation laws \eref{1.1a} are evolved in time by solving the following system of ODEs
\begin{equation}\label{3.1}
\frac{\mathrm{d} \bmu_{i,j}}{\mathrm{d}t} = -\frac{\bm{{\cal F}}_{i,\jph} - \bm{{\cal{F}} }_{i,\jmh}}{\dx}.
\end{equation}
Here,  $\bm{{\cal F}}_{i,\jph}=\bm{{\cal F}}_{i,\jph}\big(\bmu^{-,*}_{i,\jph}, \bmu^{+,*}_{i,\jph}\big)$ stand for the FV LLF fluxes (see \cite{KTcl,Rus61}):
\begin{equation}\label{3.2}
  \bm{{\cal{F}}}_{i,\jph} = \frac{ \mf(\bmu^{-,*}_{i,\jph})+\mf(\bmu^{+,*}_{i,\jph})}{2}
                         -\frac{a_{i,\jph}}{2}\big(\bmu^{+,*}_{i,\jph}-\bmu^{-,*}_{i,\jph} \big),
\end{equation}
where $\bmu^{\pm,*}_{i,\jph}$ are the one-sided point values of  $\bmu^*_{i,j}=\Delta V\sum\limits_{\ell=1}^{N_{\bmu}}\bm z_\ell  \, \mu^*_{i,j,\ell}$ at the cell interface $(\xi, x)=(\xi_i, x_\jph)$  estimated using a piecewise linear interpolant:
\begin{equation}\label{3.3}
\widetilde \bmu^*(x)=\,\bmu^*_{i,j}+(\bmu^*_x)_{i,j}(x-x_j),
\end{equation}
which leads to
\begin{equation}\label{3.4}
\bmu^{-,*}_{i,\jph}=\,\bmu^*_{i,j}+\frac{\dx}{2}(\bmu^*_x)_{i,j},\quad \bmu^{+,*}_{i,\jmh}=\,\bmu^*_{i,j}-\frac{\dx}{2}(\bmu^*_x)_{i, j}.
\end{equation}
Here, the values of $\mu^*_{i,j,\ell}:=\mu^{*,n}_{i,j,\ell}$  are obtained by solving the linear programming \eref{linprog-a}--\eref{linprog-e}.

In order to ensure the non-oscillatory nature of the reconstruction \eref{3.4}--\eref{3.5}, the slopes $(\bmu^*_x)_{i,j}$ in \eref{3.3} need to be computed with the help of a nonlinear limiter.  In all of the numerical experiments reported in \S\ref{sec5}, we have used a generalized minmod limiter
\cite{lie03,Nessyahu90,Sweby84}:
\begin{equation}
(\bmu^*_x)_{i,j}={\rm minmod}\left(\theta\frac{\,\bmu^*_{i,j}-\,\bmu^*_{i,j-1}}{\dx},\,\frac{\,\bmu^*_{i,j+1}-\,\bmu^*_{i,j-1}}{2\dx},\,
\theta\frac{\,\bmu^*_{i,j+1}-\,\bmu^*_{i,j}}{\dx}\right),\quad\theta\in[1,2],
\label{3.5}
\end{equation}
applied in a component-wise manner. Here, the minmod function is defined as
\begin{equation}\label{3.6}
{\rm minmod}(z_1,z_2,\ldots):=\begin{cases}
\min_j\{z_j\}&\mbox{if}~z_j>0\quad\forall\,j,\\
\max_j\{z_j\}&\mbox{if}~z_j<0\quad\forall\,j,\\
0            &\mbox{otherwise}.
\end{cases}
\end{equation}
The parameter $\theta$ in \eref{3.5} is used to control the amount of numerical viscosity present in the resulting scheme. In general, larger values of $\theta$ correspond to sharper but, in general, more oscillatory reconstructions. In all of the numerical examples reported in \S \ref{sec5}, we take $\theta=1.5$.

Finally, the local speed of propagation $a_{i,\jph}$ is estimated using the eigenvalues $\lambda_1(A)\le\ldots\le\lambda_d(A)$ of the Jacobian
$A=\frac{\partial \mf}{\partial \bmu}$, defined by 
\begin{equation*}
a_{i,\jph}=\max\Big\{\big|\lambda_1\big(\bmu^{-,*}_{i,\jph})\big)\big|,\,\big|\lambda_1\big(\bmu^{+,*}_{i,\jph})\big)\big|,\,
\big|\lambda_d\big(\bmu^{-,*}_{i,\jph})\big)\big|,\,\big|\lambda_d\big(\bmu^{+,*}_{i,\jph})\big)\big|\Big\}.
\end{equation*}

\begin{rmk}
The LLF fluxes given by \eref{3.2} can also be replaced by 
\begin{equation}\label{3.7}
  \bm{{\cal{F}}}_{i,\jph} = \frac{\Delta V \sum_{\ell=1}^{N_{\bmu}} \mf(\bm z_\ell)\mu^{-,*}_{i,j,\ell}+\Delta V \sum_{\ell=1}^{N_{\bmu}}\mf(\bm z_\ell)\mu^{+,*}_{i,j,\ell}}{2}
                         -\frac{a_{i,\jph}}{2}\big(\bmu^{+,*}_{i,\jph}-\bmu^{-,*}_{i,\jph} \big),
\end{equation}
where
\begin{equation}\label{3.8}
  \bmu^{-,*}_{i,\jph}=\Delta V \sum_{\ell=1}^{N_{\bmu}}\bm z_\ell  \mu^{-,*}_{i,j,\ell}, \quad {\rm and} \quad \bmu^{+,*}_{i,\jph}=\Delta V \sum_{\ell=1}^{N_{\bmu}}\bm z_\ell \mu^{+,*}_{i,j,\ell},
\end{equation}
and $\mu^{\pm,*}_{i,j,\ell}$ are obtained by solving the linear programming \eref{linprog-a}--\eref{linprog-e} but with $\bmu_{i,j}$ in \eref{linprog-e}  replaced by $\bmu^{-}_{i,\jph}$ and $\bmu^{+}_{i,\jph}$, respectively. The values of $\bmu^{\pm}_{i,\jph}$ are the one-sided point values of $\bmu_{i,j}$ at the cell interface $(\xi, x)=(\xi_i, x_\jph)$  estimated using the piecewise linear interpolant \eref{3.3}--\eref{3.6}.

It is noted that the numerical results computed by the LLF fluxes \eref{3.2} and \eref{3.7} coincide, which will be demonstrated in one example in \S \ref{sec5}. However, one can clearly see that the computations of LLF fluxes \eref{3.2} are computationally much cheaper than \eref{3.7} since it only needs to solve the linear programming problem once, which can save around half of the CPU times. 
\end{rmk}

\begin{rmk}
Notice that applying the minmod limiter to the piecewise linear interpolation \eref{3.3}--\eref{3.6} with large values of $\theta$ directly may lead to relatively large oscillations in the computed solution. In order to suppress these spurious oscillations, we apply the minmod limiter to the local characteristic variables, which are obtained using the LCD based linear interpolation (see, e.g., \cite{Jiang13,Qiu02,Shu20,CKM2025,CHK25} and references therein); see Appendix A for details.
\end{rmk}

\section{Fifth-Order Local Lax–Friedrichs (LLF) Scheme}\label{sec4}
In this section, we extend the second-order scheme introduced in \S \ref{sec3} to the fifth order of accuracy in the framework of the FD A-WENO scheme introduced in \cite{Jiang13} (see also \cite{Liu17}), which has been proven to be a powerful tool for generalizing low-order FV schemes to higher-order FD ones.

According to \cite{Jiang13}, the point values $\bmu_{i,j}$ are evolved in time by solving the following system of ODEs:
\begin{equation}
\frac{{\rm d}\bmu_{i,j}}{{\rm d}t}=-\frac{{\bm{{\cal H}}_{i,\jph}}-{\bm{{\cal H}}_{i,\jmh}}}{\dx},
\label{4.1}
\end{equation}
where $\bm{{\cal H}}_{i,\jph}$ is the (fifth-order accurate) numerical flux defined by
\begin{equation}\label{4.2}
\bm{{\cal H}}_{i,\jph}=\bm{{\cal F}}^{\rm }_{i,\jph}-\frac{1}{24}(\dx)^2(\mF_{xx})_{i,\jph}+\frac{7}{5760}(\dx)^4(\mF_{xxxx})_{i,\jph}.
\end{equation}
Here, $\bm{{\cal F}}^{\rm }_{i,\jph}$ is the FV numerical flux as in \eref{3.2},  $(\mF_{xx})_{i,\jph}$ and $(\mF_{xxxx})_{i,\jph}$ are the higher-order correction term computed by the fourth- and second-order accurate FDs, respectively; see, e.g., \cite{CKX23,CHT25,Gao20,CKX22}. In this paper, we have used the following new more efficient higher-order correction terms from \cite{CKX23}:
\begin{equation}\label{4.3}
\begin{aligned}
&(\mF_{xx})_{i,\jph} =\frac{1}{12(\dx)^2}\Big[-\bm{{\cal F}}_{i,j-\frac{3}{2}}+16\bm{{\cal F}}_{i,\jmh}-30\bm{{\cal F}}_{i,\jph}
+16\bm{{\cal F}}_{i,j+\frac{3}{2}}-\bm{{\cal F}}_{i,j+\frac{5}{2}}\Big],\\
&(\mF_{xxxx})_{i,\jph}=\frac{1}{(\dx)^4}\Big[\bm{{\cal F}}_{i,j-\frac{3}{2}}-4\bm{{\cal F}}_{i,\jmh}+6\bm{{\cal F}}_{i,\jph}-
4\bm{{\cal F}}_{i,j+\frac{3}{2}}+\bm{{\cal F}}_{i,j+\frac{5}{2}}\Big],
\end{aligned}
\end{equation}
which are computed using the FV numerical fluxes $\bm{{\cal F}}^{\rm }_{i,\jph}$, and have been proved to be more efficient than the old version (see \eref{4.4} below), while affecting neither the accuracy nor the quality of resolution. 

\begin{rmk}
As in the second-order scheme in \S \ref{sec3}, one can also use the FV numerical fluxes $\bm{{\cal F}}^{\rm }_{i,\jph}$ defined by \eref{3.7}--\eref{3.8}, but it is more computationally expensive.  At the same time, to ensure that the obtained numerical scheme \eref{4.1}--\eref{4.3} is of the fifth order of accuracy, the one-sided point values $\bmu^{-,*}_{i,\jph}$ and $\bmu^{+,*}_{i,\jph}$, or $\bmu^{-}_{i,\jph}$ and $\bmu^{+}_{i,\jph}$, employed to compute the numerical flux $\bm{{\cal F}}^{\rm }_{i,\jph}$ need to be at least fifth order accurate. This can be done by using a certain nonlinear limiting procedure like the fifth-order WENO-Z interpolation from \cite{Gao20,Jiang13,Liu17} applied to the local characteristic variables; see Appendix \ref{appb} for details.
\end{rmk}

\begin{rmk} 
The high-order correction terms $(\mF_{xx})_{i,\jph}$ and $(\mF_{xxxx})_{i,\jph}$  in \eref{4.2} can also be computed by 
\begin{equation}\label{4.4}
\begin{aligned}
&(\mF_{xx})_{i,j+\frac{1}{2}}=\frac{1}{48 (\Delta x)^2}(-5\mF_{i,j-2}+39\mF_{i,j-1}-34 \mF_{i,j}-34\mF_{i,j+1}+39\mF_{i,j+2}-5\mF_{i,j+3}) \\
&(\mF_{xxxx})_{i,j+\frac{1}{2}}=\frac{1}{2(\Delta x)^4}(\mF_{i,j-2}-3\mF_{i,j-1}+2\mF_{i,j}+2\mF_{i,j+1}-3\mF_{i,j+2}+\mF_{i,j+3}),
\end{aligned}
\end{equation}
where ${\mF}_{i,j}$ are defined by \eref{2.4}.

Notice that the main difference between \eref{4.3} and \eref{4.4} is that in \eref{4.3}, we compute the high-order correction terms $(\mF_{xx})_{i,\jph}$ and
$({\mF_{xxxx}})_{i,\jph}$ terms with the help of the FV numerical fluxes $\bm{{\cal F}}^{\rm }_{i,\jph}$ instead of the point values $\mF_{i,j}$. One can clearly see that  \eref{4.3} is computationally much cheaper than \eref{4.4} as there is no need to compute the point values $\mF_{i,j}$, where the term $\mu_{i,\ell}^*(\bm{u}_j^n )$ is obtained by solving the linear programming problem \eref{linprog-a}--\eref{linprog-e}.
\end{rmk}

\begin{rmk}
In this paper, we extend the second order scheme to the fifth order of accuracy in the framework of FD A-WENO methods, but it is also easy to extend it to even higher-order accuracy in this framework. 
\end{rmk}

\section{Numerical Examples} \label{sec5}
In this section, we apply the studied first-, second-, and fifth-order schemes to several numerical examples to demonstrate the superiority and robustness of the studied high-order schemes. For the sake of brevity, the tested first-, second-, and fifth-order schemes will be referred to as the 1-, 2-, and 5-Order schemes, respectively. As in \cite{CHLZ25}, we present the results obtained using both the Young-measure and stochastic collocation (which will be referred to as collocation) methods, as a comparison. We numerically integrate semi-discrete ODE systems \eref{3.1} and \eref{4.1} using the SSP RK3 method and use the CFL number $0.45$ in all examples. We take $\lambda_F=1$ in Examples 1--6 and 0.05 in Example 7.

\subsection{Numerical Examples for Random Hyperbolic Conservation Laws}
We first apply the studied Young-measure and collocation schemes to a number of numerical examples for random hyperbolic conservation laws, including the Burgers equation and isentropic Euler system to show the validity of the Young-measure formulation.

\subsubsection{Burgers Equation}
In this section, we consider the 1-D Burgers equation
\begin{equation}\label{5.1}
u_t+ \Big( \frac{u^2}{2} \Big)_x =0.
\end{equation}
The corresponding entropy function $\eta$ in \eref{linprog-a} is defined by $\eta=\frac{1}{2} u^2$.

\begin{figure}[ht!]
\centerline{\includegraphics[trim=0.1cm 0.2cm 1.cm 0.8cm, clip, width=5.cm]{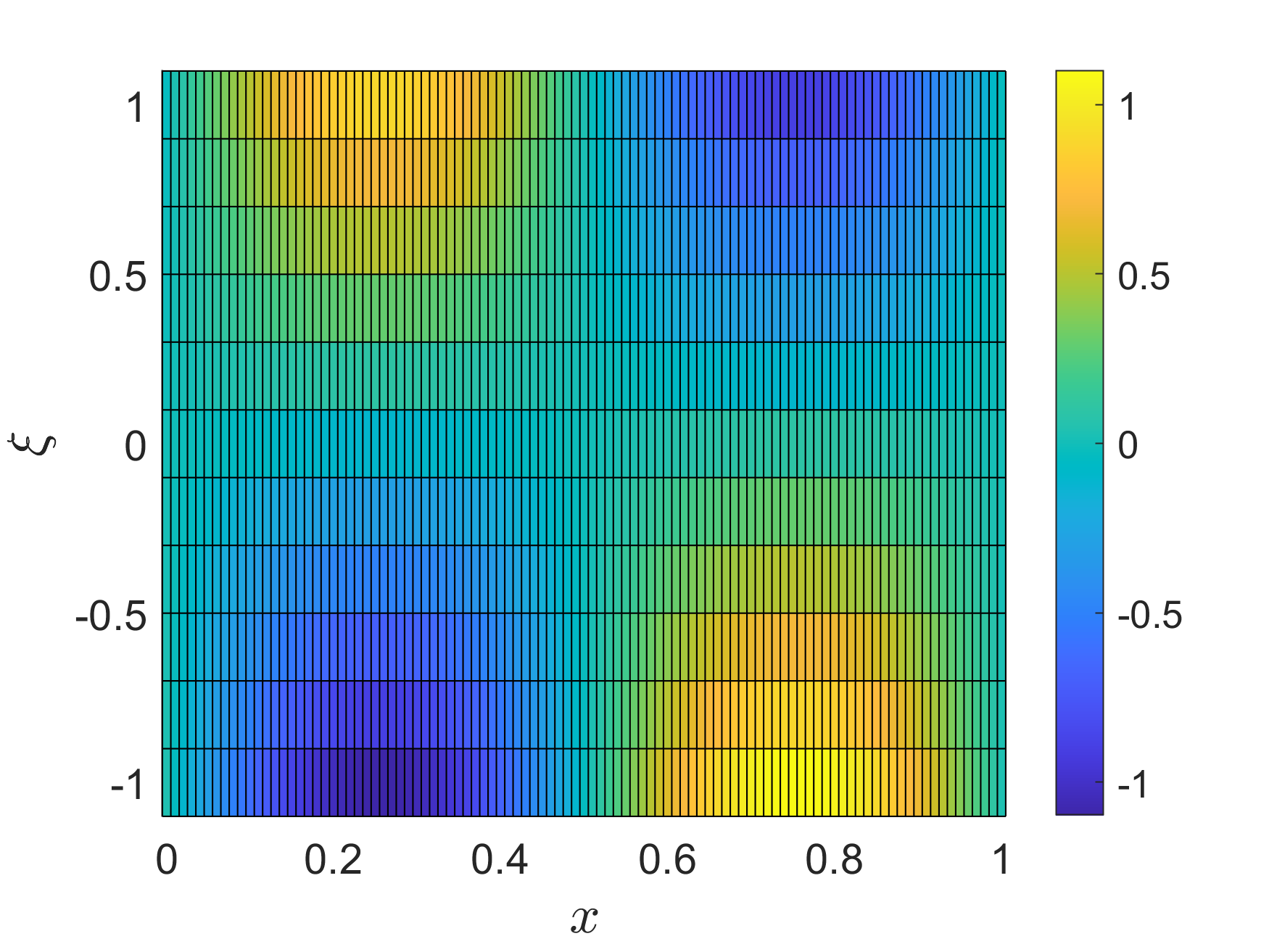}}
\caption{\sf Example 1: Initial data.\label{fig1.1}}
\end{figure}

\paragraph*{Example 1.}
In the first example taken from \cite{CHLZ25}, we consider the Burgers equation \eref{5.1} with the following initial data 
\begin{equation*}
u(x,\xi,0) = \xi \sin\left( 2 \pi x\right),
\end{equation*}
prescribed in the computational domain $[0,1]\times[-1,1]$ subject to the free boundary conditions. We plot the initial data in Figure \ref{fig1.1} and note that even though the initial data are smooth, a shock will be formed at later time $t>0$ for any  $\xi>0$. For $\xi<0$, the solution remains smooth.

\begin{figure}[ht!]
\centerline{\includegraphics[trim=0.1cm 0.2cm 0.9cm 0.8cm, clip, width=5.cm]{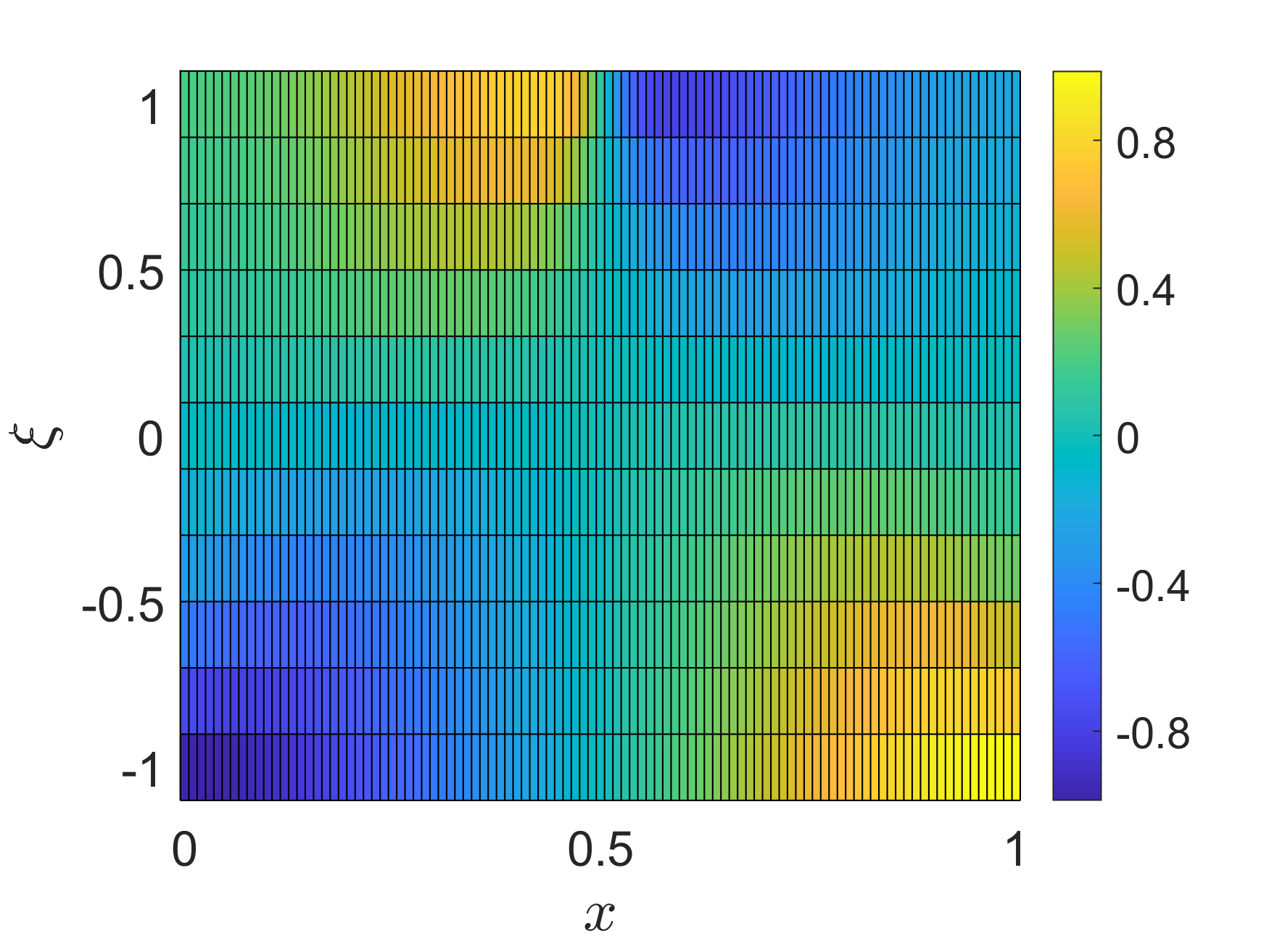} \hspace{0.5cm}
            \includegraphics[trim=0.1cm 0.2cm 0.9cm 0.8cm, clip, width=5.cm]{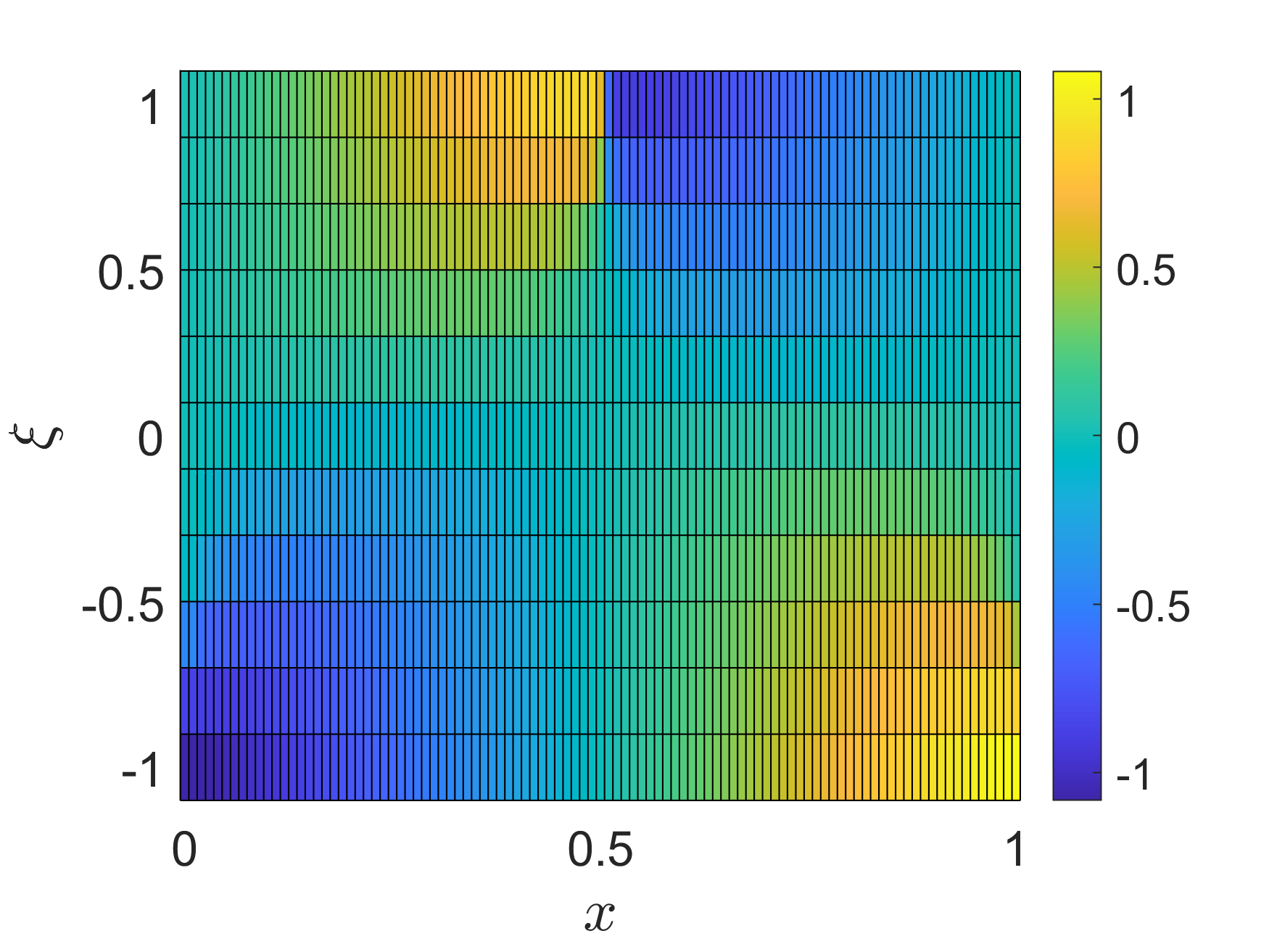} \hspace{0.5cm}
            \includegraphics[trim=0.1cm 0.2cm 0.9cm 0.8cm, clip, width=5.cm]{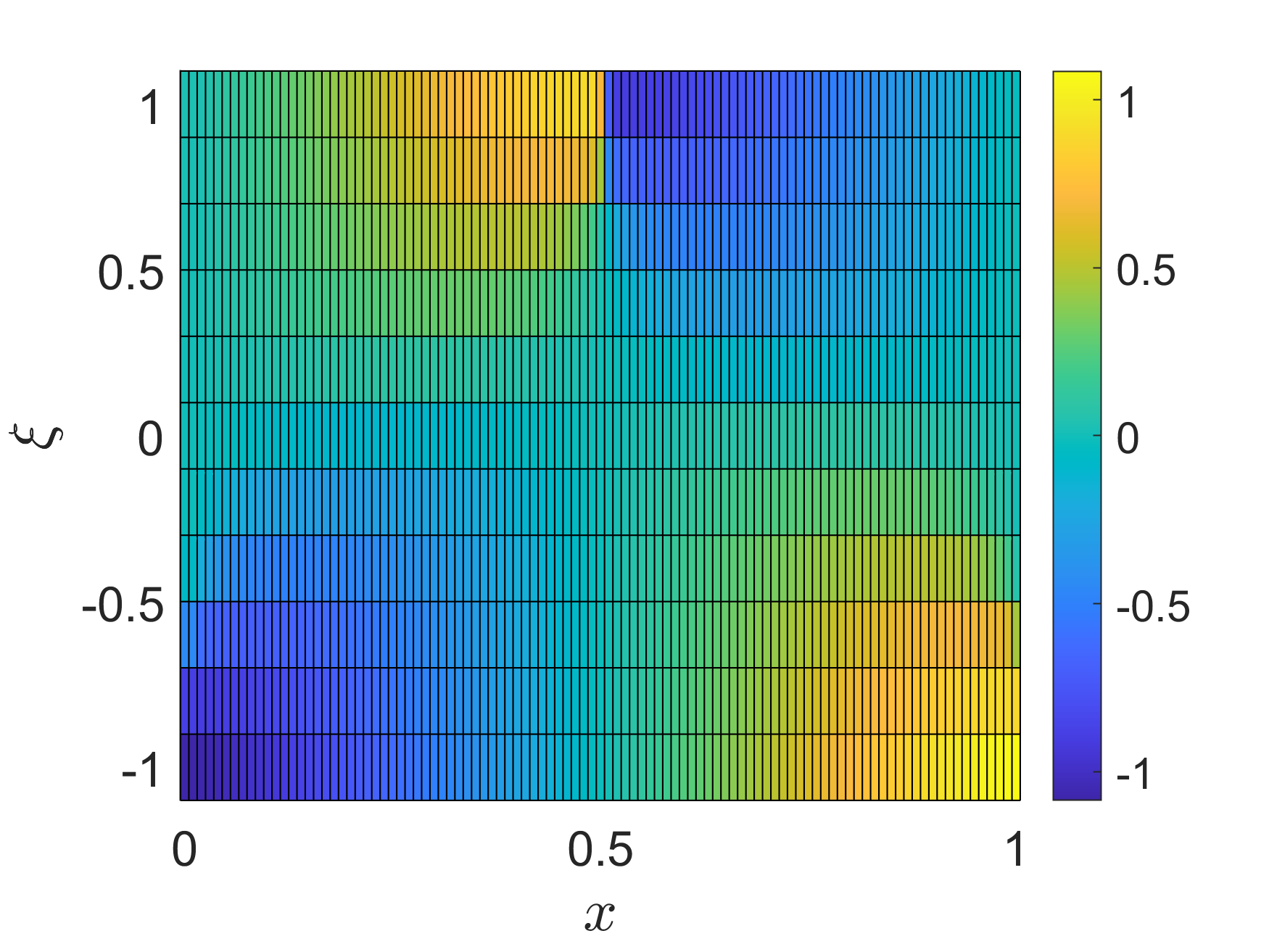}}
\vskip 10pt 
\centerline{\includegraphics[trim=0.1cm 0.2cm 0.9cm 0.8cm, clip, width=5.cm]{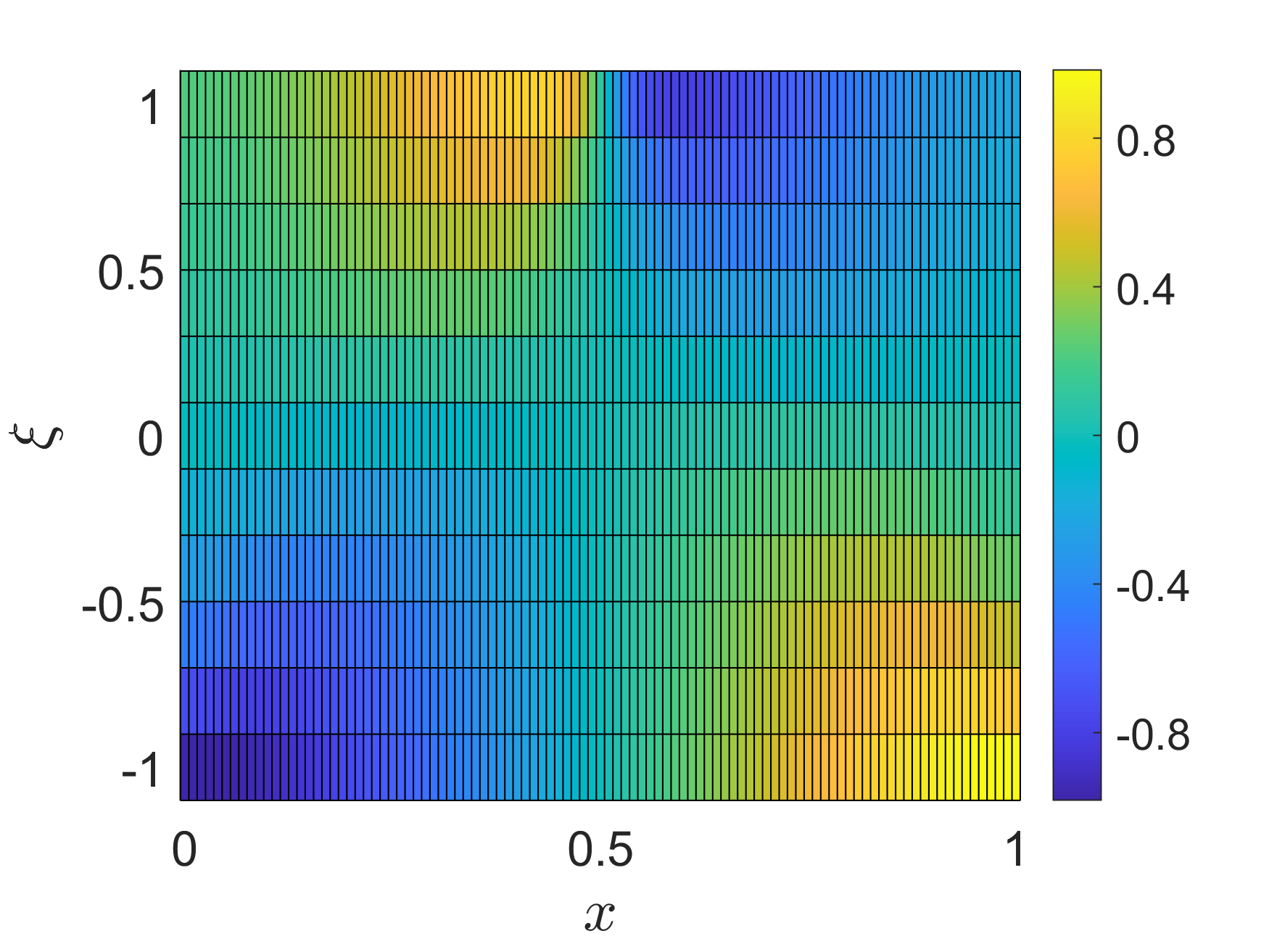}\hspace{0.5cm}
            \includegraphics[trim=0.1cm 0.2cm 0.9cm 0.8cm, clip, width=5.cm]{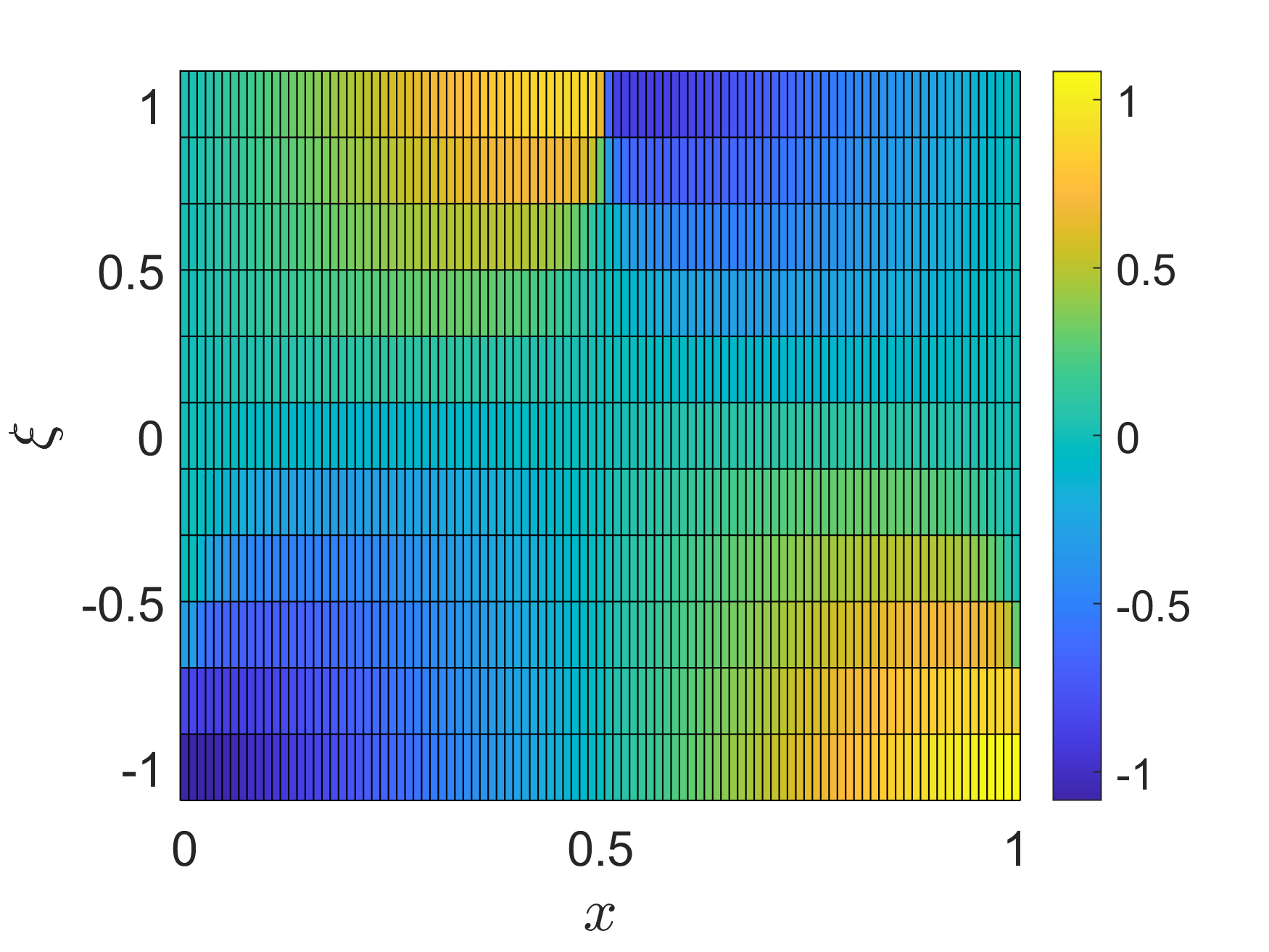} \hspace{0.5cm}
            \includegraphics[trim=0.1cm 0.2cm 0.9cm 0.8cm, clip, width=5.cm]{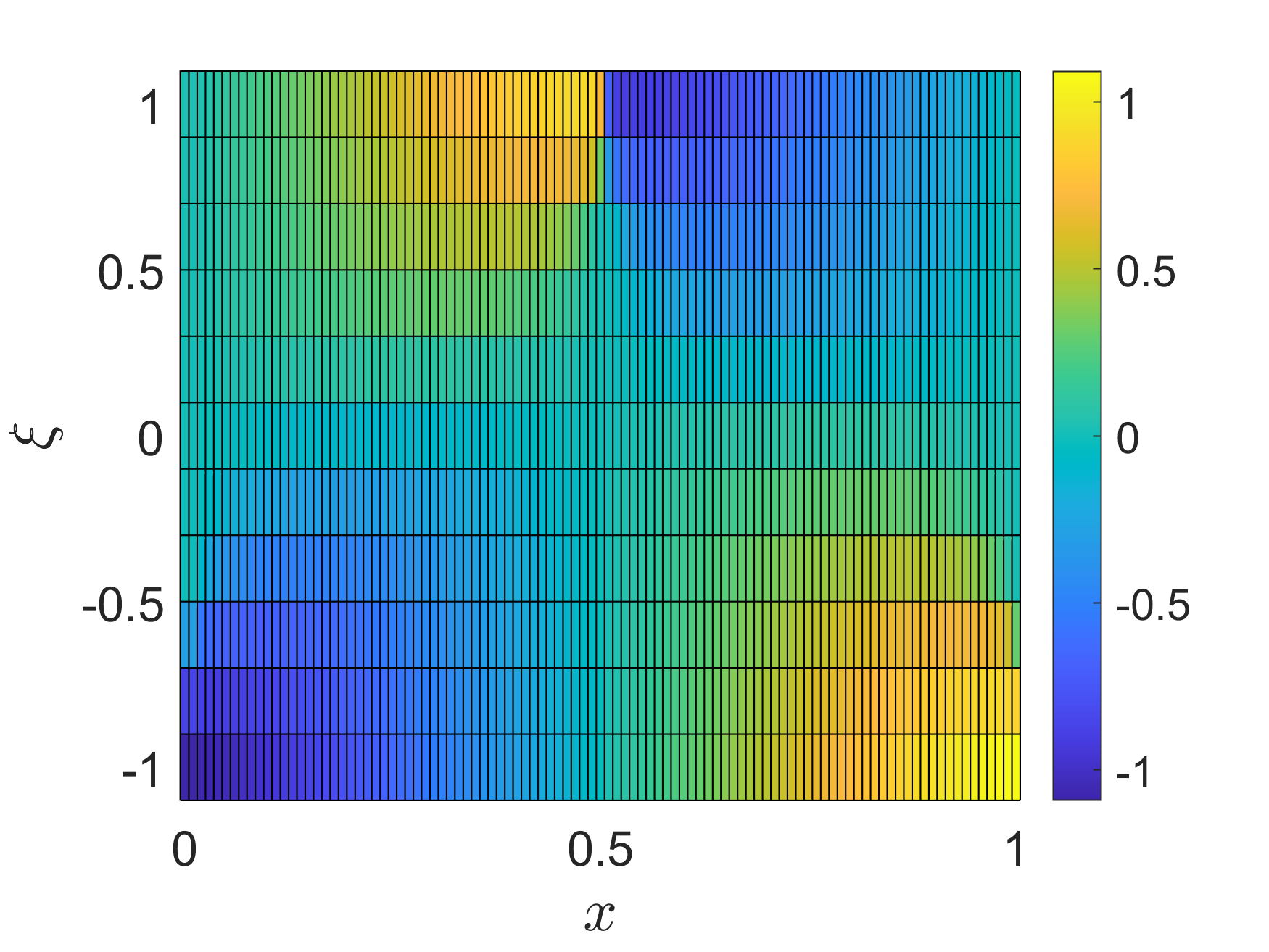}}
\caption{\sf Example 1: Numerical results computed by the 1-Order (left column), 2-Order (middle column), and 5-Order (right column) collocation (top row) and Young-measure (bottom row) schemes.}
\label{fig1.1a}
\end{figure}

We compute the numerical solutions until the final time $t=0.25$ by the studied 1-Order, 2-Order, and 5-Order Young-measure schemes on the uniform mesh with $N_x=100$ and $N_\xi=10$. For the discretization in phase space, we set $u \in [-1.5,1.5]$ and $N_{u} = 100$. The numerical results obtained are presented in Figure \ref{fig1.1a} together with the numerical results computed by the corresponding collocation approaches, where we present the results from a top-view perspective to enhance visualization of the resolution near the shock wave. One can clearly see that the resolution of the computed results near the shock waves improves significantly with the use of high-order schemes, especially when transiting from the 1-Order scheme to the 2-Order one.

\begin{figure}[ht!]
\centerline{\includegraphics[trim=0.9cm 0.4cm 1.3cm 0.2cm, clip, width=5.cm]{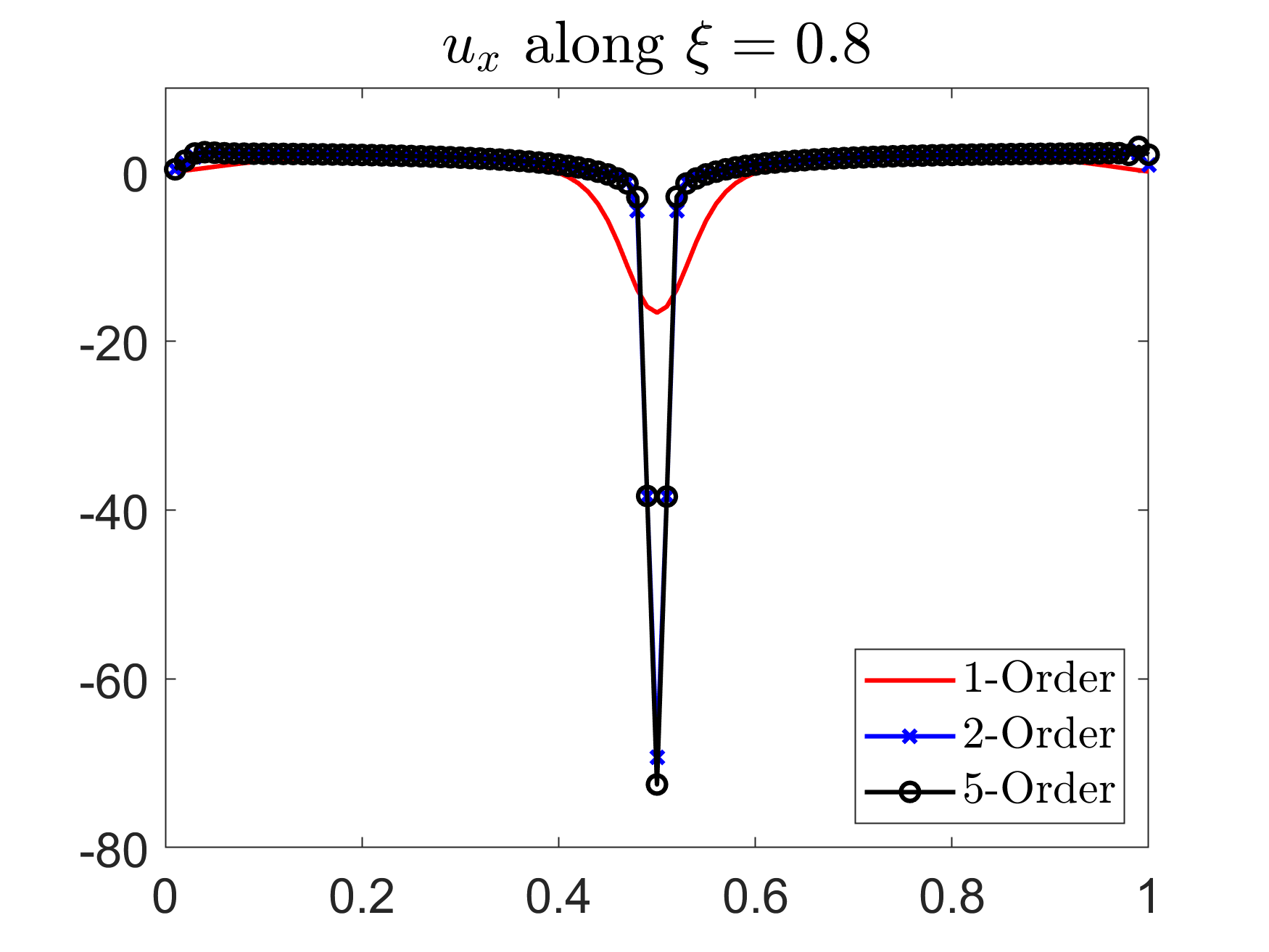}}
\caption{\sf Example 1: $u_x$ along $\xi=0.8$ computed by the 1-Order, 2-Order, and 5-Order Young-measure schemes.}
\label{fig1.1b}
\end{figure}

In order to show the differences between the results computed by the studied schemes, we compute $u_x$ along the line $\xi=0.8$ using the central difference scheme  and plot the obtained results in Figure \ref{fig1.1b}. One can clearly see that the 2-Order and 5-Order schemes achieve much higher resolution than the 1-Order one, while the 5-Order scheme is slightly better than the 2-Order one. 

We also estimate the $L^1$-, $L^2$-, and $L^\infty$-errors between the numerical results computed by the Young-measure methods with the reference solution, which is computed by the 1-Order collocation method on a much finer mesh with $N_x=4000$ and $N_\xi=400$. The obtained numerical results are reported in Table \ref{tab1a}, which once again show that the 2-Order and 5-Order schemes are much more accurate than the 1-Order scheme.

\begin{table}[ht!]
\centering
\begin{tabular}{cccccccccc}
\toprule 
 \multicolumn{3}{c}{$L^1$-error} & \multicolumn{3}{c}{$L^2$-error}& \multicolumn{3}{c}{$L^\infty$-error}\\
\cline{1-9}
 1-Order & 2-Order & 5-Order & 1-Order & 2-Order & 5-Order & 1-Order & 2-Order & 5-Order\\
\hline
7.90e-02 & 2.62e-02 & 2.53e-02 & 1.10e-01 & 2.78e-02 & 2.74e-02 & 6.99e-01 & 1.95e-01 & 1.91e-01\\
\bottomrule
\end{tabular}
\caption{\sf Example 1: The $L^1$-, $L^2$-, and $L^\infty$-errors of $u$ between the numerical results computed by the 1-Order, 2-Order, and 5-Order Young-measure schemes and the reference solution.\label{tab1a}}
\end{table}

\begin{figure}[ht!]
\centerline{\includegraphics[trim=0.1cm 0.2cm 0.9cm 0.8cm, clip, width=5.cm]{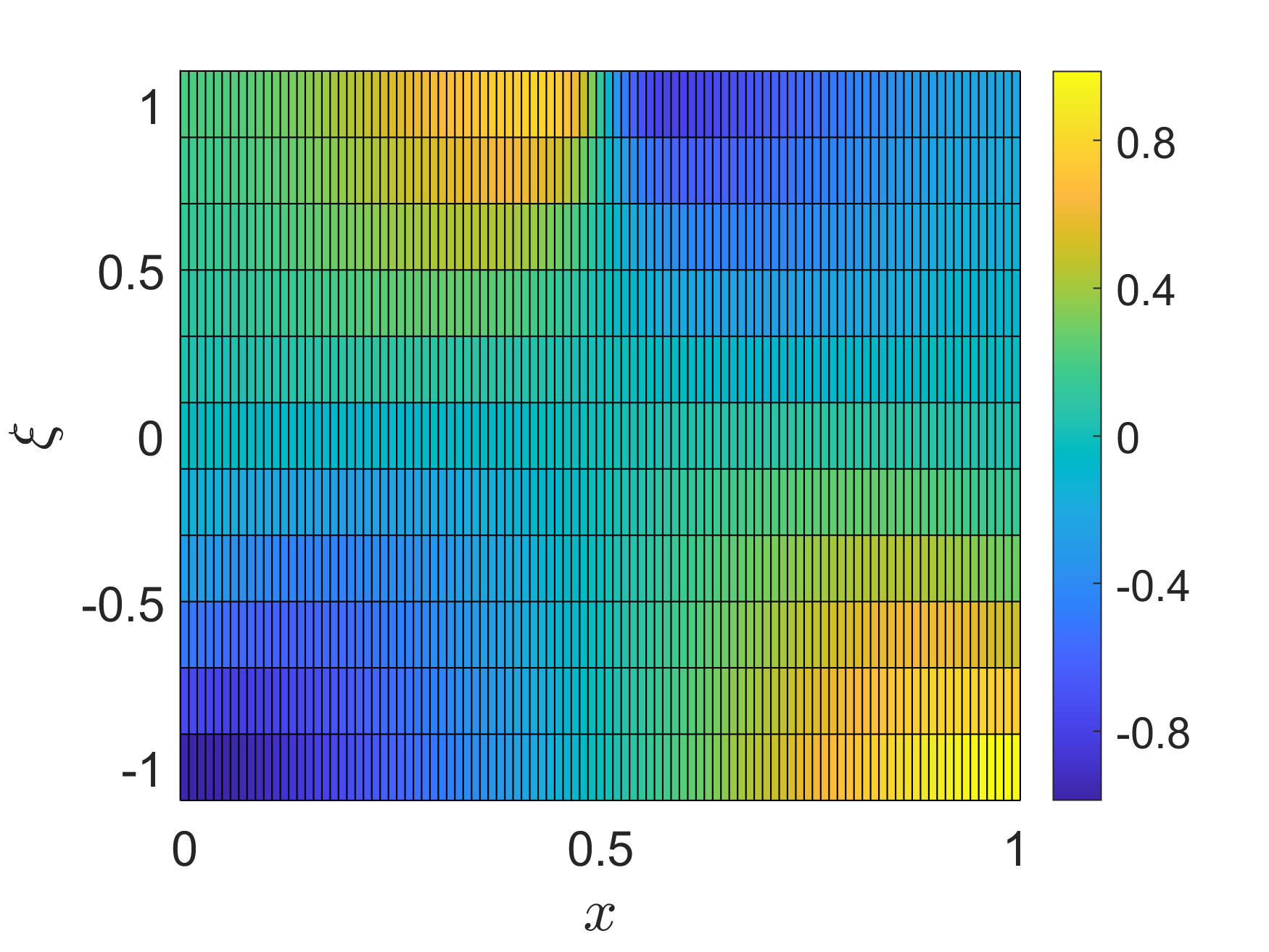}\hspace{0.5cm}
            \includegraphics[trim=0.1cm 0.2cm 0.9cm 0.8cm, clip, width=5.cm]{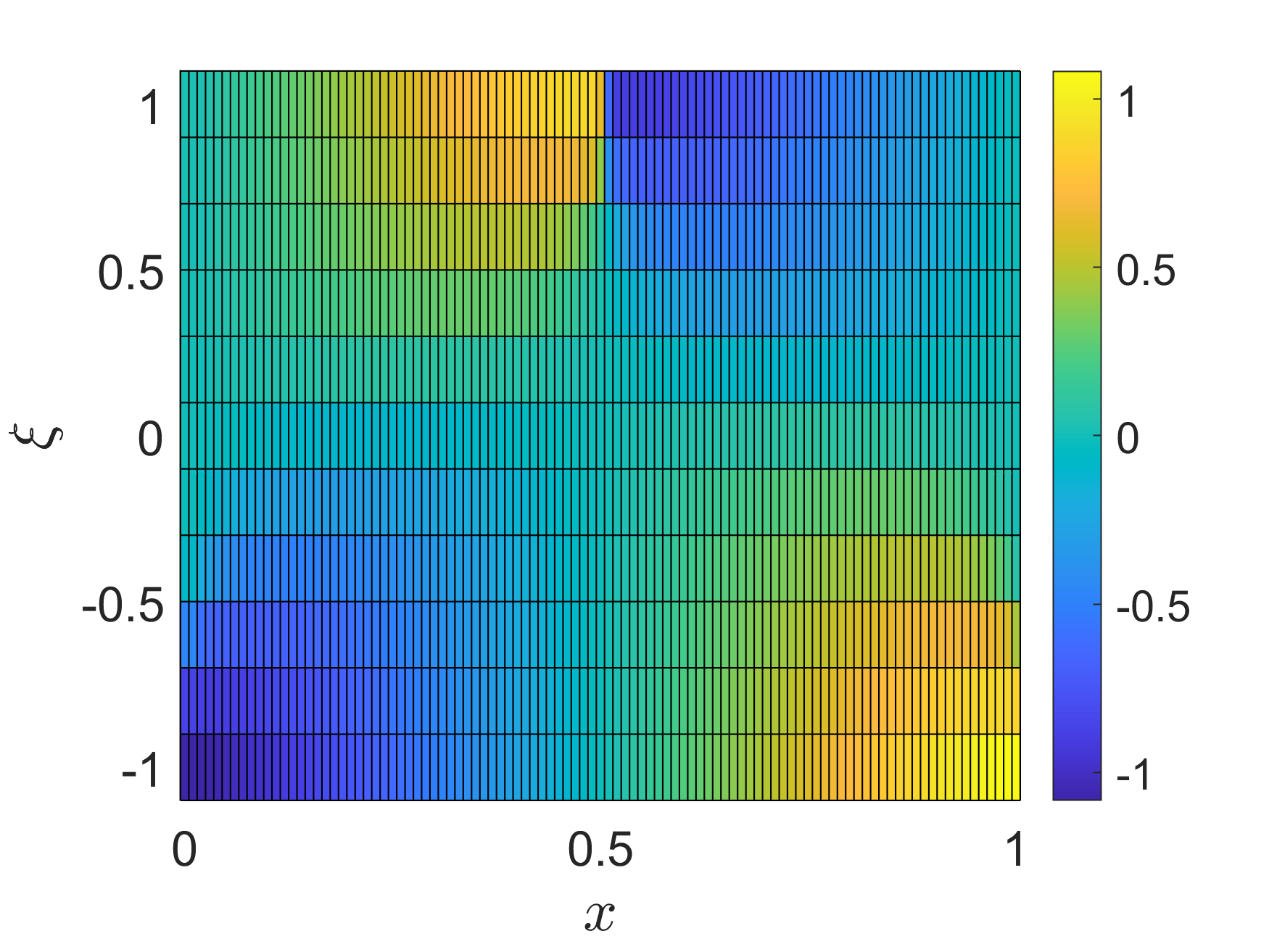} \hspace{0.5cm}
            \includegraphics[trim=0.1cm 0.2cm 0.9cm 0.8cm, clip, width=5.cm]{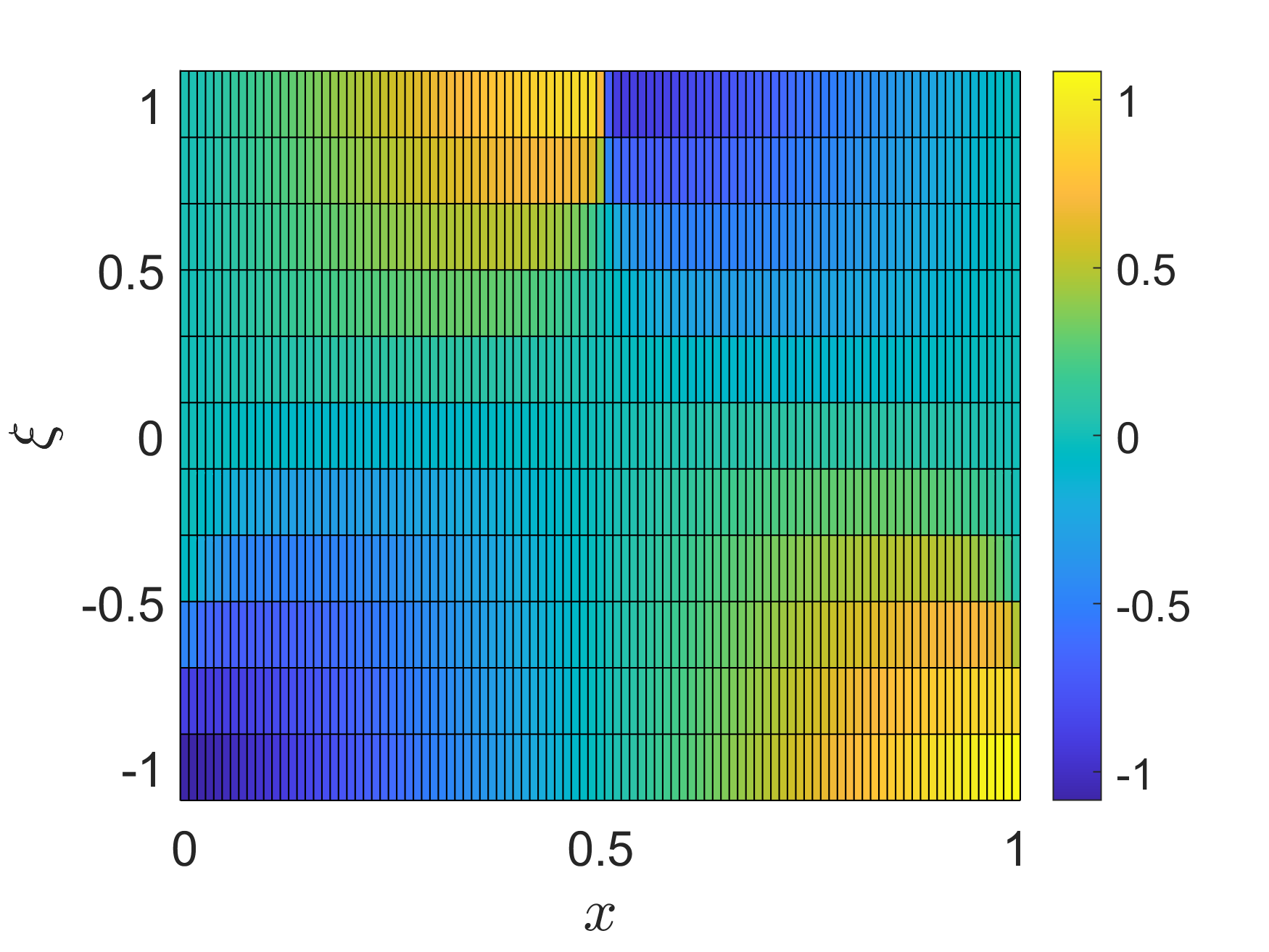}}
\caption{\sf Example 1: Numerical results computed by the 1-Order (left column), 2-Order (middle column), and 5-Order (right column) Young-measure (bottom row) schemes using the KKT-based solver.}
\label{fig1.1c}
\end{figure}

As mentioned in Remark 2.2, we propose an alternative method for solving the linear programming problem using the KKT-based solver. In this example, we show the results computed by the alternative method in Figure \ref{fig1.1c}, where one can see that the alternative method works well. 

\paragraph*{Example 2.}
In order to check the orders of accuracy of the studied schemes, we take the same initial data as in Example 1, but subject to the periodic boundary conditions. We compute the numerical solutions until the final time $t=0.05$ by the studied 1-Order, 2-Order, and 5-Order Young-measure schemes on the uniform mesh with $N_x \times N_\xi =20 \times 20$, $40 \times 40$, and $80 \times 80$. For the discretization in phase space, we also set $u \in [-1.5,1.5]$ and $N_{u} = 100$. We then compute the $L^1$-errors between the numerical and exact solution, which is given by 
$$
u_{\rm exact}(x,\xi,t)=\xi\sin\bigl(2\pi a(x,\xi,t)\bigr),
$$
where $a=a(x,\xi,t)$ is determined by solving $a+t\,\xi\sin(2\pi a)\equiv x$. We then estimate the experimental convergence rates. The obtained results are reported in Table \ref{tab1}, where one can clearly observe the expected orders of accuracy. 

\begin{table}[ht!]
\centering
\begin{tabular}{ccccccc}
\toprule 
\multirow{2}{*}{$N_x\times N_\xi$} & \multicolumn{2}{c}{1-Order} & \multicolumn{2}{c}{2-Order}& \multicolumn{2}{c}{5-Order}\\
\cline{2-7}
& Error & Rate & Error & Rate & Error & Rate\\
\hline 
$20\times 20$  & 4.51e-02 & ---  & 3.50e-04 & ---  & 1.01e-05 & ---\\
$40\times 40$  & 1.94e-02 & 1.22 & 7.84e-05 & 2.16 & 2.66e-07 & 5.25\\
$80\times 80$  & 8.84e-03 & 1.13 & 1.88e-05 & 2.06 & 8.26e-09 & 5.01\\
\bottomrule
\end{tabular}
\caption{\sf Example 2: The $L^1$-errors experimental convergence rates for the studied 1-Order, 2-Order, and 5-Order schemes.\label{tab1}}
\end{table}

\subsubsection{Isentropic Euler System}
In this section, we consider the 1-D isentropic Euler system, which reads as 
\begin{equation}\label{5.3}
\begin{cases}
 \rho_t +q_x=0,& \\[1.ex]
 q_t +\big(\frac{q^2}{\rho}+p\big)_x =0, &
\end{cases}
\end{equation}
where $\rho$, $v:=\frac{q}{\rho}$, and $p$ are the density, velocity, and pressure, respectively. The system is completed through the following equation of state:
\begin{equation}\label{5.4}
  p=\kappa \rho^\gamma,
\end{equation}
where $\kappa$ is a constant, and $\gamma$ represents the specific heat ratio. In all of the numerical experiments reported in this section, we have used $\kappa=1$ and $\gamma=1.5$. The corresponding entropy function $\eta$ in \eref{linprog-a} is defined by $\eta=\frac{1}{2}\rho v^2+\frac{p}{\gamma-1}$.  

\begin{rmk}
As mentioned in \S \ref{sec3} and \S \ref{sec4}, we use the LCD based piecewise linear interpolation in the 2-Order scheme and WENO-Z interpolation in the 5-Order scheme.  In Appendix \ref{appc}, we provide a detailed explanation on how the average matrix $\widehat A_{i,j+1/2}$ and the corresponding matrices $R_{i,\jph}$ and $R^{-1}_{i,\jph}$ are computed for the 1-D isentropic Euler system \eref{5.3}--\eref{5.4}. 
\end{rmk}

\begin{figure}[ht!]
\centerline{\includegraphics[trim=0.6cm 0.2cm 1.1cm 0.7cm, clip, width=5.cm]{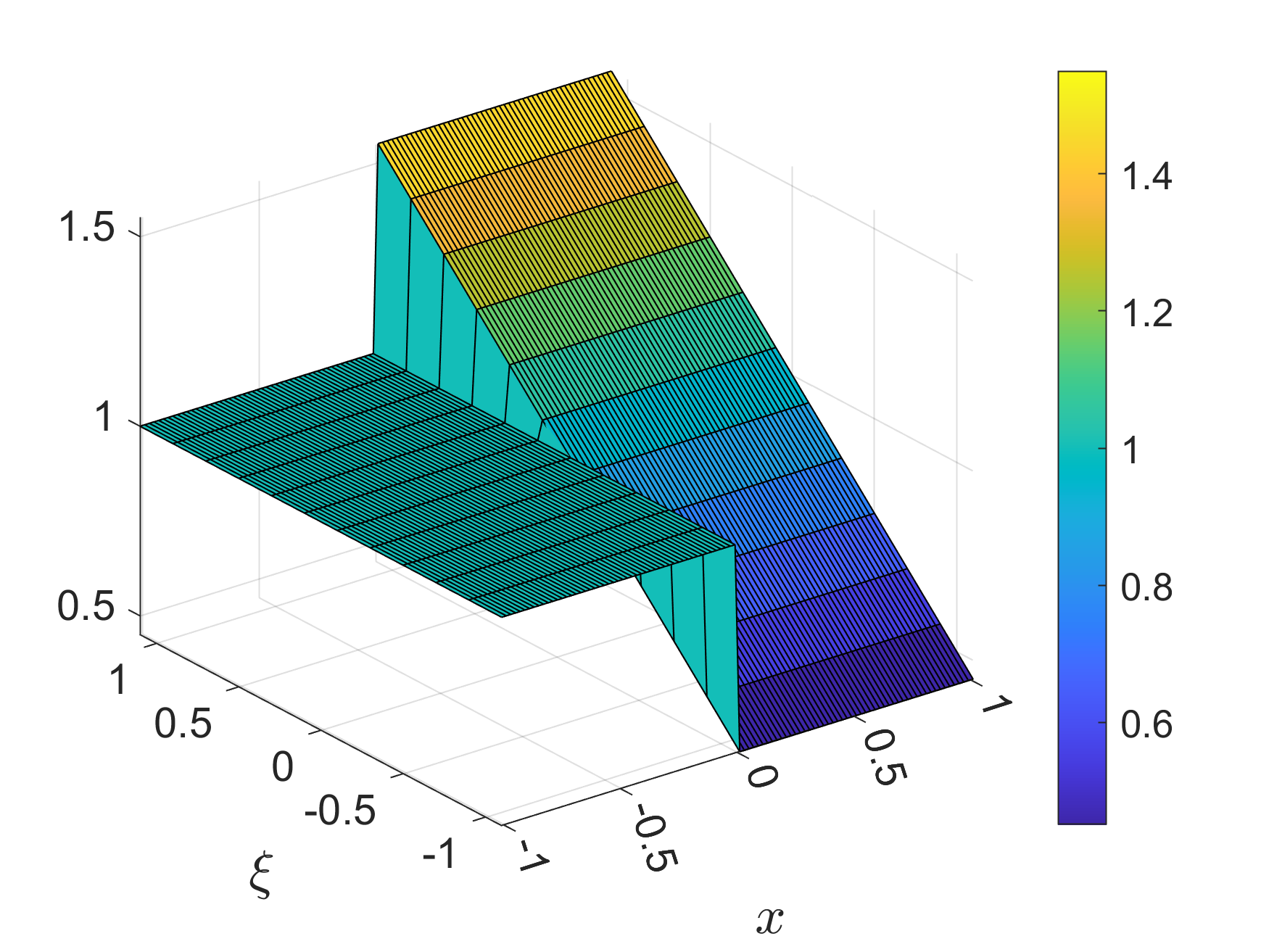} \hspace{1cm}
            \includegraphics[trim=0.6cm 0.2cm 1.1cm 0.7cm, clip, width=5.cm]{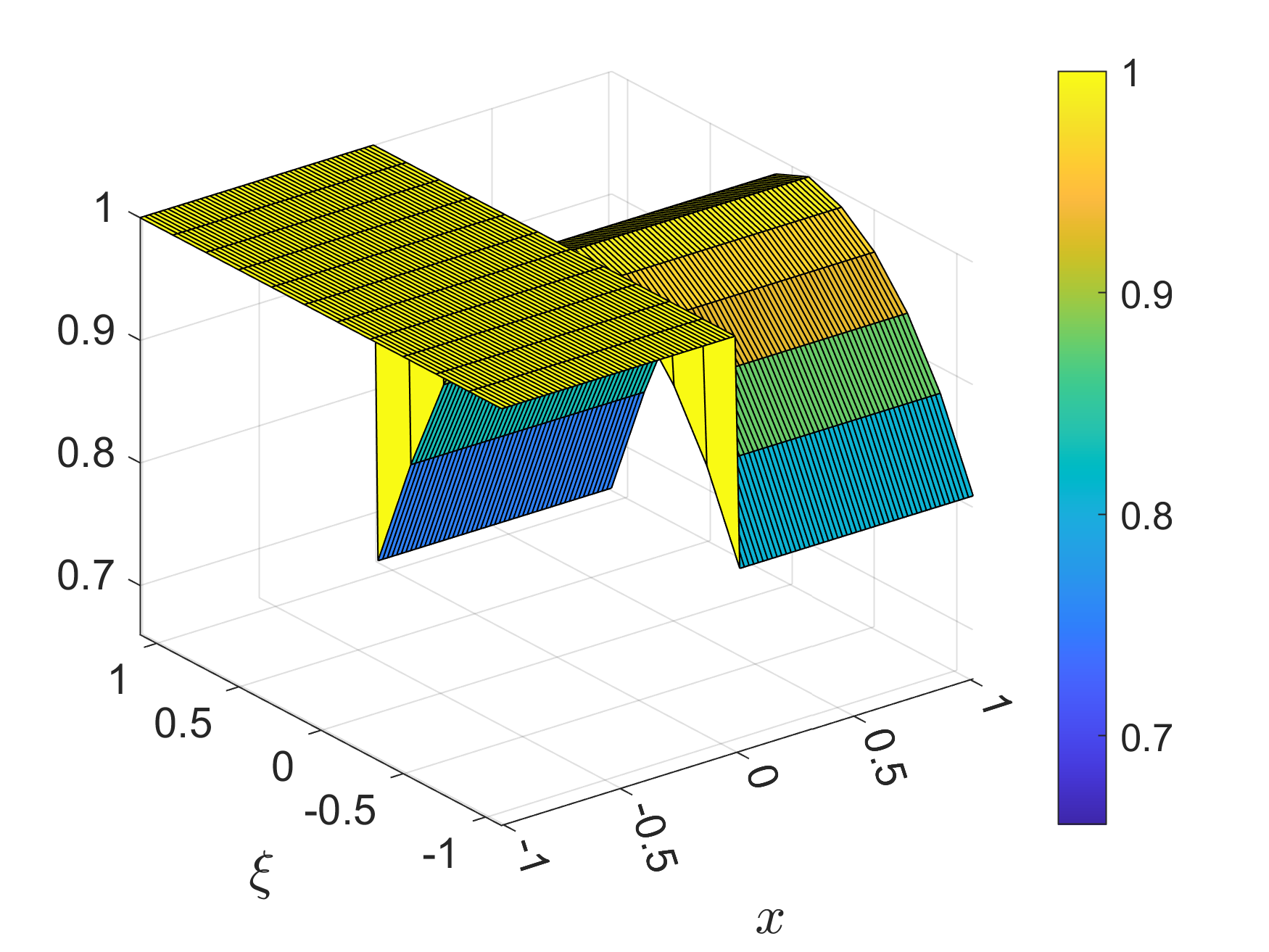}}
\caption{\sf Example 3: Initial data for $\rho$ (left) and $q$ (right).\label{fig3.3}}
\end{figure}

\paragraph*{Example 3.} In this example, taken from \cite{CHLZ25}, we consider the isentropic Euler system \eref{5.3}--\eref{5.4}. The left initial datum is deterministic and is given by
$$
\bm U_L:=(\rho_L,q_L)^\top=(1,1)^\top .
$$
The right initial datum depends on the random parameter $\xi$. More precisely, we set
$$
s(\xi):=\rho_L+\frac{1}{2}\xi=1+\frac{1}{2}\xi,
$$
and define
$$
\bm u_0(x,\xi)=\begin{cases}
\bm U_L, & x<0,\\[1mm]
\bm U_R(\xi), & x>0,
\end{cases}
$$
where
$$
\bm U_R(\xi):=\bigl(s(\xi),q_R(\xi)\bigr)^\top .
$$
The momentum \(q_R(\xi)\) is computed by
$$
q_R(\xi)=
\begin{cases}
s(\xi)\rho_L-\sqrt{\dfrac{s(\xi)}{\rho_L}\bigl(s(\xi)-\rho_L\bigr)\bigl(p(s(\xi))-p(\rho_L)\bigr)},& s(\xi)\geq \rho_L,\\[3mm]
s(\xi)\rho_L-s(\xi)\bigl(\log(s(\xi))-\log(\rho_L)\bigr),& 0<s(\xi)<\rho_L .
\end{cases}
$$
The initial conditions are prescribed in the computational domain $[-1,1]\times[-1,1]$ subject to the free boundary conditions. In Figure \ref{fig3.3}, we plot the initial data for both $\rho$ and $q$.

\begin{figure}[ht!]
\centerline{\includegraphics[trim=0.6cm 0.2cm 1.1cm 0.8cm, clip, width=5.cm]{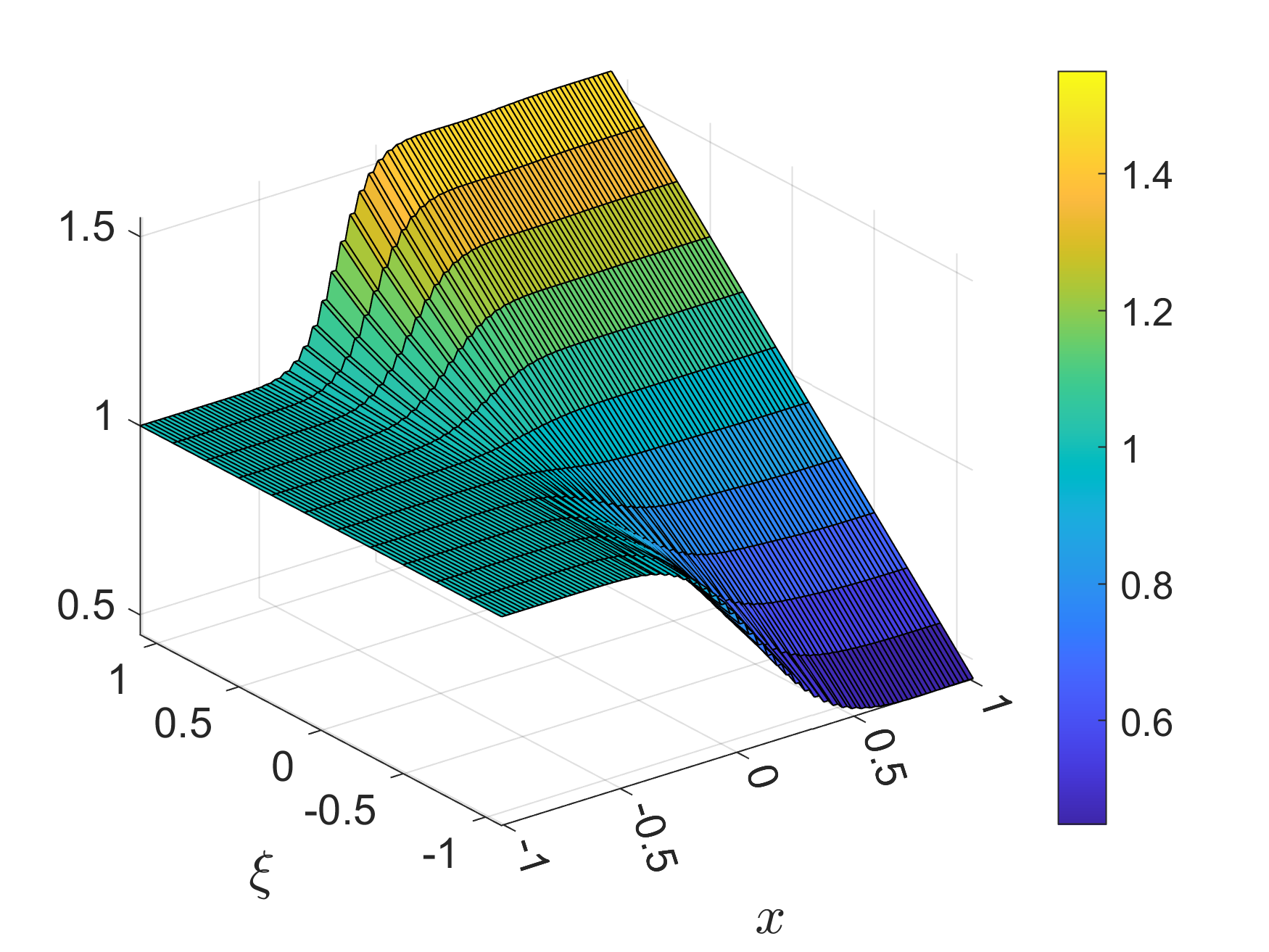} \hspace{0.5cm}
            \includegraphics[trim=0.6cm 0.2cm 1.1cm 0.8cm, clip, width=5.cm]{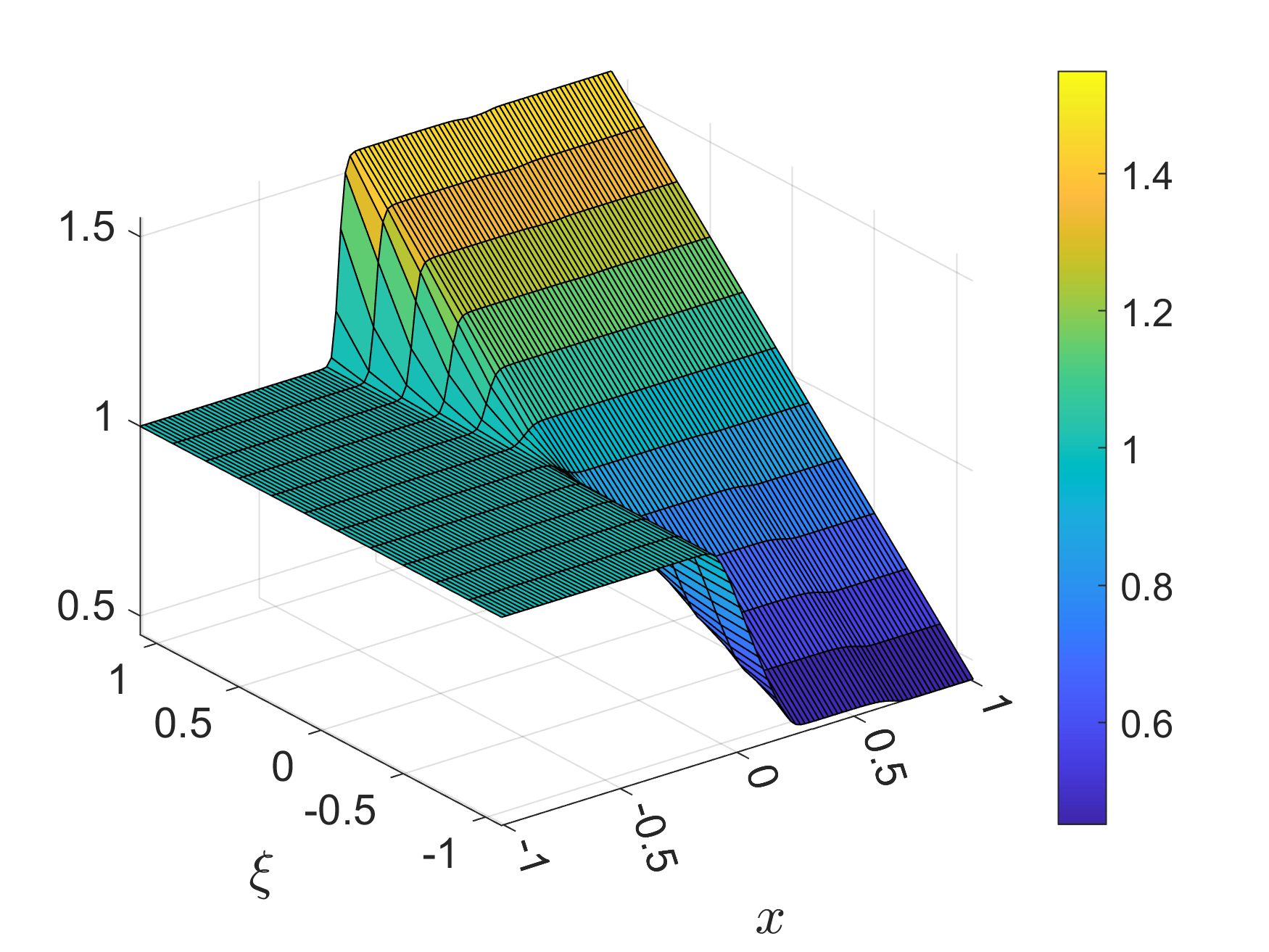} \hspace{0.5cm}
            \includegraphics[trim=0.6cm 0.2cm 1.1cm 0.8cm, clip, width=5.cm]{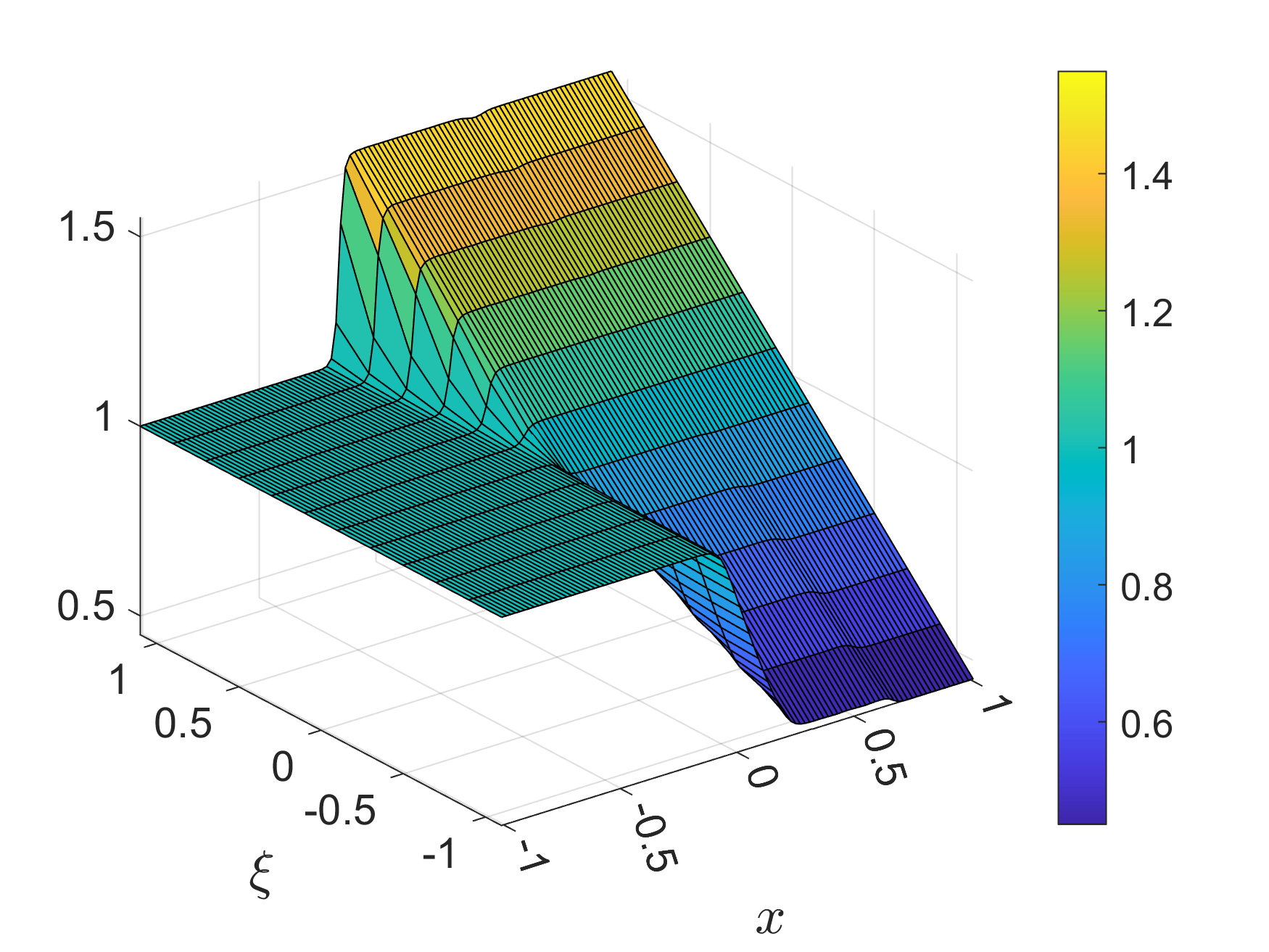}}
\vskip 10pt
\centerline{\includegraphics[trim=0.6cm 0.2cm 1.1cm 0.8cm, clip, width=5.cm]{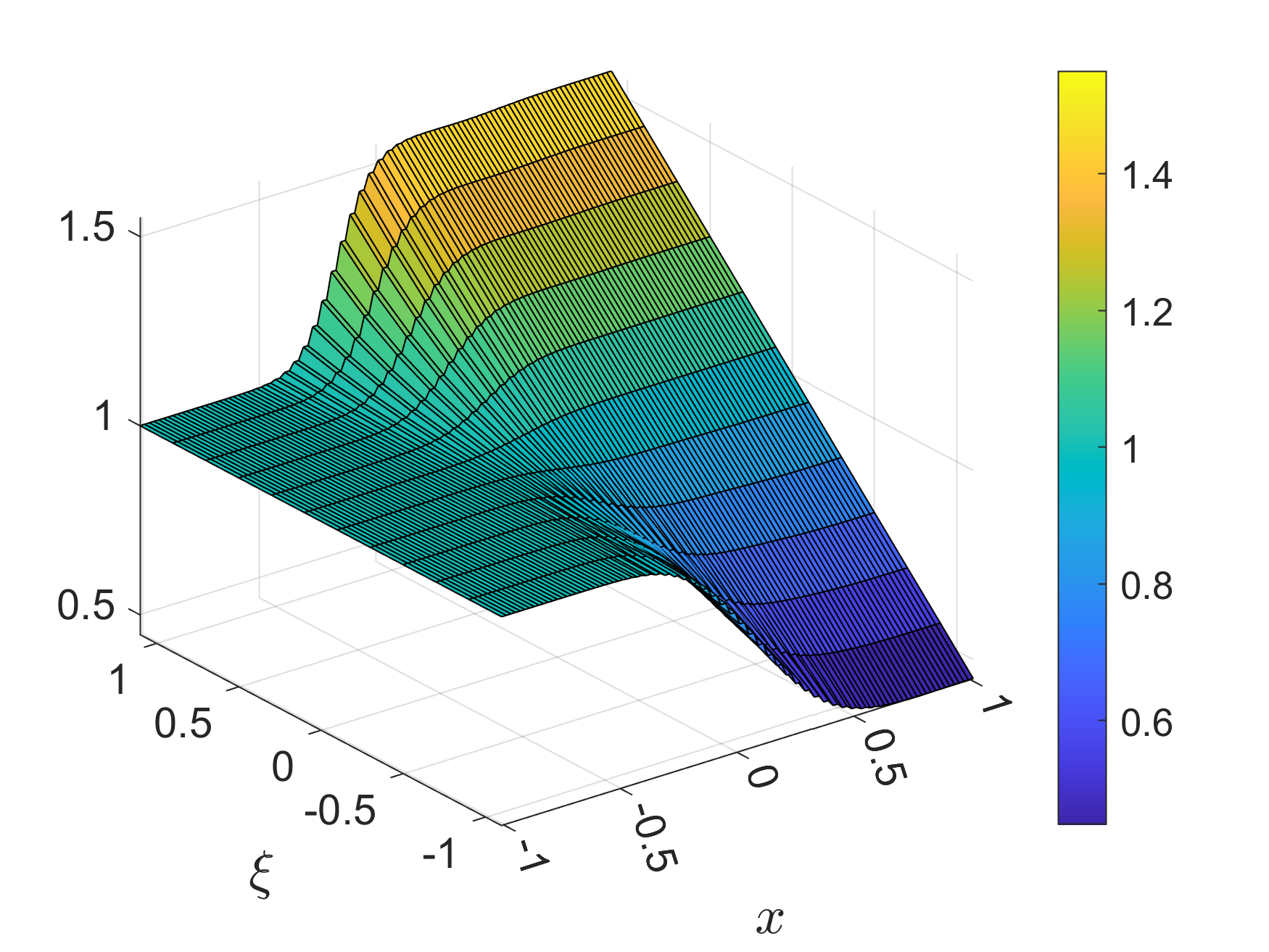} \hspace{0.5cm}
            \includegraphics[trim=0.6cm 0.2cm 1.1cm 0.8cm, clip, width=5.cm]{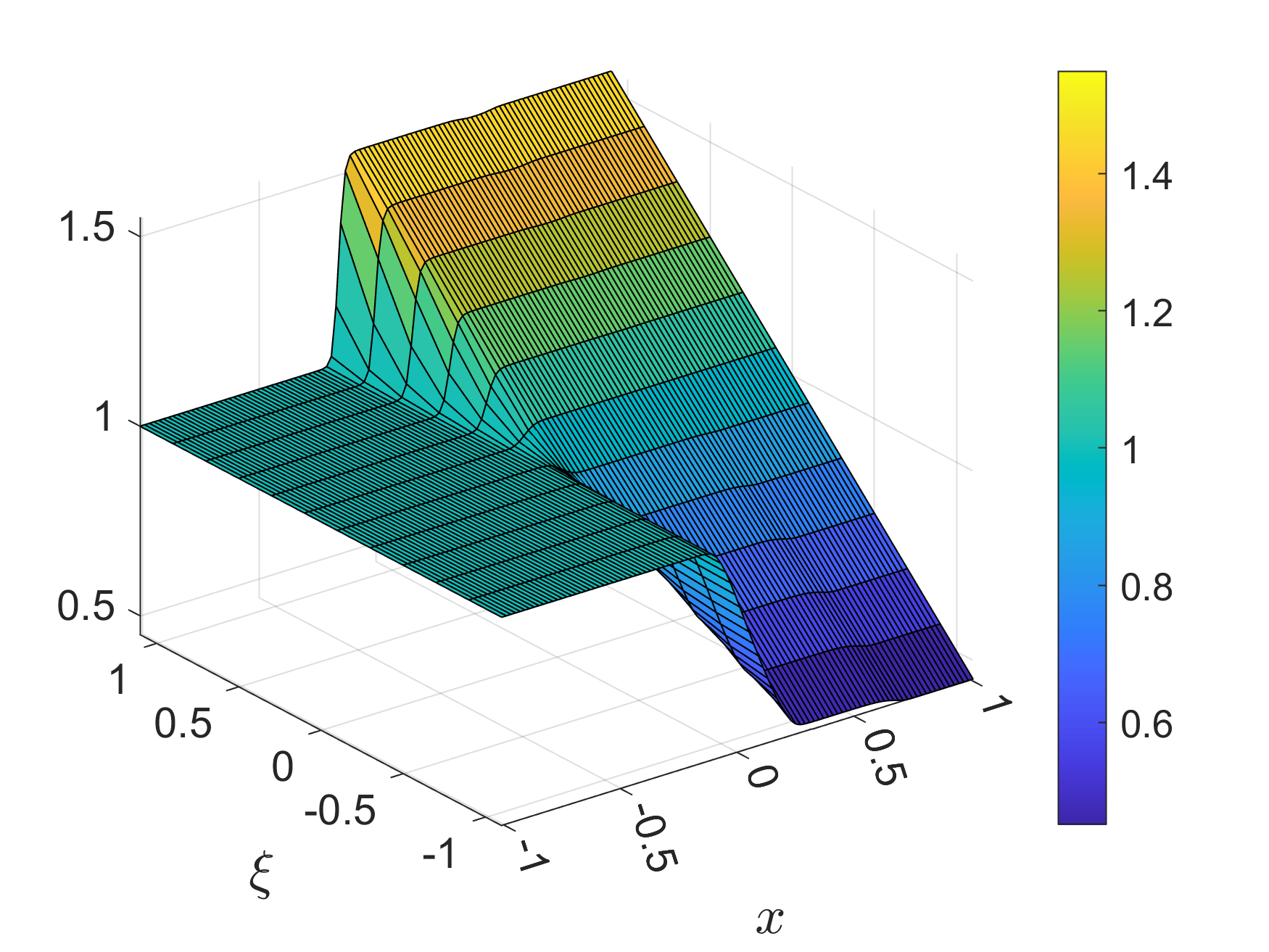} \hspace{0.5cm}
            \includegraphics[trim=0.6cm 0.2cm 1.1cm 0.8cm, clip, width=5.cm]{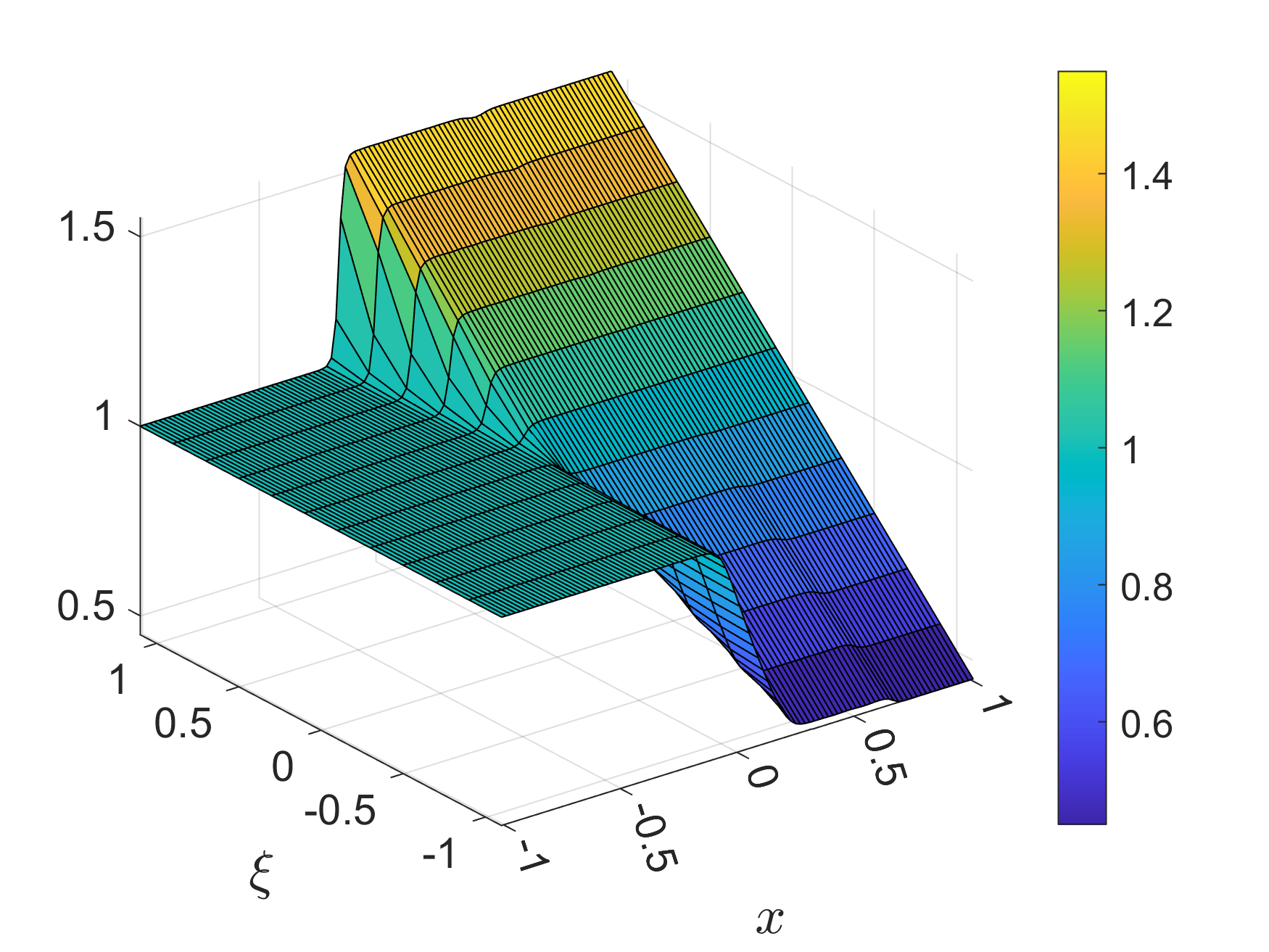}}           
\caption{\sf Example 3: Numerical results of $\rho$ computed by the 1-Order (left column), 2-Order (middle column), and 5-Order (right column) collocation (top row) and Young-measure (bottom row) schemes. \label{fig3.3a}}
\end{figure}

\begin{figure}[ht!]
\centerline{\includegraphics[trim=0.6cm 0.2cm 1.1cm 0.8cm, clip, width=5.cm]{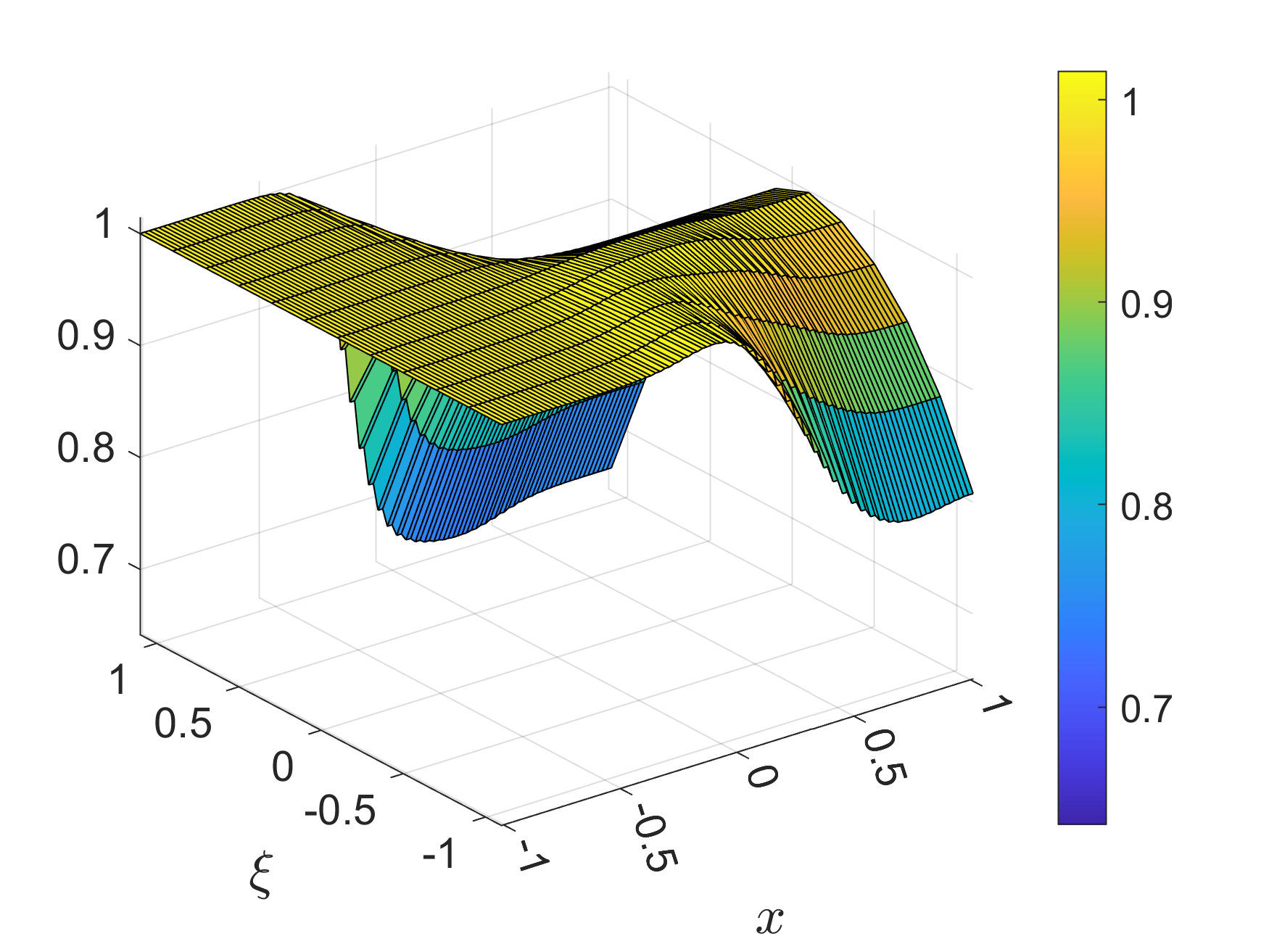} \hspace{0.5cm}
            \includegraphics[trim=0.6cm 0.2cm 1.1cm 0.8cm, clip, width=5.cm]{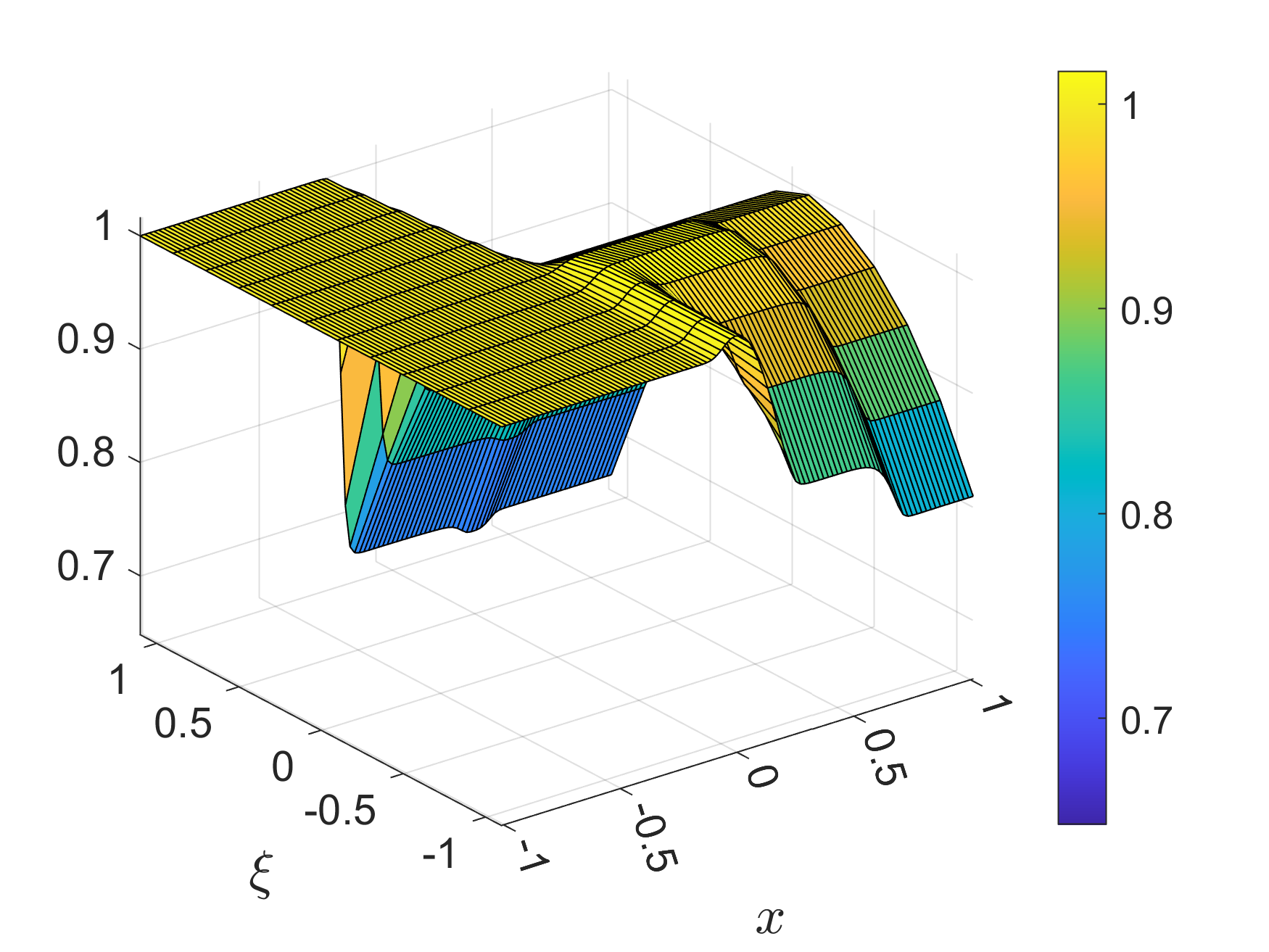} \hspace{0.5cm}
            \includegraphics[trim=0.6cm 0.2cm 1.1cm 0.8cm, clip, width=5.cm]{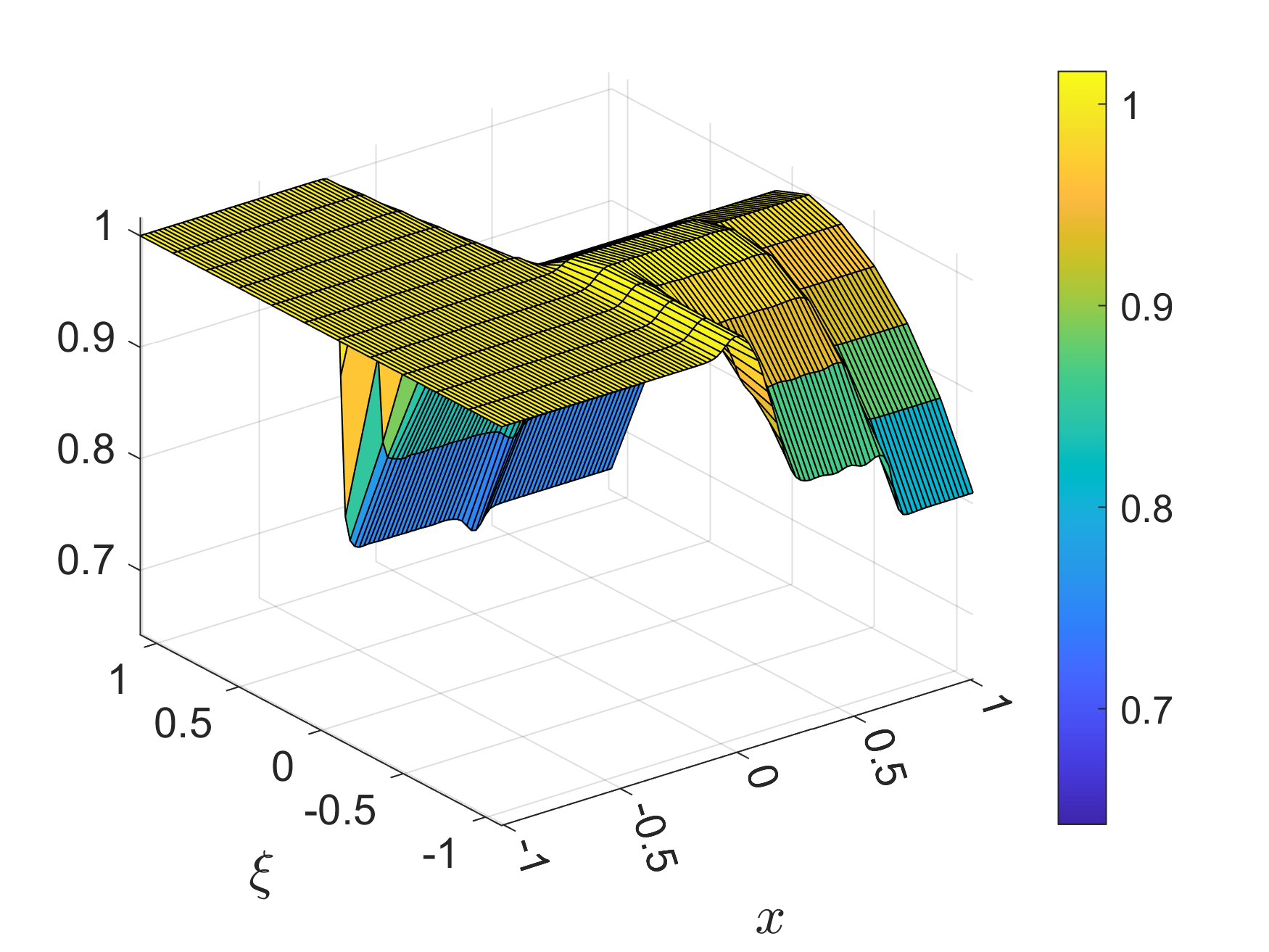}}          
\vskip 10pt
\centerline{\includegraphics[trim=0.6cm 0.2cm 1.1cm 0.8cm, clip, width=5.cm]{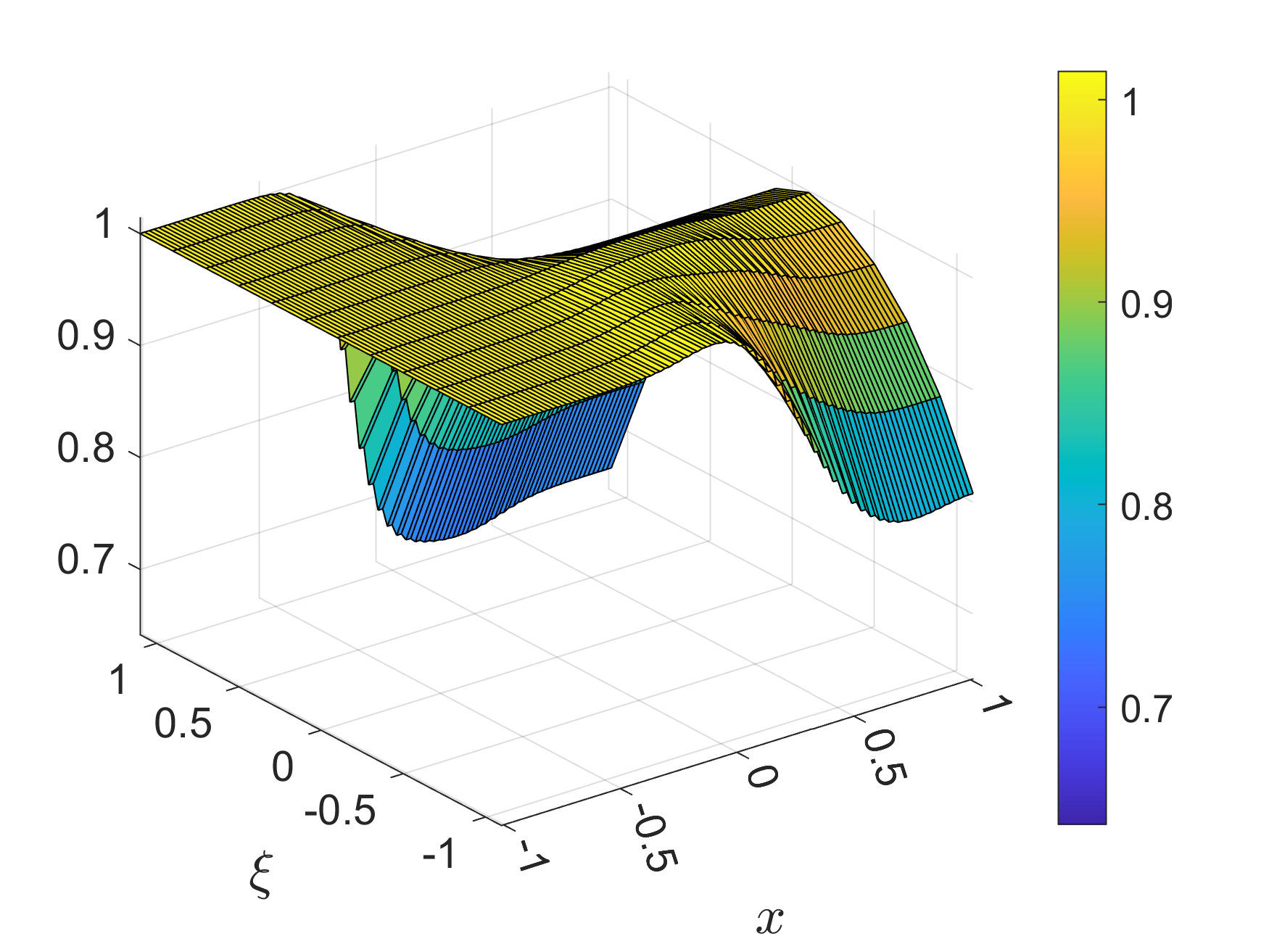}\hspace{0.5cm}
            \includegraphics[trim=0.6cm 0.2cm 1.1cm 0.8cm, clip, width=5.cm]{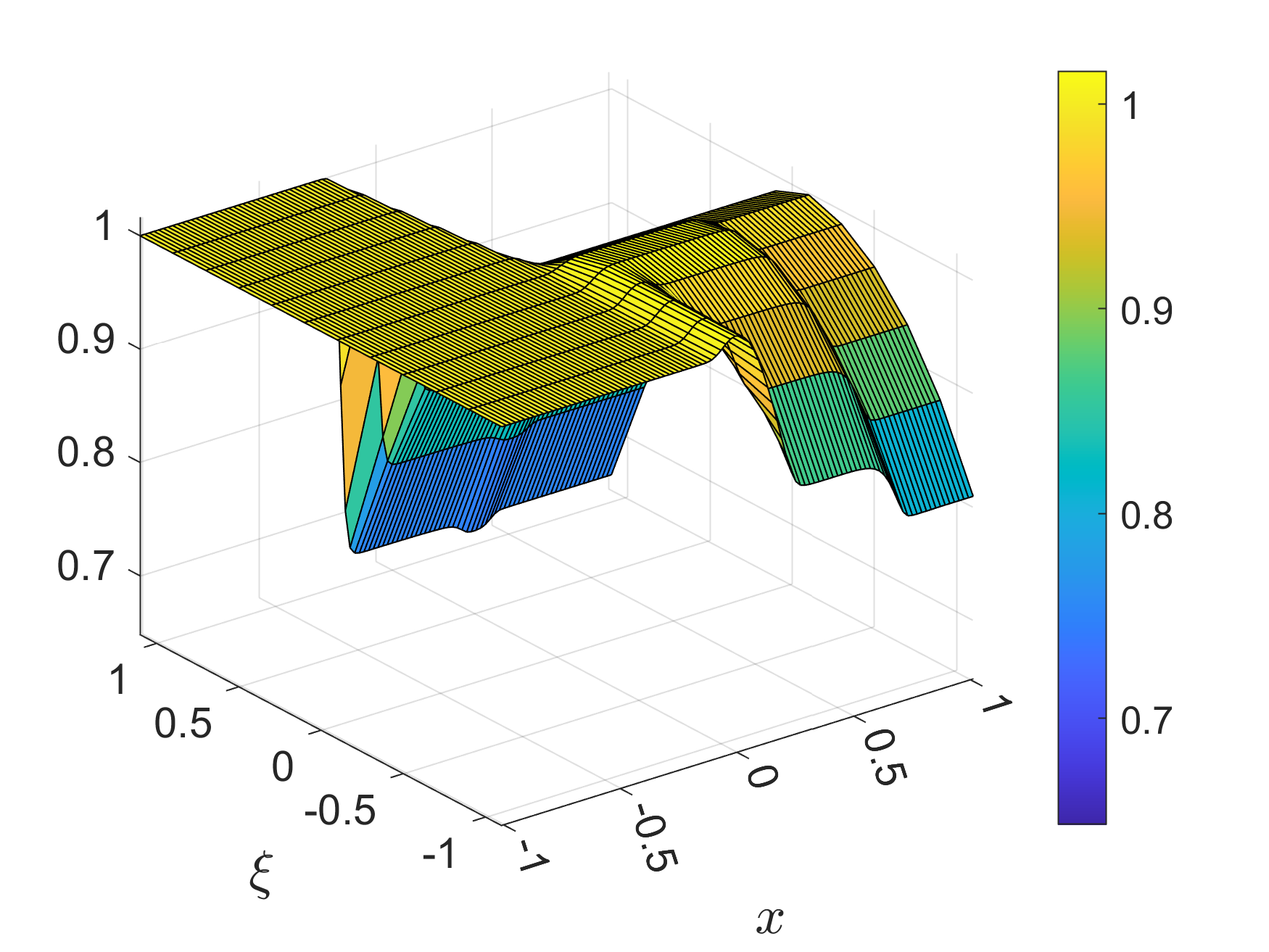}\hspace{0.5cm}
            \includegraphics[trim=0.6cm 0.2cm 1.1cm 0.8cm, clip, width=5.cm]{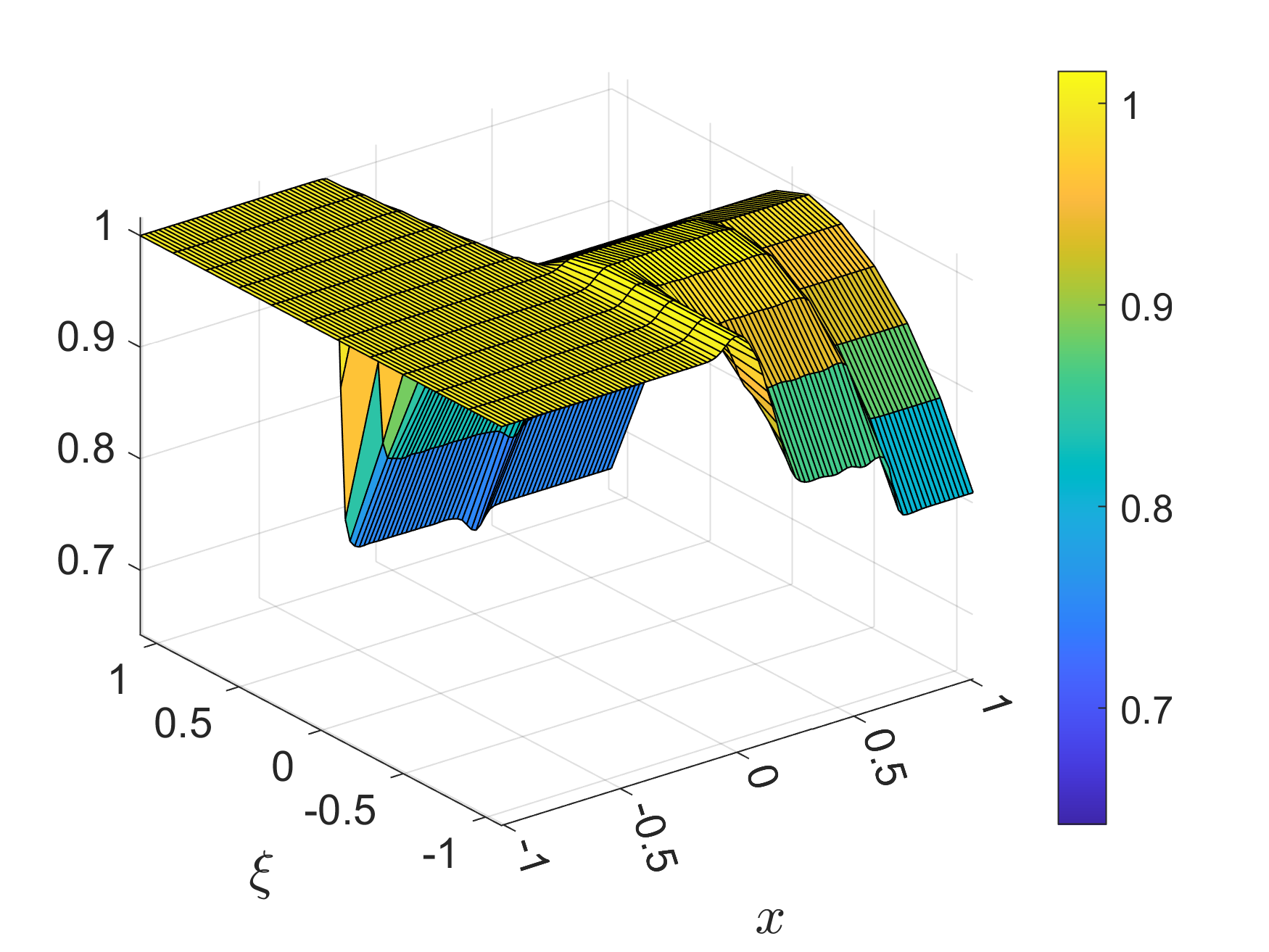}}       
\caption{\sf Example 3: Numerical results of $q$ computed by the 1-Order (left column), 2-Order (middle column), and 5-Order (right column) collocation (top row) and Young-measure (bottom row) schemes. \label{fig3.3b}}
\end{figure}

We compute the numerical solutions until the final time $t=0.25$ by the studied 1-Order, 2-Order, and 5-Order Young-measure schemes on the uniform mesh with $N_x=100$, $N_\xi=10$. For the discretization in phase space $\bmu=(\rho, q)$, we set $[0.3,1.8]\times [0.3,1.3]$ and use $N_{\bmu}=25$ grid points in each variable. We present the results computed by the collocation and Young-measure approaches in Figures \ref{fig3.3a} and \ref{fig3.3b}, where the solution exhibits a shock wave for $\xi<0$ in the first family and a rarefaction wave of the first family for all values of $\xi\geq 0$. One can clearly see that the solutions computed by the collocation and Young-measure methods agree well, and that the high-order schemes improve the resolution, particularly near the shock waves.

\medskip 

\begin{figure}[ht!]
\centerline{\includegraphics[trim=0.6cm 0.2cm 1.1cm 0.8cm, clip, width=5.cm]{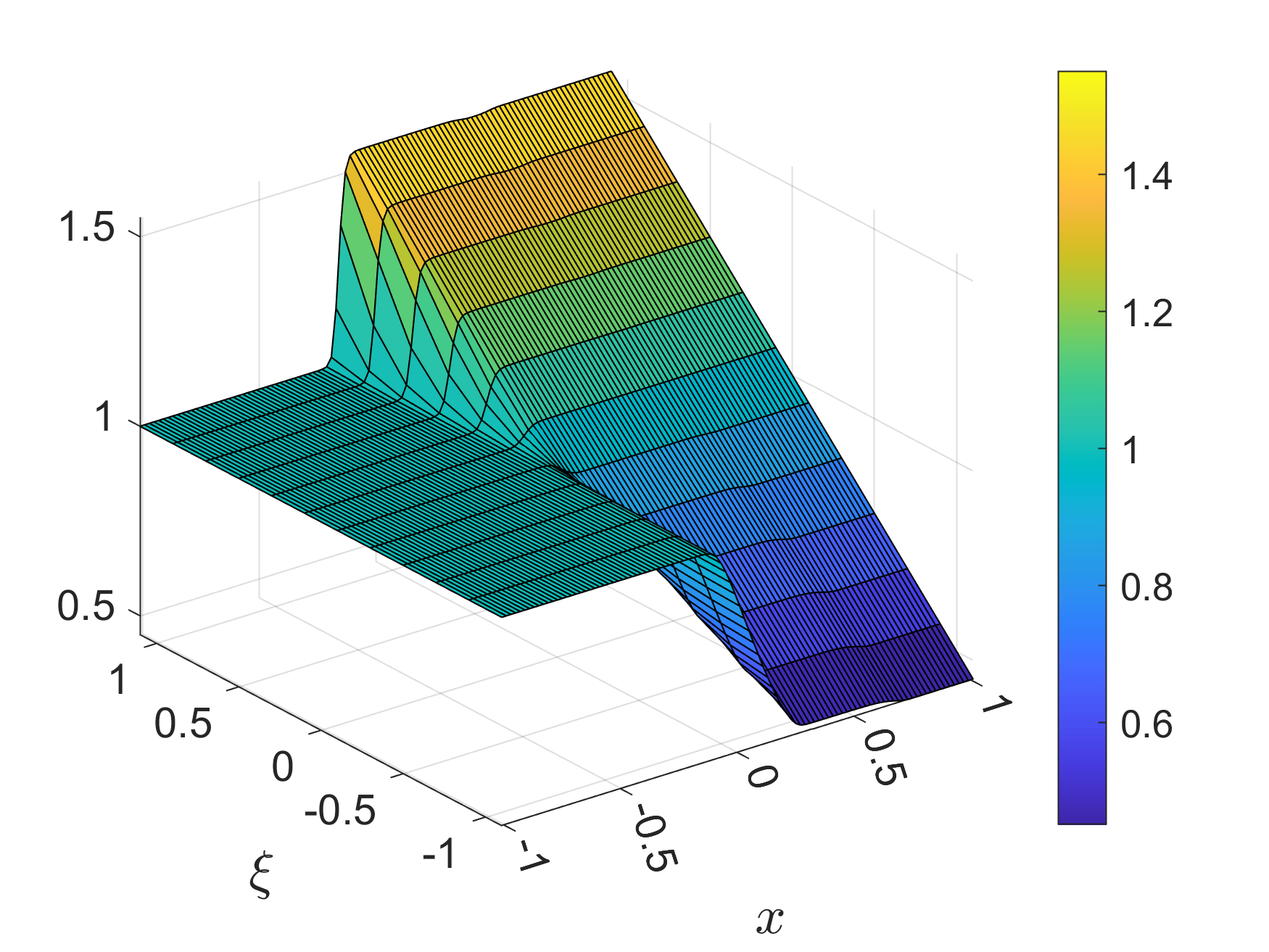} \hspace{1cm}
            \includegraphics[trim=0.6cm 0.2cm 1.1cm 0.8cm, clip, width=5.cm]{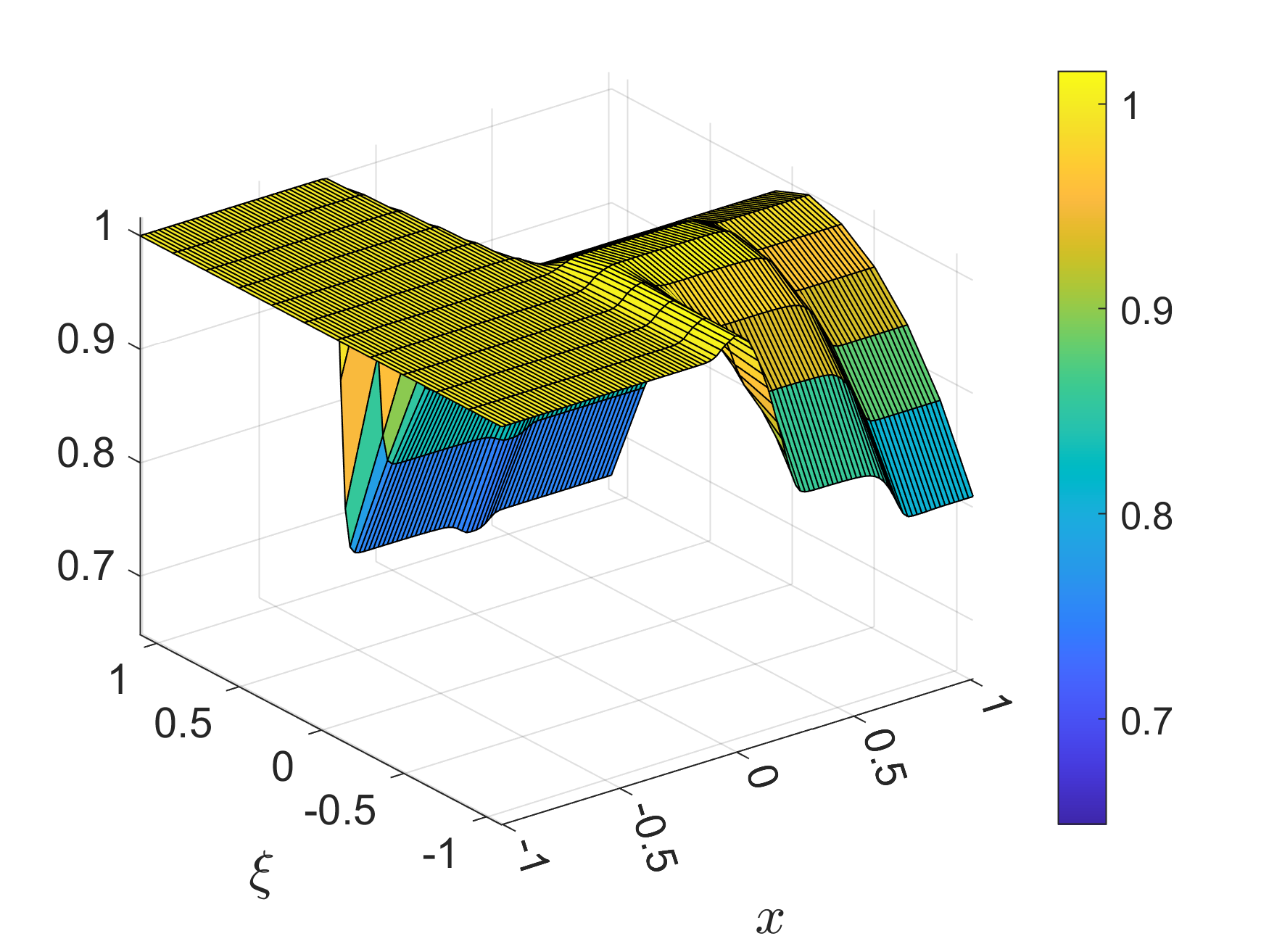}}
\vskip 10pt
\centerline{\includegraphics[trim=0.6cm 0.2cm 1.1cm 0.8cm, clip, width=5.cm]{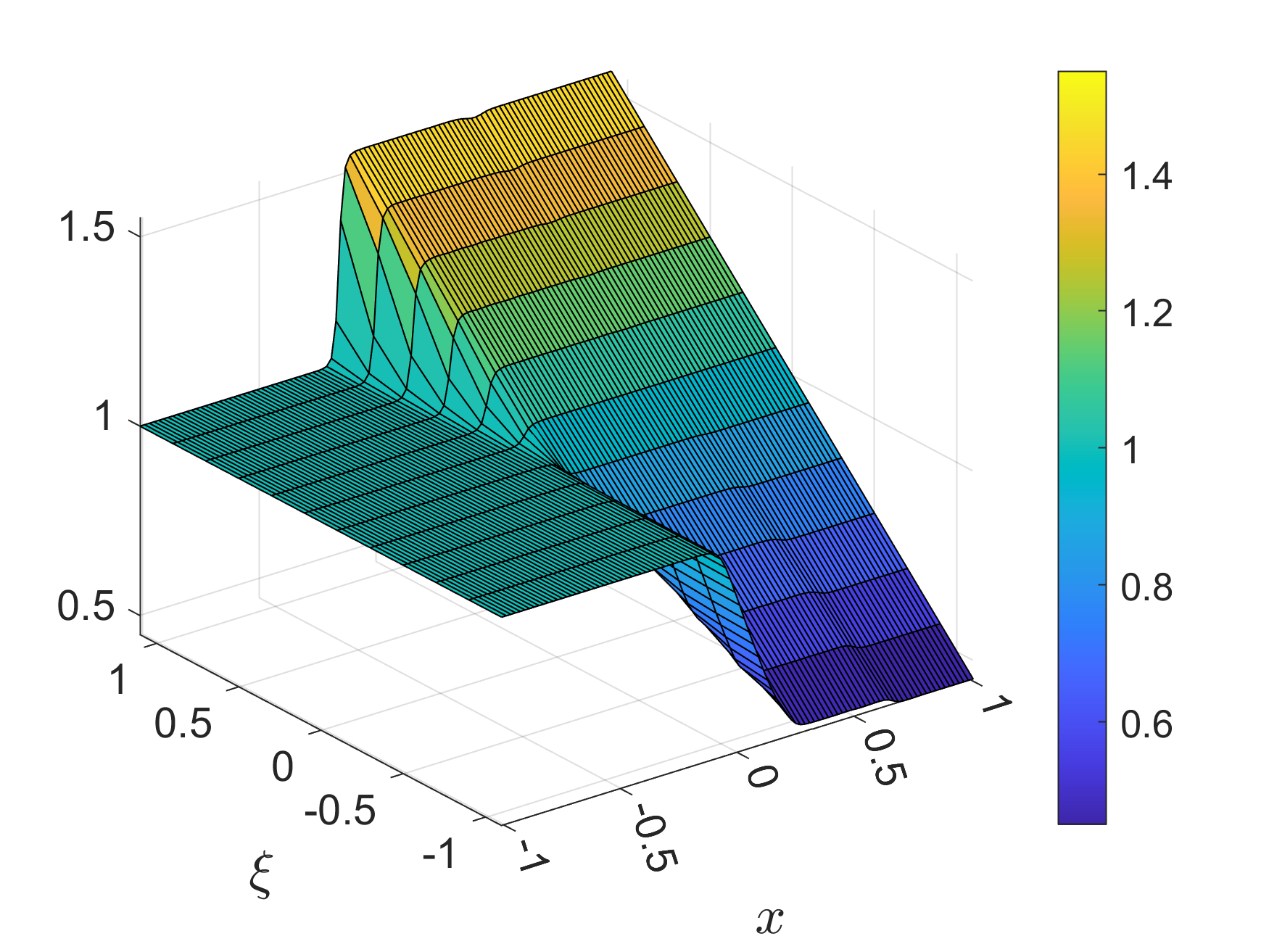} \hspace{1cm}
            \includegraphics[trim=0.6cm 0.2cm 1.1cm 0.8cm, clip, width=5.cm]{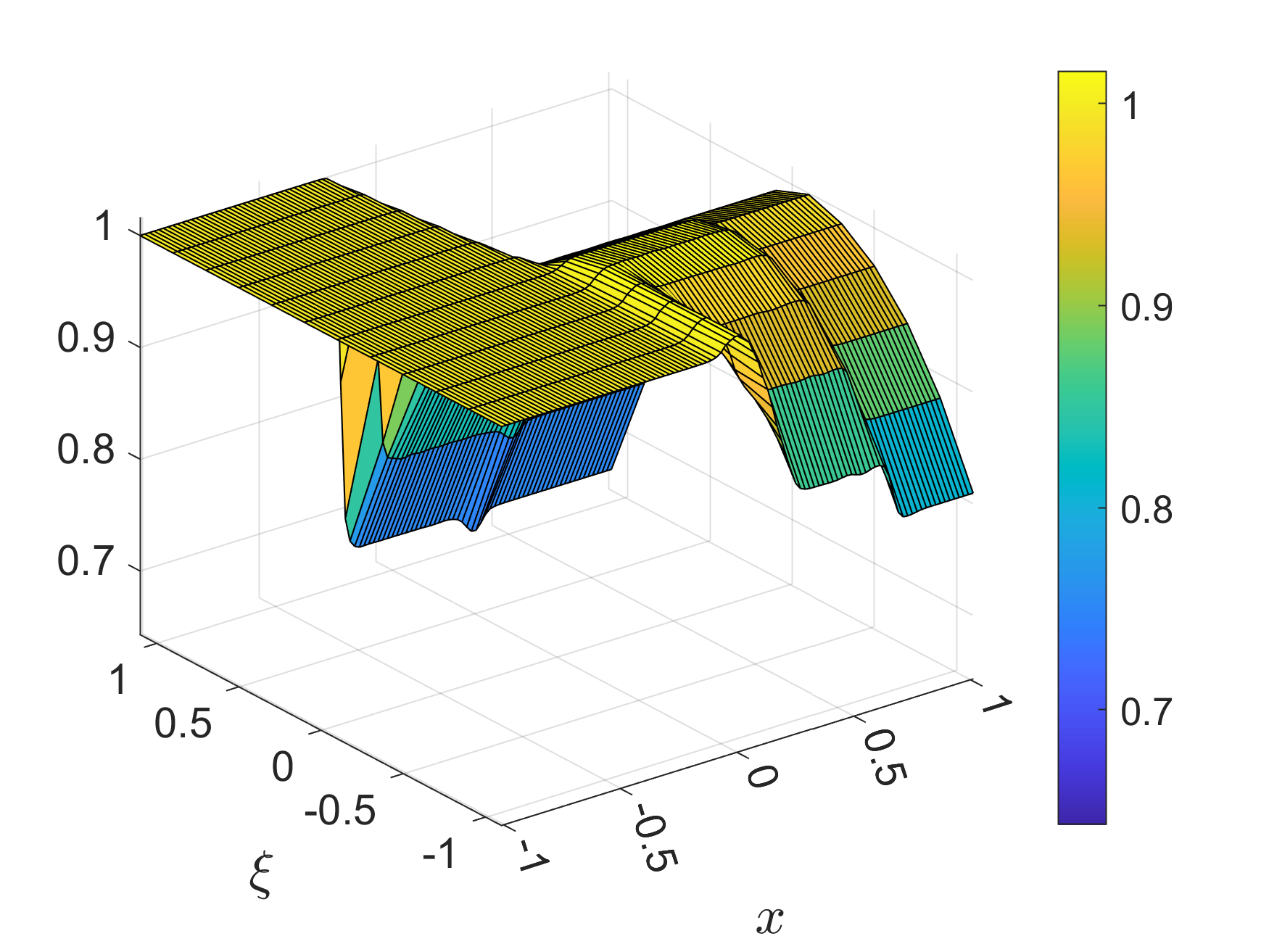}}
\caption{\sf Example 3: Numerical results of $\rho$ (left column) and $q$ (right column) computed by the 2-Order (top row) and 5-Order (bottom row) Young-measure
schemes using the numerical fluxes defined in \eref{3.7}--\eref{3.8}. \label{fig3.9c}}
\end{figure}

As mentioned in Remarks 3.1 and 4.2, one can also use the numerical fluxes defined in  \eref{3.7}--\eref{3.8}, and the obtained numerical results are essentially identical to those computed by the numerical fluxes given by \eref{3.2} in this test. Figure \ref{fig3.9c} shows that the computed results are visually indistinguishable from those plotted in Figures \ref{fig3.3a} and \ref{fig3.3b}.

\medskip 

\begin{figure}[ht!]
\centerline{\includegraphics[trim=1.2cm 0.4cm 1.3cm 0.2cm, clip, width=5.cm]{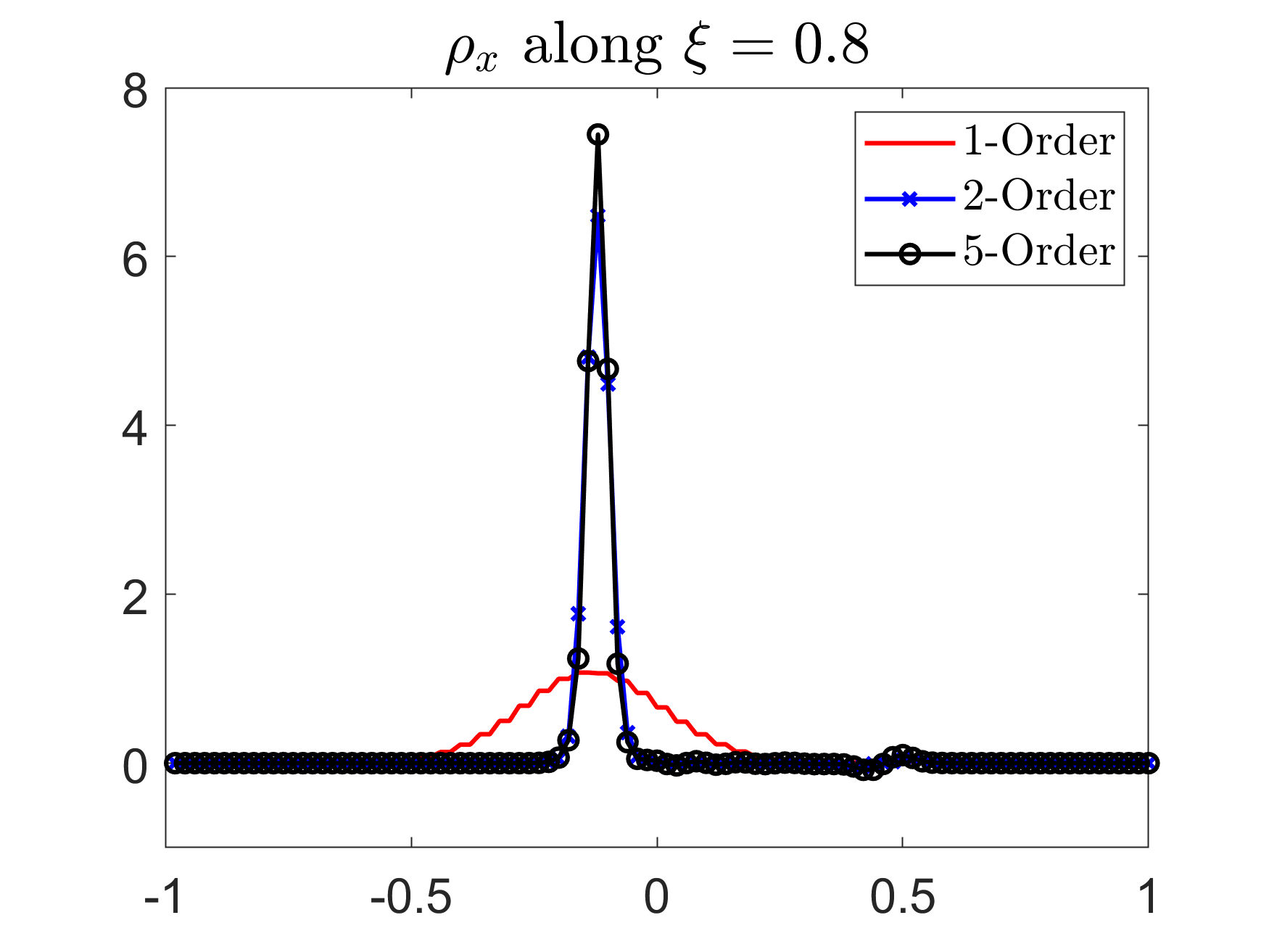}\hspace{1cm}
            \includegraphics[trim=1.2cm 0.4cm 1.3cm 0.2cm, clip, width=5.cm]{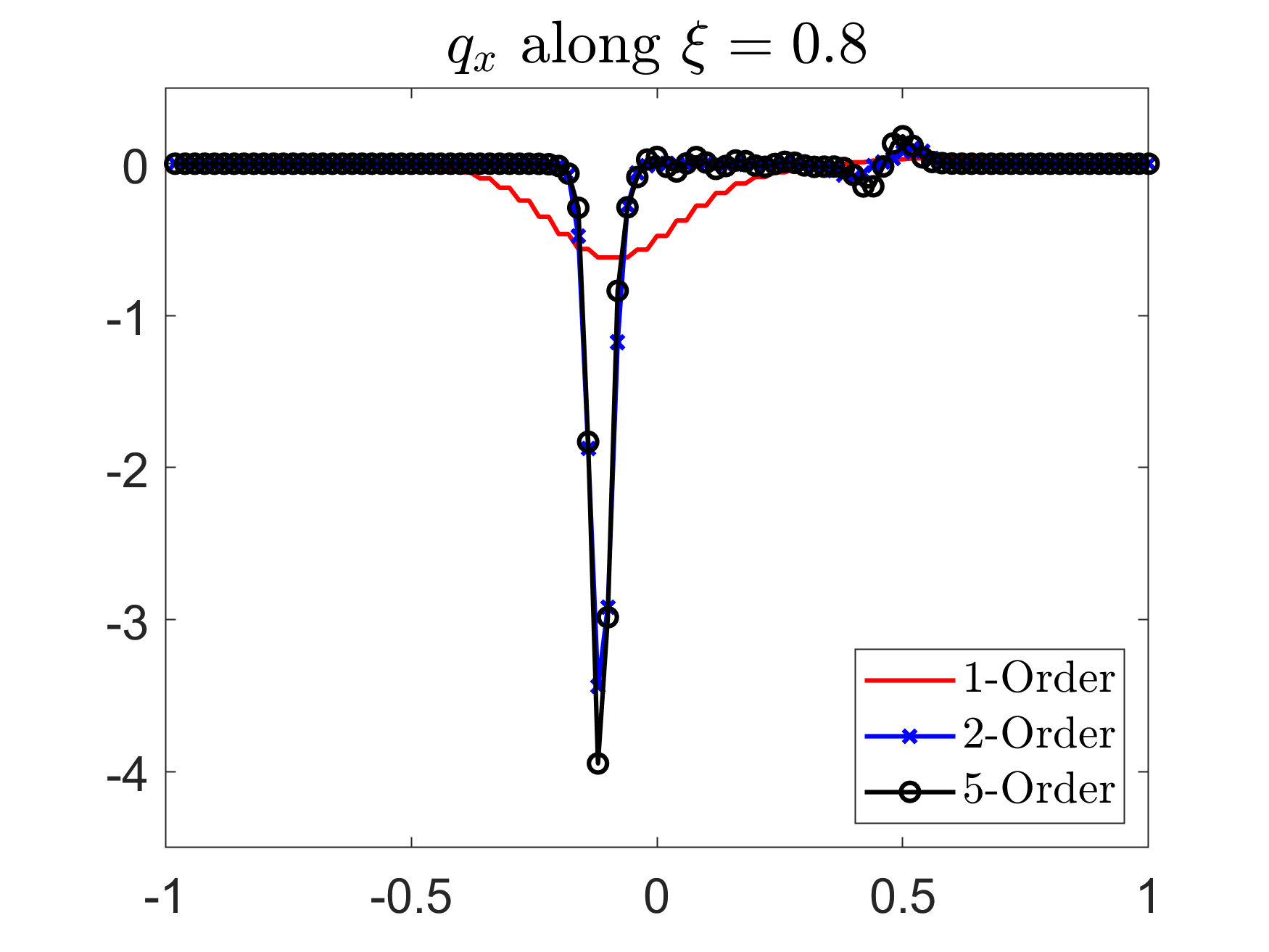}}
\caption{\sf Example 3: $\rho_x$ and $q_x$ along $\xi=0.8$ computed by the 1-Order, 2-Order, and 5-Order Young-measure schemes.}
\label{fig3.1b}
\end{figure}

In order to show the differences between the results computed by the studied schemes, similar to Example 1, we compute $\rho_x$ and $q_x$ along the line $\xi=0.8$ using the central difference scheme  and plot the obtained results in Figure \ref{fig3.1b}. One can clearly see that the 2-Order and 5-Order schemes achieve much higher resolution than the 1-Order one, while the 5-Order scheme is slightly better than the 2-Order one. 

We also estimate the $L^1$-, $L^2$-, and $L^\infty$-errors between the numerical results computed by the Young-measure methods with the reference solution, which is computed by the 1-Order collocation method on a much finer mesh with $N_x=4000$ and $N_\xi=400$. The obtained numerical results are reported in Table \ref{tab3a}, which once again show that the 2-Order and 5-Order schemes are much more accurate than the 1-Order scheme.

\begin{table}[ht!]
\centering
\begin{tabular}{cccccccccc}
\toprule 
 & \multicolumn{3}{c}{$L^1$-error} & \multicolumn{3}{c}{$L^2$-error}& \multicolumn{3}{c}{$L^\infty$-error}\\
\cline{1-10}
& 1-Order & 2-Order & 5-Order & 1-Order & 2-Order & 5-Order & 1-Order & 2-Order & 5-Order\\
\hline
$\rho$  & 6.88e-02 & 7.74e-03 & 7.07e-03 & 7.32e-02 & 2.08e-02 & 1.97e-02 & 2.25e-01 & 2.00e-01 & 1.89e-01\\
$q$     & 2.93e-02 & 4.07e-03 & 3.87e-03 & 3.43e-02 & 1.27e-02 & 1.26e-02 & 1.64e-01 & 1.39e-01 & 1.28e-01\\
\bottomrule
\end{tabular}
\caption{\sf Example 3: The $L^1$-, $L^2$-, and $L^\infty$-errors of $\rho$ and $q$ between the numerical results computed by the 1-Order, 2-Order, and 5-Order Young-measure schemes and the reference solution.\label{tab3a}}
\end{table}

\begin{figure}[ht!]
\centerline{\includegraphics[trim=1.2cm 0.4cm 1.3cm 0.2cm, clip, width=5.cm]{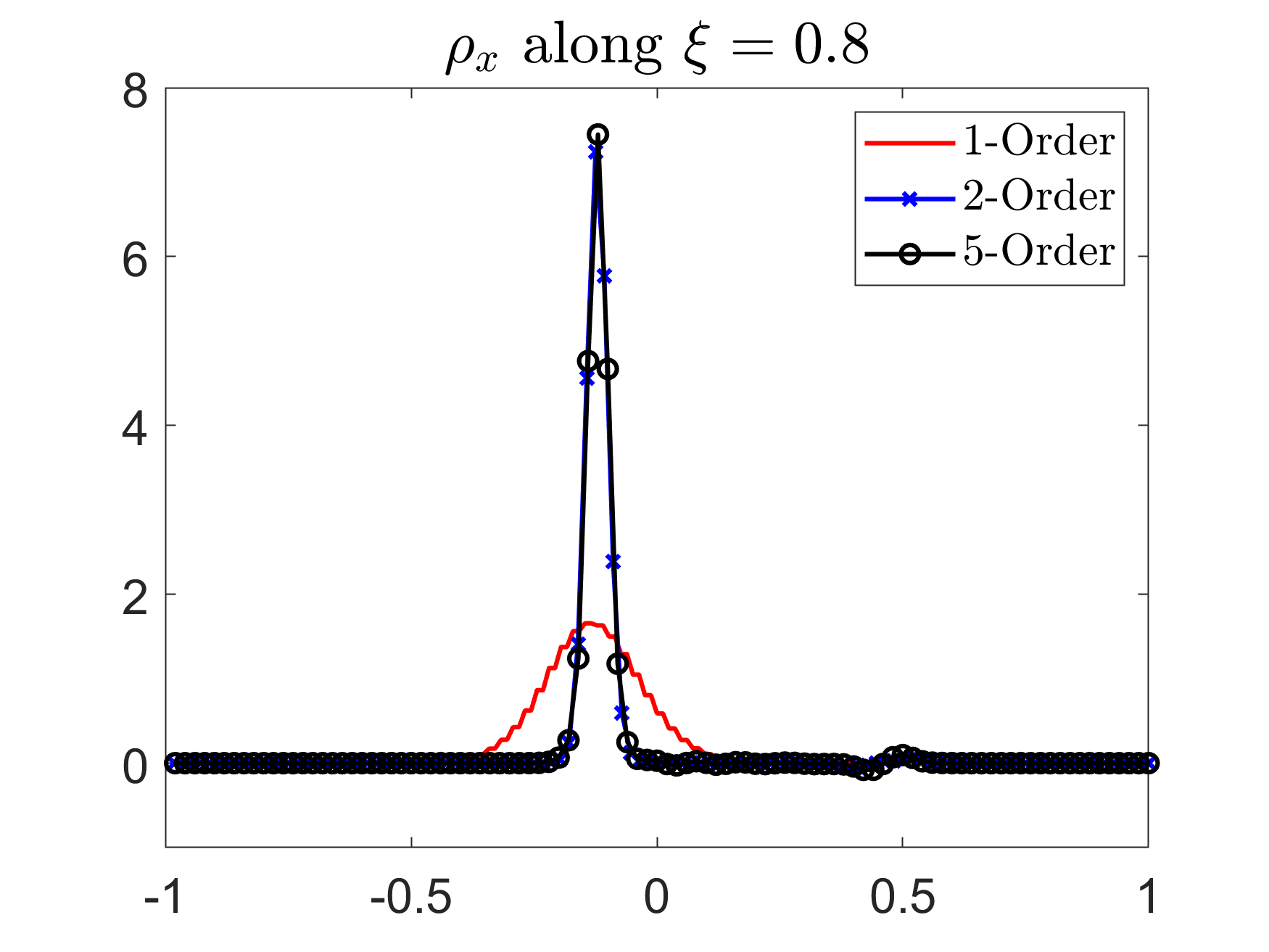}\hspace{1cm}
            \includegraphics[trim=1.2cm 0.4cm 1.3cm 0.2cm, clip, width=5.cm]{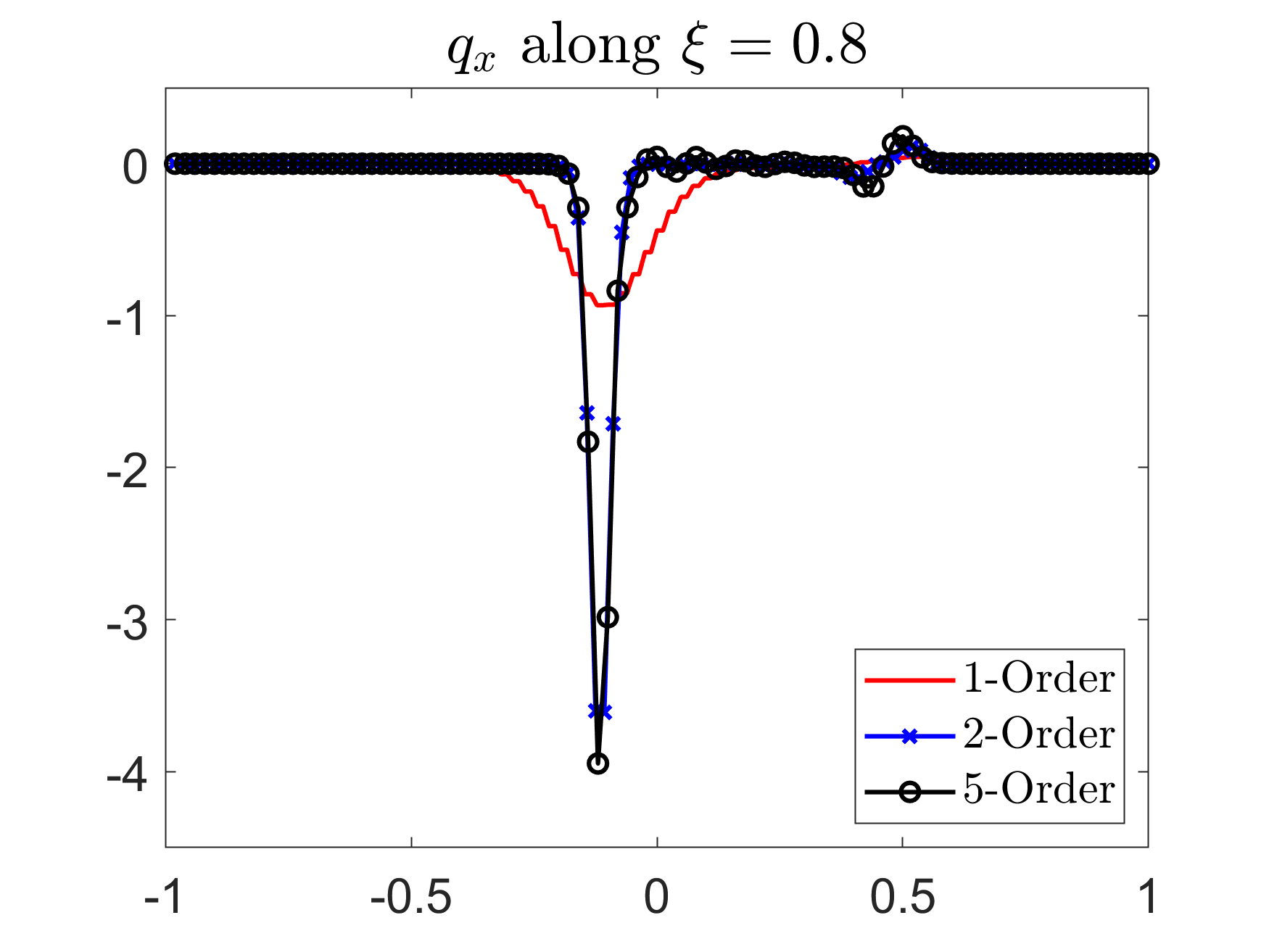}}
\caption{\sf Example 3: $\rho_x$ and $q_x$ along $\xi=0.8$ computed by the 1-Order (with $N_x=164$), 2-Order (with $N_x=112$), and 5-Order (with $N_x=100$) Young-measure schemes.}
\label{fig3.1c}
\end{figure}

In this example, we also check the efficiency of the studied 1-Order, 2-Order, and 5-Order schemes. To this end, we first measure the CPU time consumed by the 5-Order scheme, which is more expensive than the 1-Order and 2-Order schemes because of its greater complexity. We then refine the spatial meshes used by the 1-Order and 2-Order schemes so that the CPU times of all three schemes are comparable. The corresponding grids are $N_x=164$, $112$, and $100$ for the 1-Order, 2-Order, and 5-Order schemes, respectively.   
Finally, we plot the obtained $\rho_x$ and $q_x$ along the line $\xi=0.8$ in Figure \ref{fig3.1c}, where one can see that the 2-Order and 5-Order schemes still achieve much higher resolution than the 1-Order scheme. At the same time, the 5-Order scheme produces the highest resolution. This clearly demonstrates the higher efficiency of the 5-Order scheme.

\subsection{Numerical Examples for Deterministic Hyperbolic Conservation Laws}
In this section, we apply the studied Young-measure and collocation schemes to several numerical examples of deterministic hyperbolic conservation laws, including conservation laws with discontinuous fluxes, the pressureless gas dynamics system, and the Burgers equation with non-atomic support. The Young-measure formulation enables the use of different objectives in the reconstruction step. For conservation laws with discontinuous coefficients, the particular choice of entropy may affect the selected numerical approximation.

\subsubsection{Conservation Laws with Discontinuous Flux}
It is instructive to investigate how the proposed Young-measure schemes respond to different entropy objectives for conservation laws with discontinuous
coefficients.  To this end, we consider the scalar conservation law with a discontinuous flux function:
$$
u_t+F(u)_x=0,
$$
where $F(u)=(1-H(x))g(u)+H(x)f(u)$, and $H(x)$ is the Heaviside function. We refer the readers to \cite{Mishra2017} for more details.

We emphasize that, for scalar conservation laws with discontinuous fluxes, the admissible solution is determined not only by the conservation law away from the interface, but also by the interface entropy condition, or equivalently by the selection criterion imposed at the discontinuity of the flux. Therefore, uniqueness should be understood only after such an interface condition has been specified. In the present Young-measure formulation, the entropy used in the objective function of the linear-programming problem acts as a selection criterion in the reconstruction step. Consequently, changing the objective function may lead to different selected numerical approximations. This should not be interpreted as producing two different solutions for the same completely specified interface entropy condition, but rather as illustrating how different entropy-based selection criteria can be incorporated into the Young-measure reconstruction.

\paragraph*{Example 4.} In this example taken from \cite{Mishra2017}, we take 
$$
g(u)=u(1-u), \quad  f(u)=1.1u(1-u),
$$ 
and consider the following initial data:
\begin{equation*}
 u(x,\xi,0) =
\begin{cases}
  0.65, & \mbox{if } x<0 \\
  0.35, & \mbox{otherwise},
\end{cases}
\end{equation*}
prescribed in the computational domain $[-4,4]$ subject to the free boundary conditions. The corresponding entropy functions $\eta$ in \eref{linprog-a} are $\eta(u)=\frac{1}{2}u^2$ and $\eta(u)=|u-c|$, where $c$ is a constant satisfying $c\in [u_{\min}, u_{\max}]$.

\begin{figure}[ht!]
\centerline{\includegraphics[trim=0.2cm 0.2cm 0.8cm 0.1cm, clip, width=5.cm]{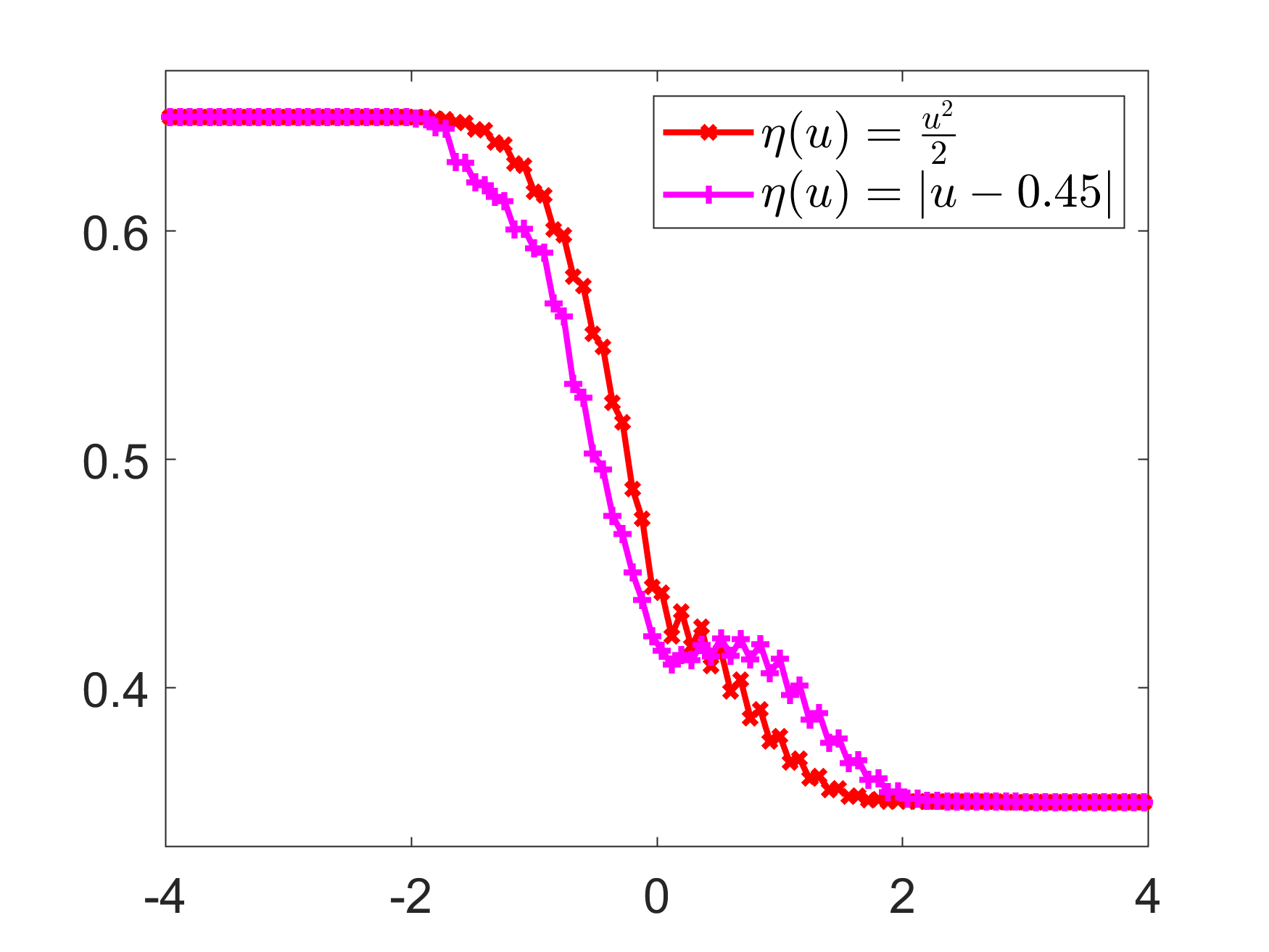}\hspace{0.5cm}
            \includegraphics[trim=0.2cm 0.2cm 0.8cm 0.1cm, clip, width=5.cm]{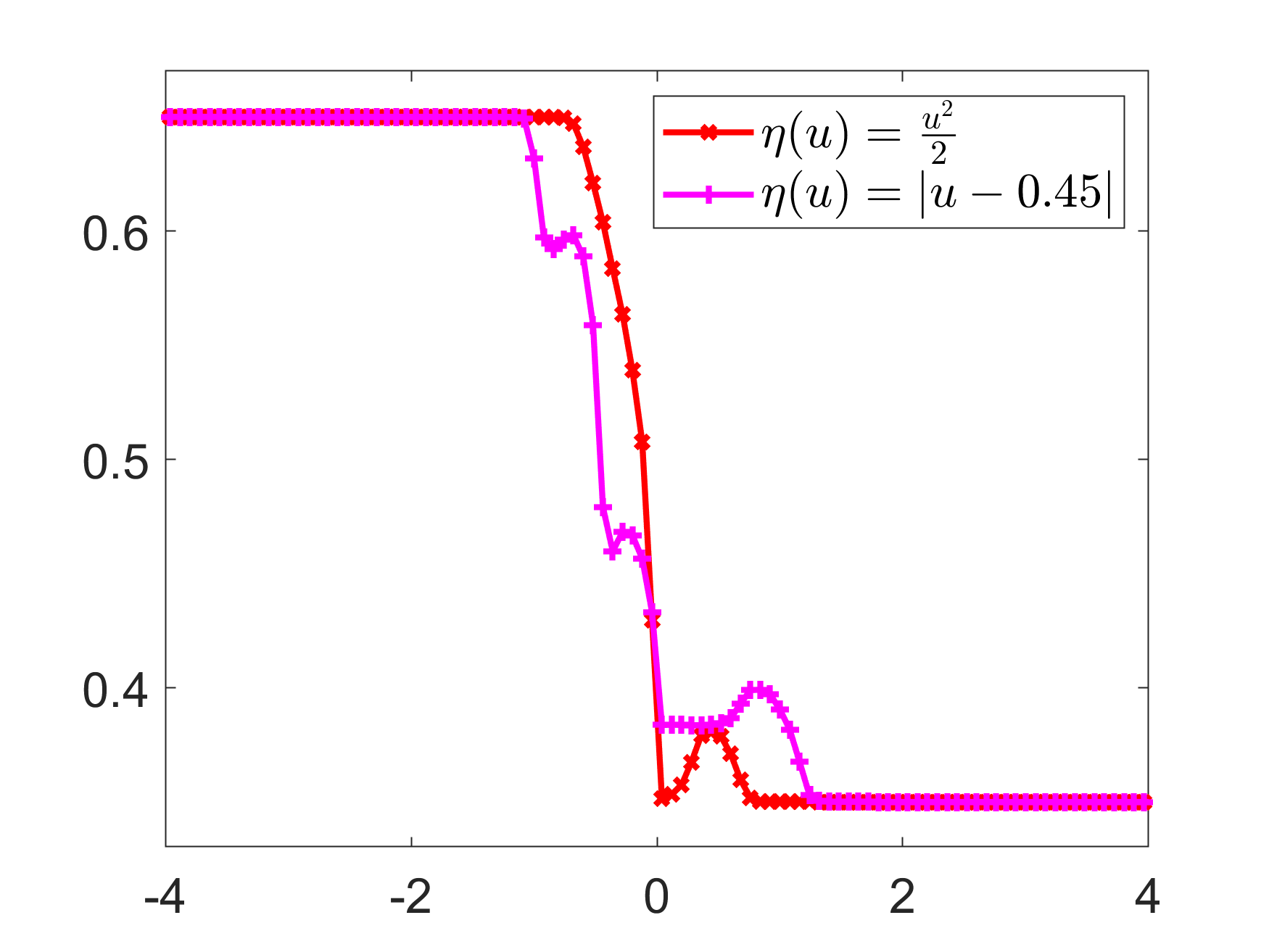}\hspace{0.5cm}
            \includegraphics[trim=0.2cm 0.2cm 0.8cm 0.1cm, clip, width=5.cm]{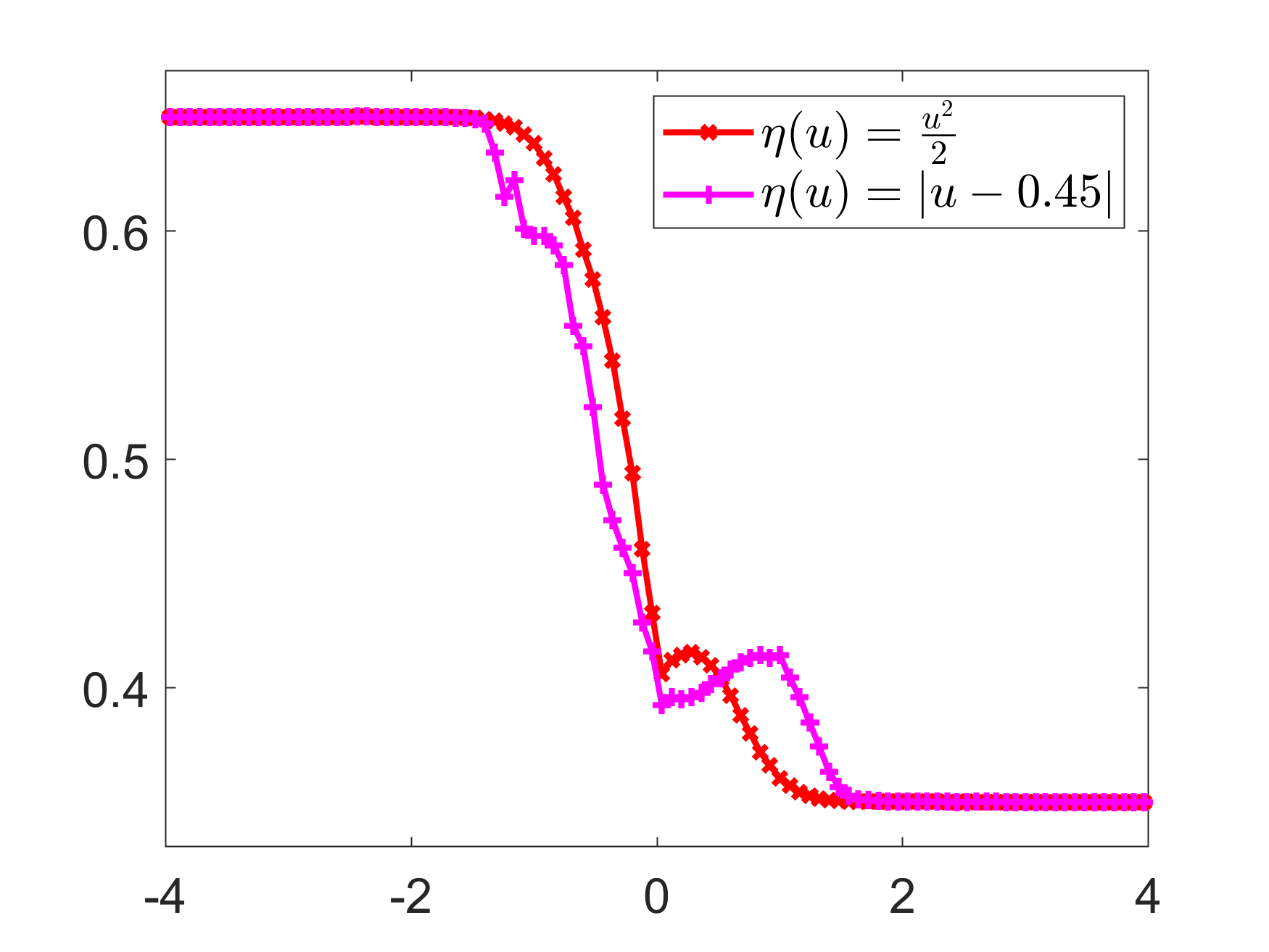}}
\caption{\sf Example 4: Numerical results computed by the 1-Order (left), 2-Order (middle), and 5-Order (right) collocation and Young-measure methods with $\eta(u)=\frac{1}{2}u^2$ and $\eta(u)=|u-0.45|$. \label{fig3.10a}}
\end{figure}

\begin{figure}[ht!]
\centerline{\includegraphics[trim=0.2cm 0.7cm 0.8cm 0.1cm, clip, width=5.cm]{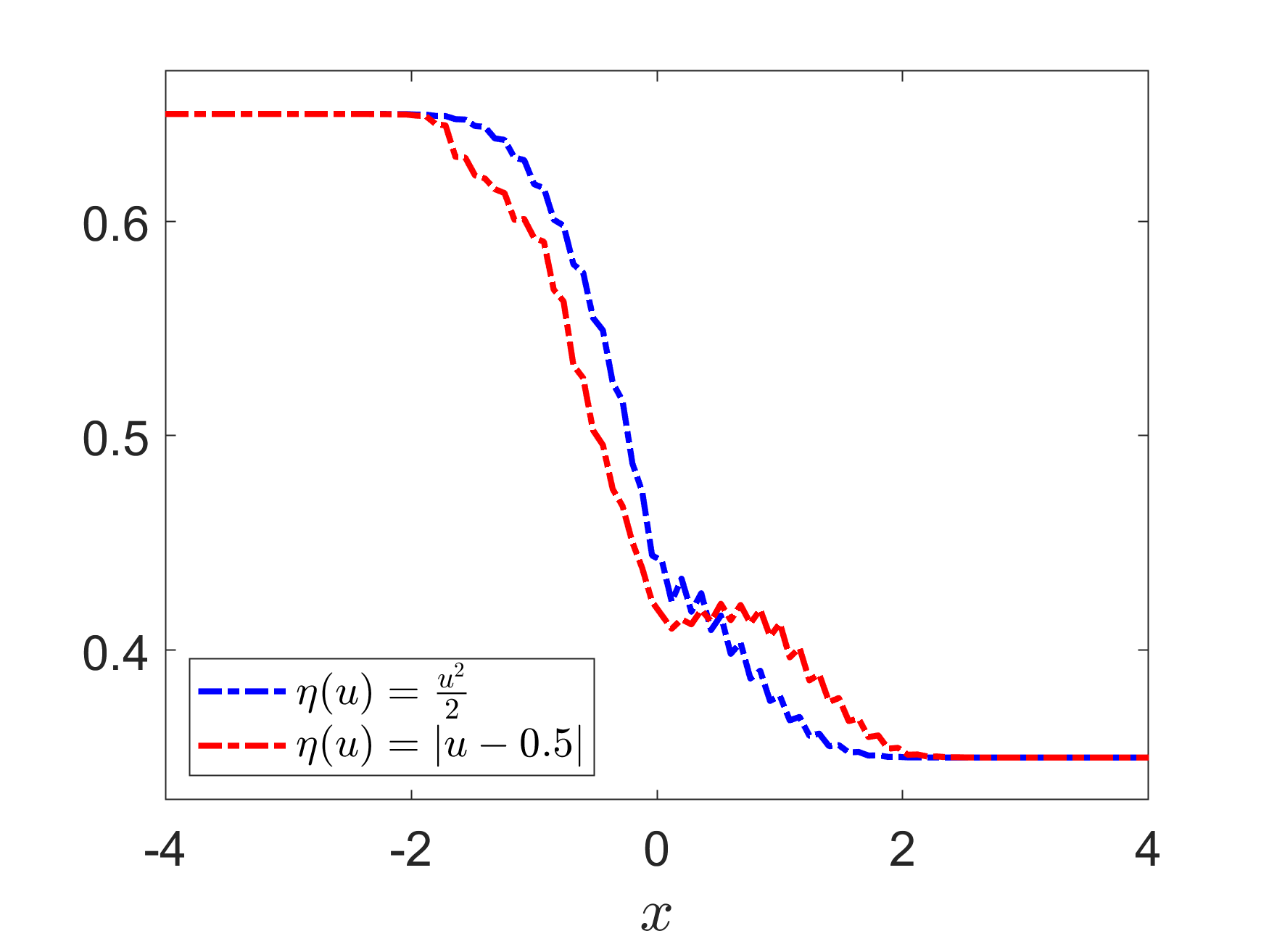}\hspace{0.5cm}
            \includegraphics[trim=0.2cm 0.7cm 0.8cm 0.1cm, clip, width=5.cm]{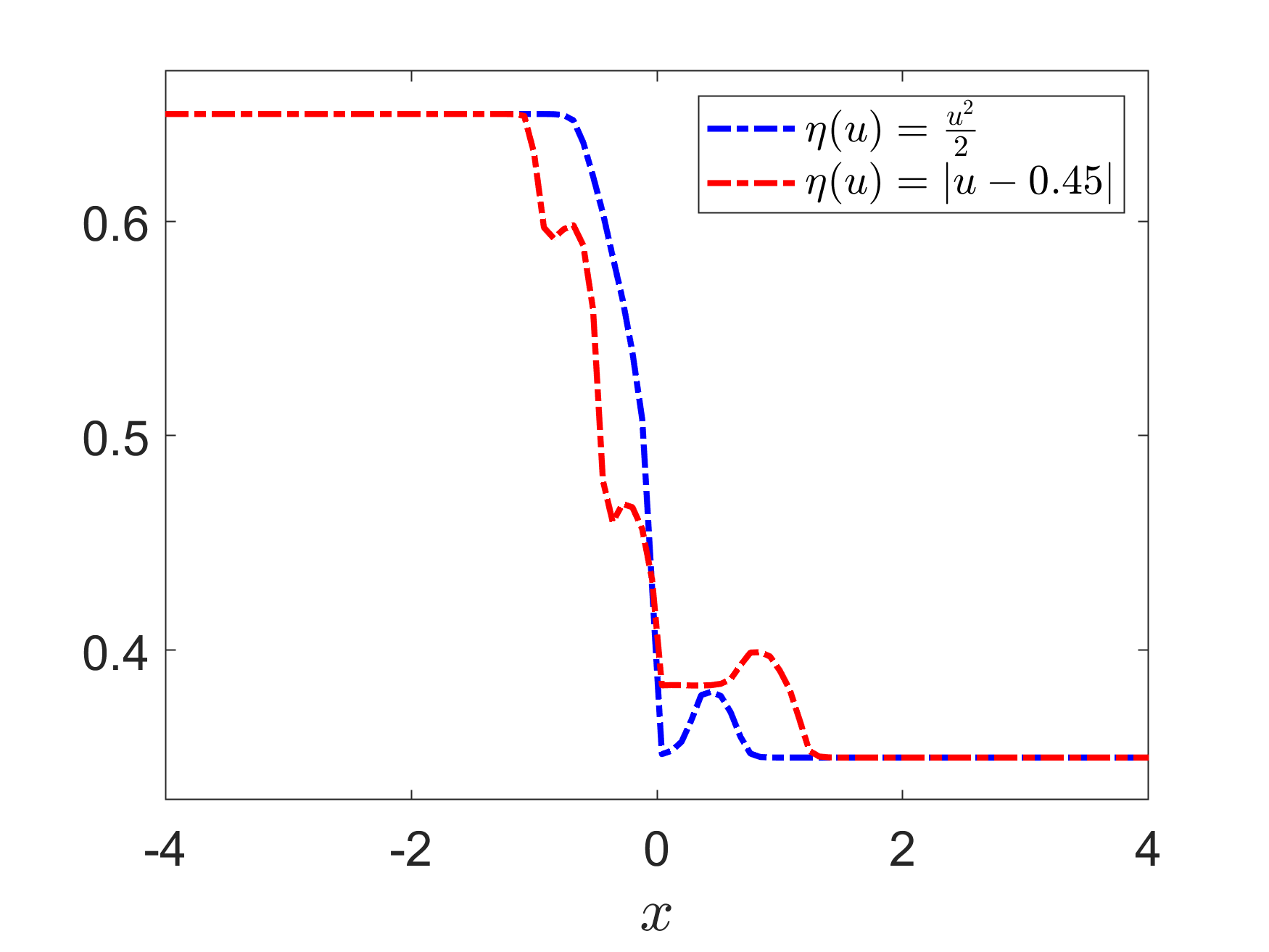}\hspace{0.5cm}
            \includegraphics[trim=0.2cm 0.7cm 0.8cm 0.1cm, clip, width=5.cm]{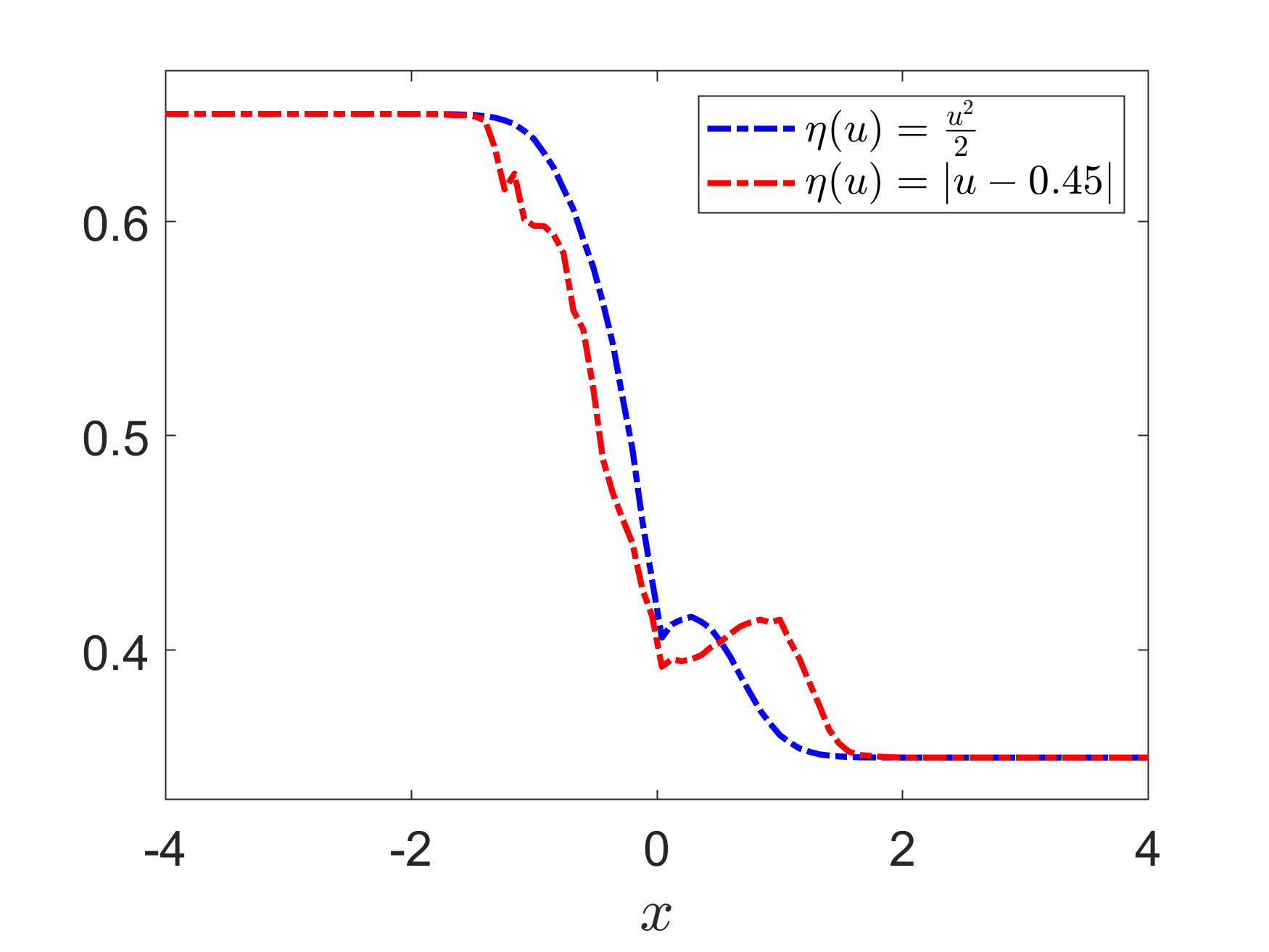}}
\caption{\sf Example 4: Moments computed by the 1-Order (left), 2-Order (middle), and 5-Order (right) Young-measure methods at the final time $t=2$ with $\eta(u)=\frac{1}{2}u^2$ and $\eta(u)=|u-0.45|$. \label{fig3.10b}}
\end{figure}

\begin{figure}[ht!]
\centerline{\includegraphics[trim=0.2cm 0.2cm 0.8cm 0.1cm, clip, width=5.cm]{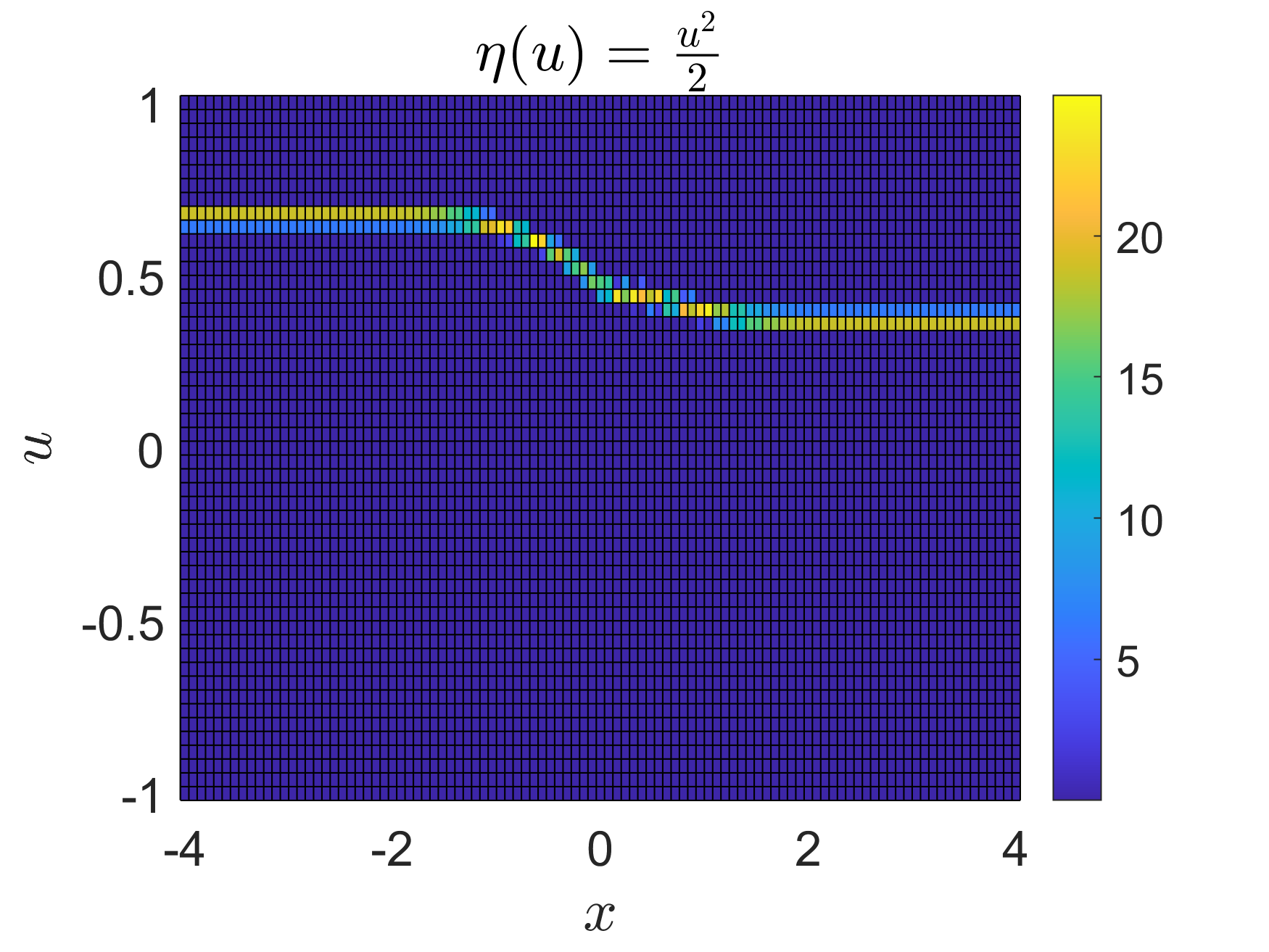}\hspace{0.5cm}
            \includegraphics[trim=0.2cm 0.2cm 0.8cm 0.1cm, clip, width=5.cm]{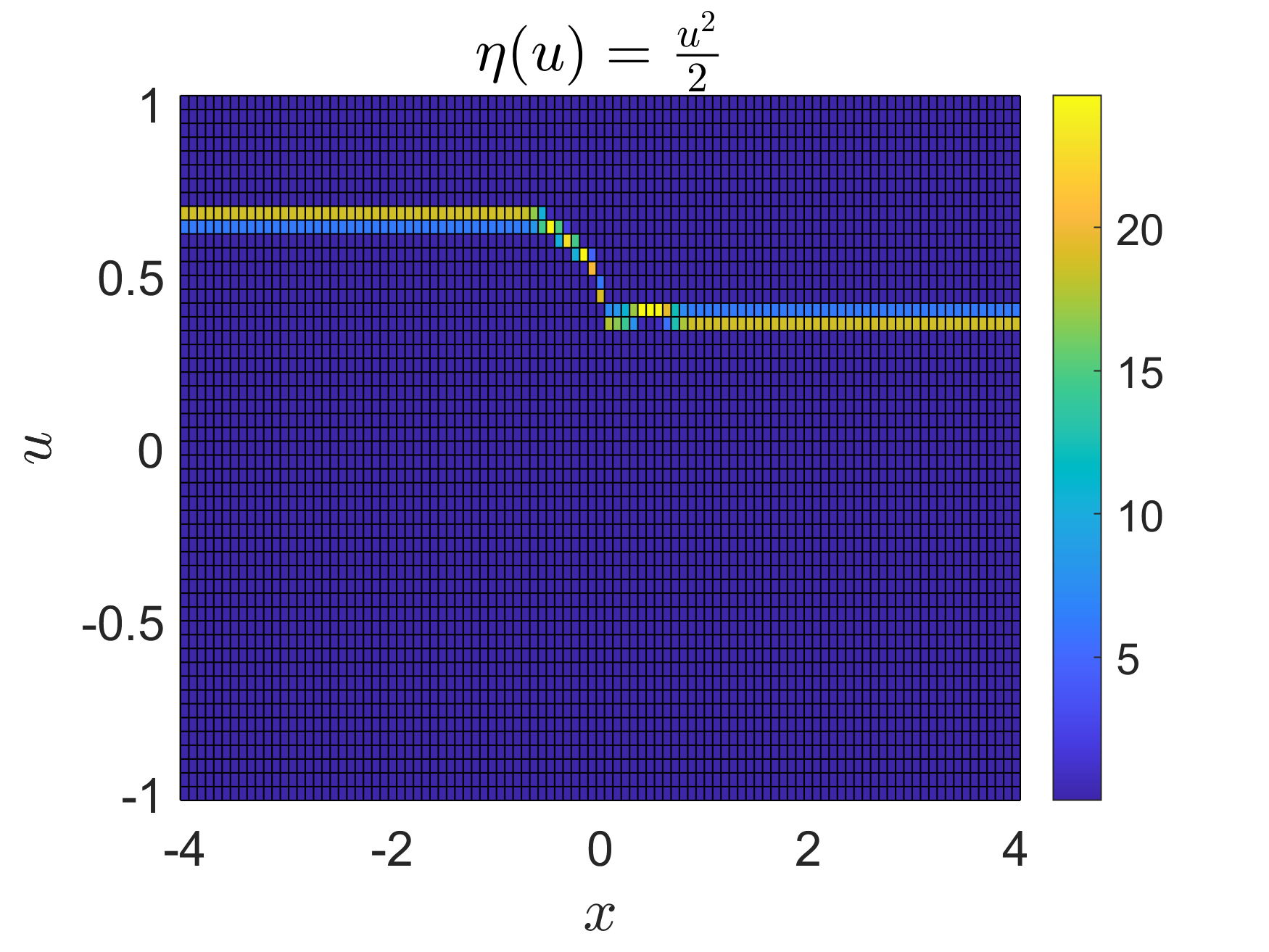}\hspace{0.5cm}
            \includegraphics[trim=0.2cm 0.2cm 0.8cm 0.1cm, clip, width=5.cm]{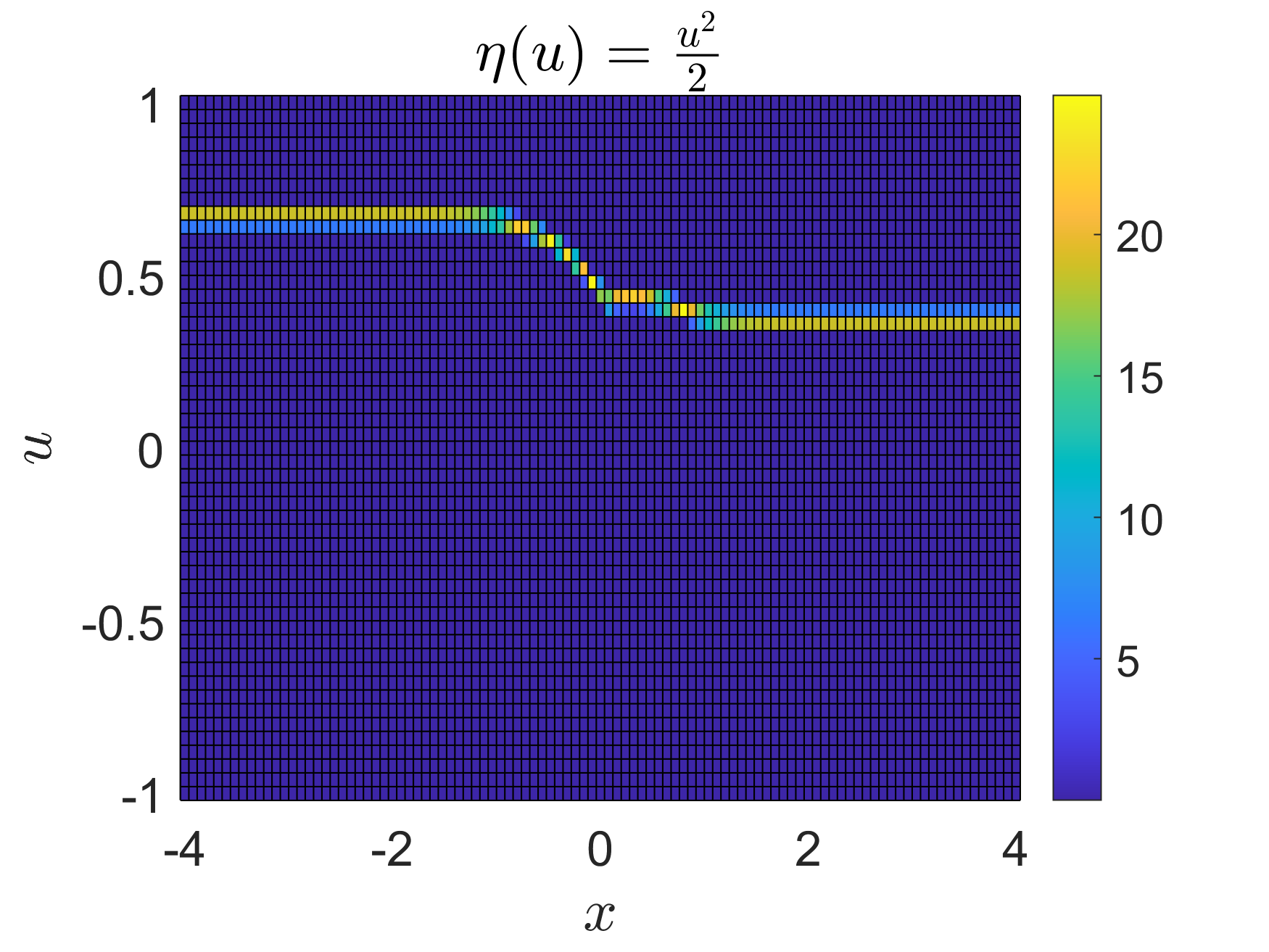}}
\vskip 15pt
\centerline{\includegraphics[trim=0.2cm 0.2cm 0.8cm 0.1cm, clip, width=5.cm]{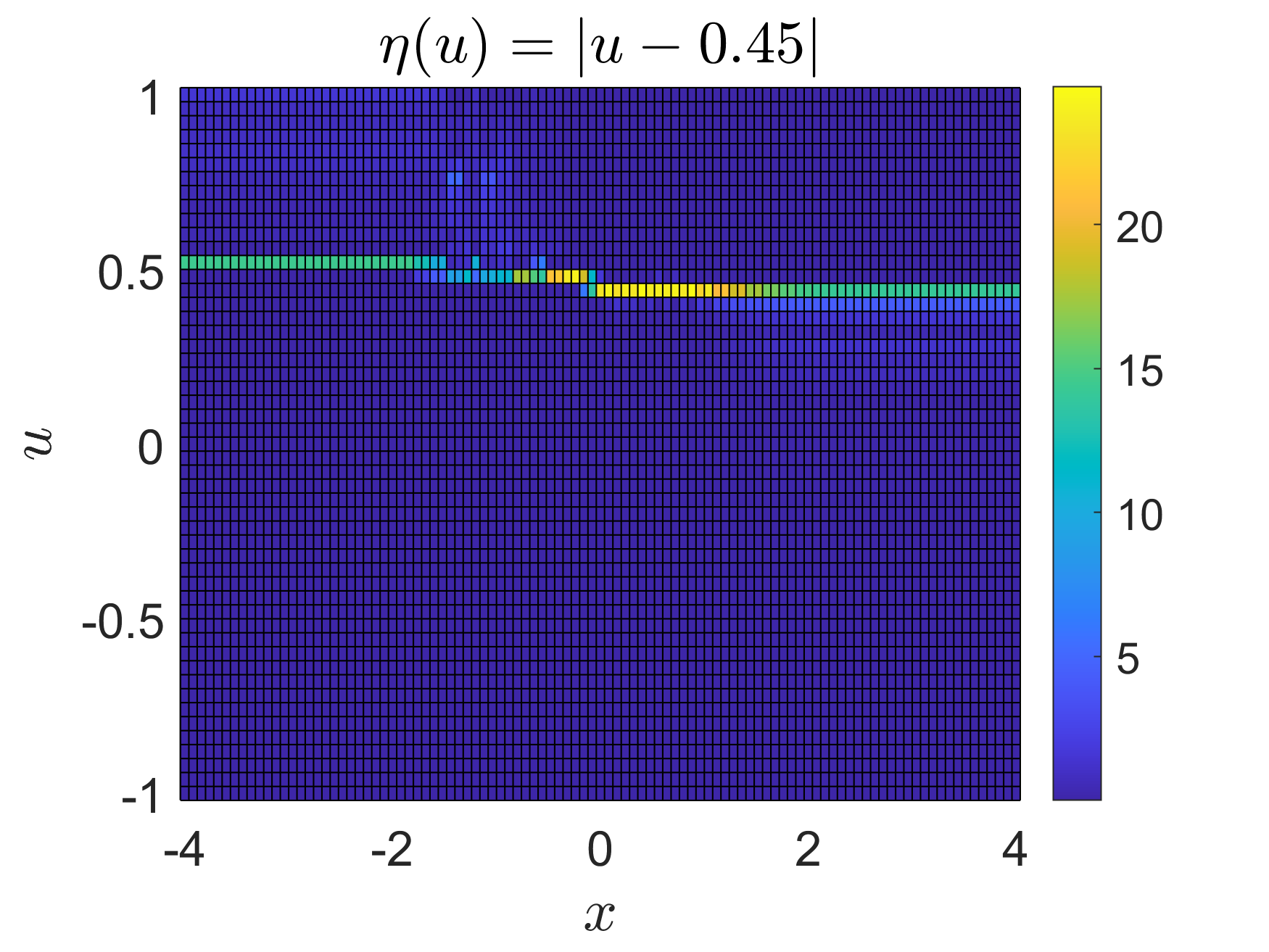}\hspace{0.5cm}
            \includegraphics[trim=0.2cm 0.2cm 0.8cm 0.1cm, clip, width=5.cm]{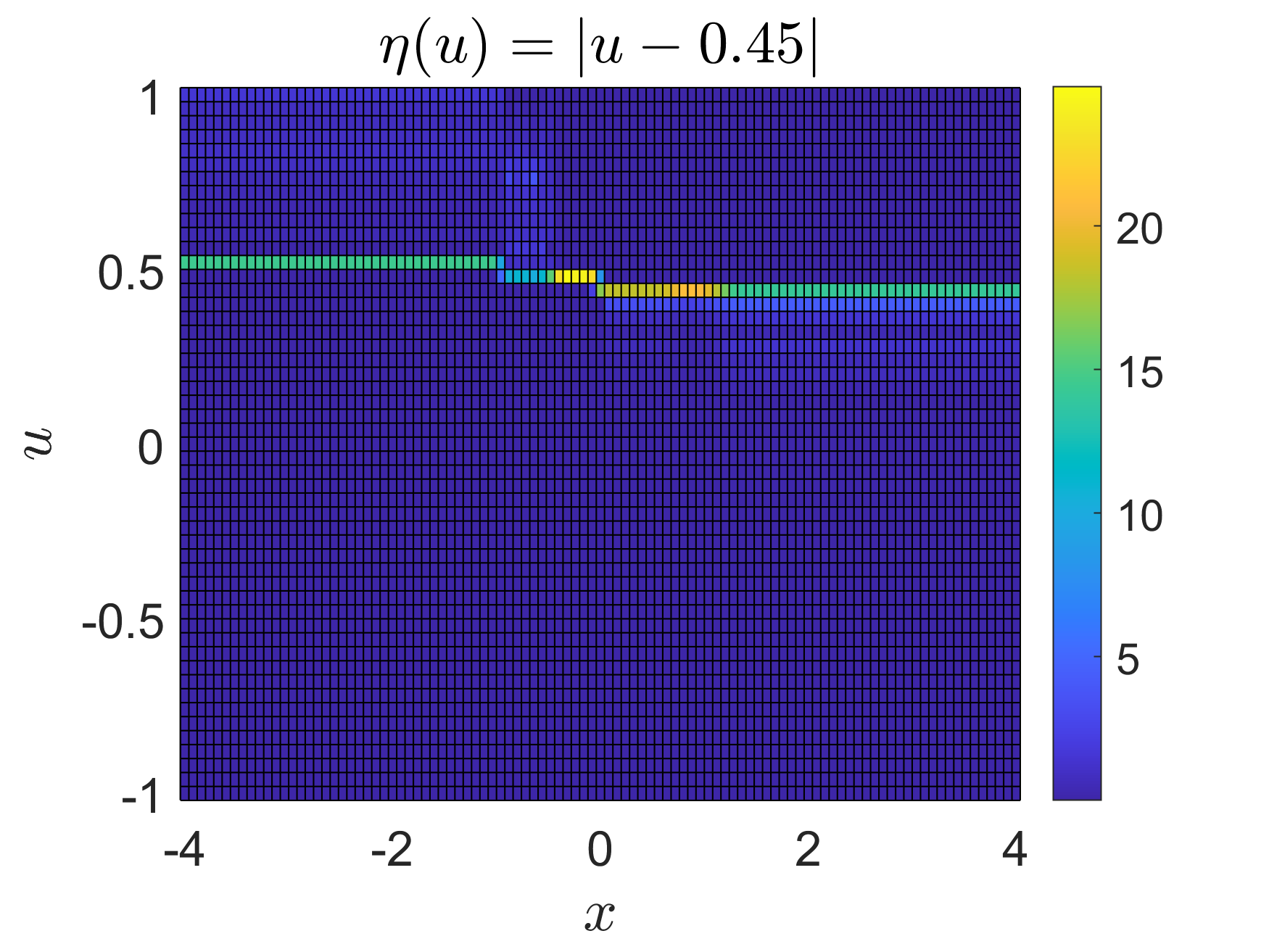}\hspace{0.5cm}
            \includegraphics[trim=0.2cm 0.2cm 0.8cm 0.1cm, clip, width=5.cm]{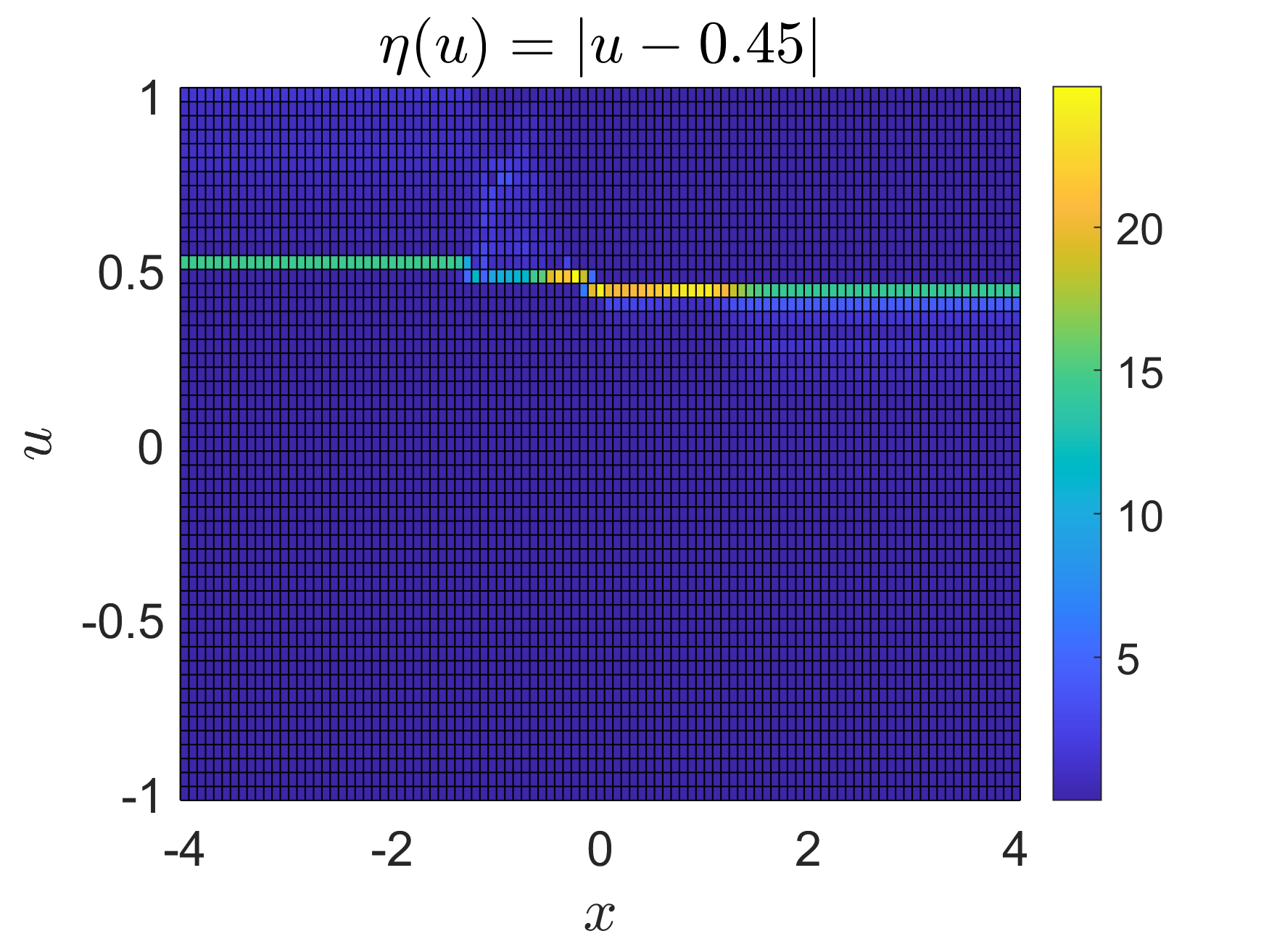}}
\caption{\sf Example 4: Supports of the Young measures computed by the 1-Order (left column), 2-Order (middle column), and 5-Order (right column) Young-measure method at the final time $t=2$ with $\eta(u)=\frac{1}{2}u^2$ (top row) and $\eta(u)=|u-0.45|$ (bottom row). \label{fig3.10c}}
\end{figure}

We compute the numerical solutions until the final time $t=2$ by the studied 1-Order, 2-Order, and 5-Order schemes on the uniform mesh with $N_x=100$, $N_\xi=1$. For the discretization in phase space, we set $u\in[-1,1]$ and use $N_{u}=50$. We plot the obtained numerical results in Figure \ref{fig3.10a}, where one can see that the numerical results computed by the studied Young-measure schemes are different when utilizing different entropy functions. The obtained results are in a good agreement with the numerical results reported in \cite{Mishra2017}, where various schemes are used to produce different numerical solutions.  We also present the corresponding moments and supports of the Young measures  in Figures \ref{fig3.10b}--\ref{fig3.10c}. One can notice that the moments and supports of the Young measures also vary when different entropy functions are utilized.

\medskip 

\paragraph*{Example 5.} In this example taken from \cite{Mishra2017}, we take 
$$
g(u)=\frac{2u(1-u)}{1+u}, \quad  f(u)=\frac{2u(1-u)}{2-u},
$$
and consider the following initial data:
\begin{equation*}
 u(x,\xi,0) =0.5,
\end{equation*}
prescribed in the computational domain $[-4,4]$ subject to the free boundary conditions. The corresponding entropy functions $\eta$ in \eref{linprog-a} are $\eta(u)=\frac{1}{2}u^2$ and $\eta(u)=|u-0.5|$.

We compute the numerical solutions until the final time $t=3$ by the studied 1-Order, 2-Order, and 5-Order schemes on the uniform mesh with $N_x=100$, $N_\xi=1$. For the discretization in phase space, we set $u\in[-1,1]$ and use $N_{u}=50$. The obtained numerical results are presented in Figures 5.14--5.16. In this example, the two entropy objectives lead to essentially the same first moment, while the support plots illustrate the corresponding reconstructed Young measures. Together with Example 4, this demonstrates how the Young-measure formulation can incorporate different entropy-based selection mechanisms through the objective function in the linear-programming problem.

\begin{figure}[ht!]
\centerline{\includegraphics[trim=0.2cm 0.2cm 0.8cm 0.1cm, clip, width=5.cm]{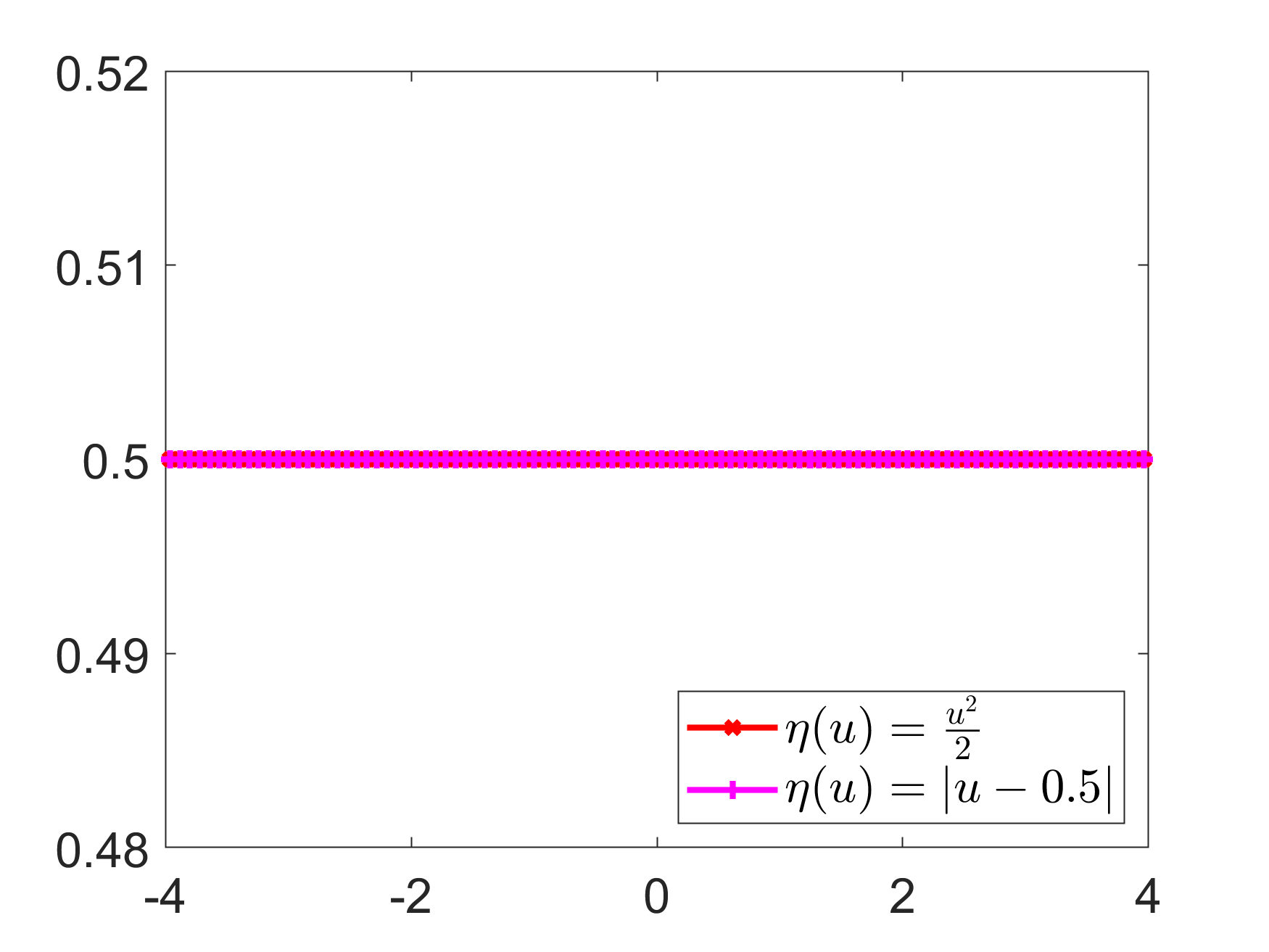}\hspace{0.5cm}
            \includegraphics[trim=0.2cm 0.2cm 0.8cm 0.1cm, clip, width=5.cm]{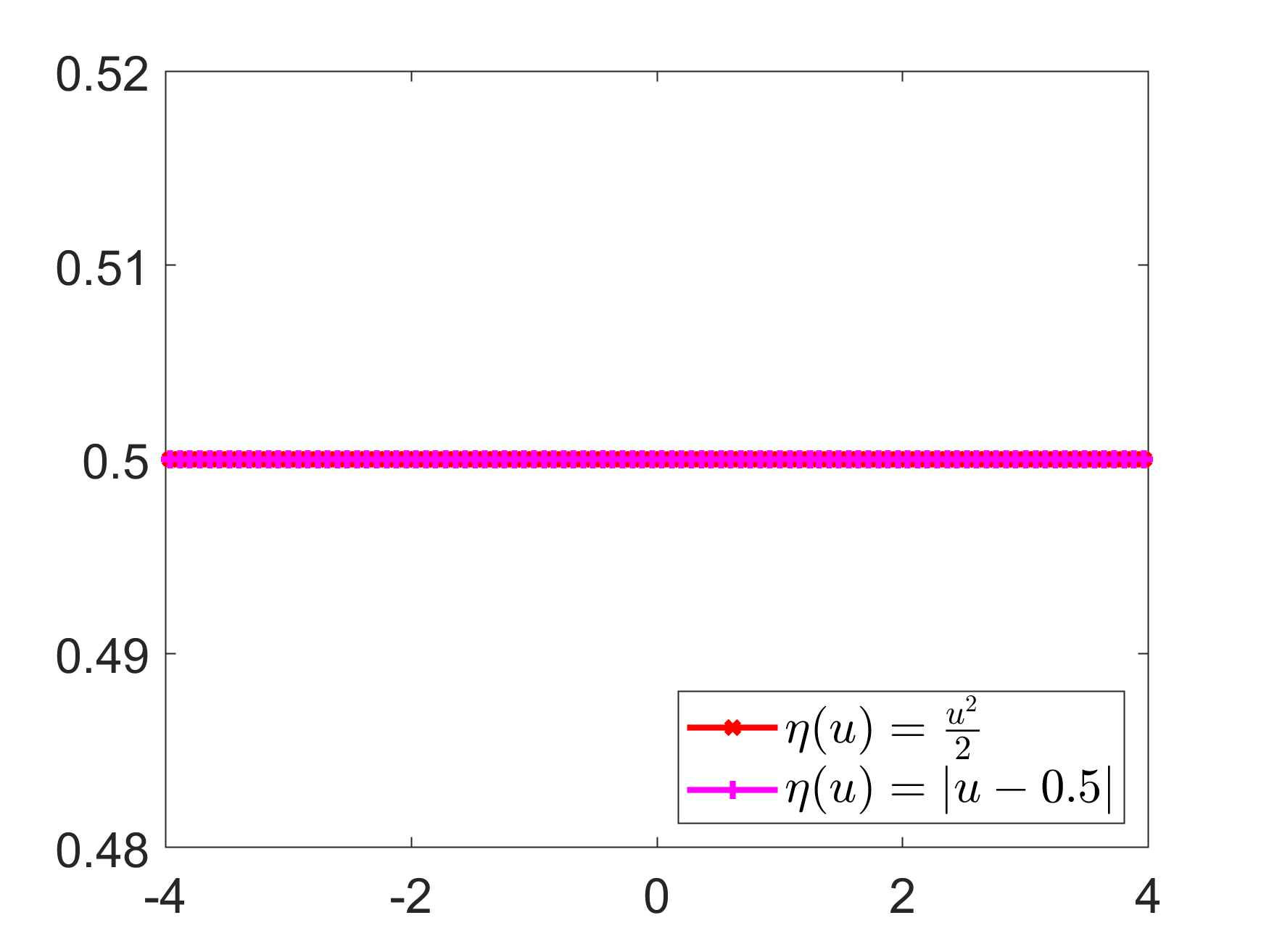}\hspace{0.5cm}
            \includegraphics[trim=0.2cm 0.2cm 0.8cm 0.1cm, clip, width=5.cm]{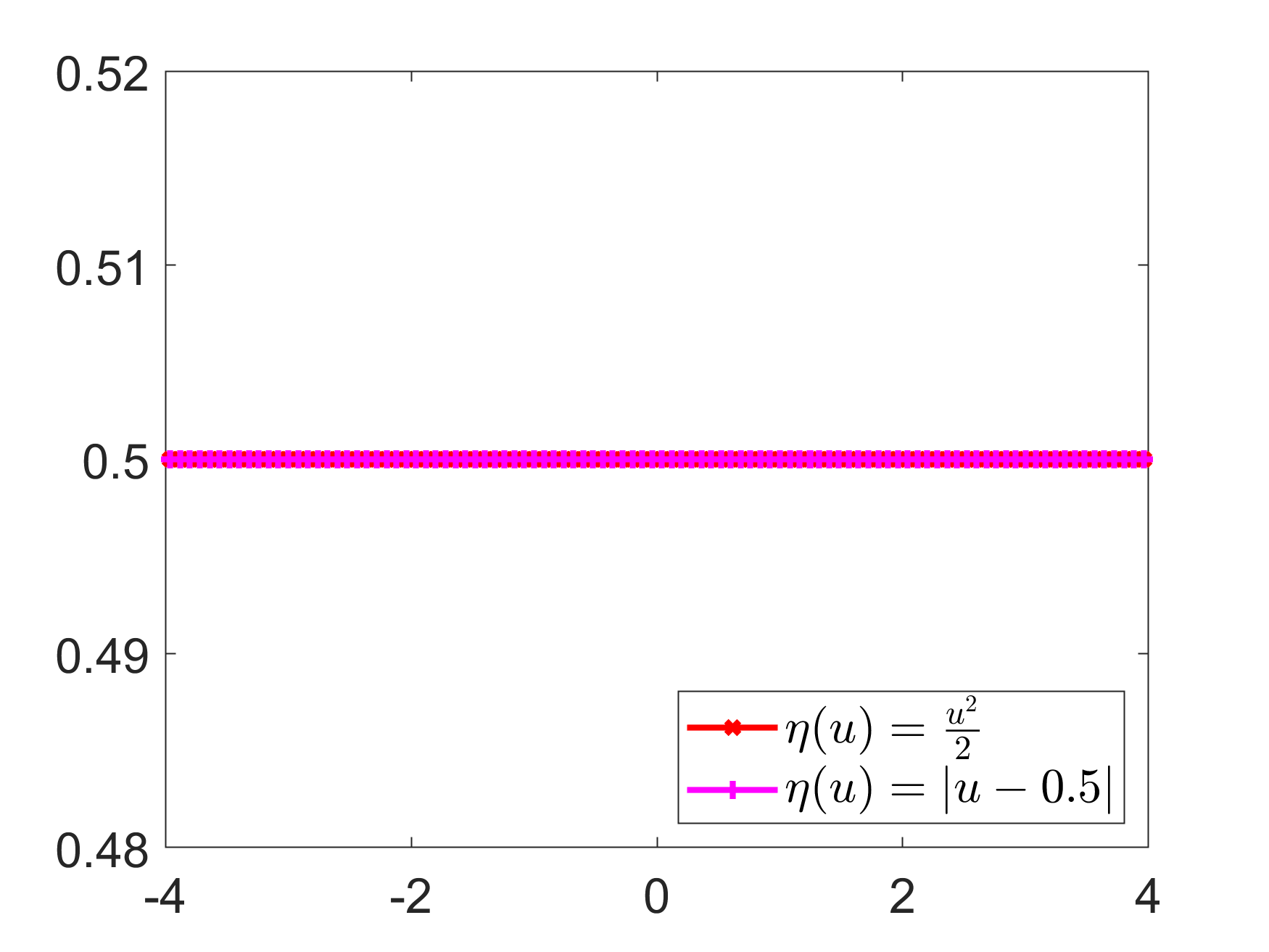}}
\caption{\sf Example 5: Numerical results computed by the 1-Order (left), 2-Order (middle), and 5-Order (right) Young-measure methods with $\eta(u)=\frac{1}{2}u^2$ and $\eta(u)=|u-0.5|$. \label{fig3.11a}}
\end{figure}

\begin{figure}[ht!]
\centerline{\includegraphics[trim=0.2cm 0.7cm 0.8cm 0.1cm, clip, width=5.cm]{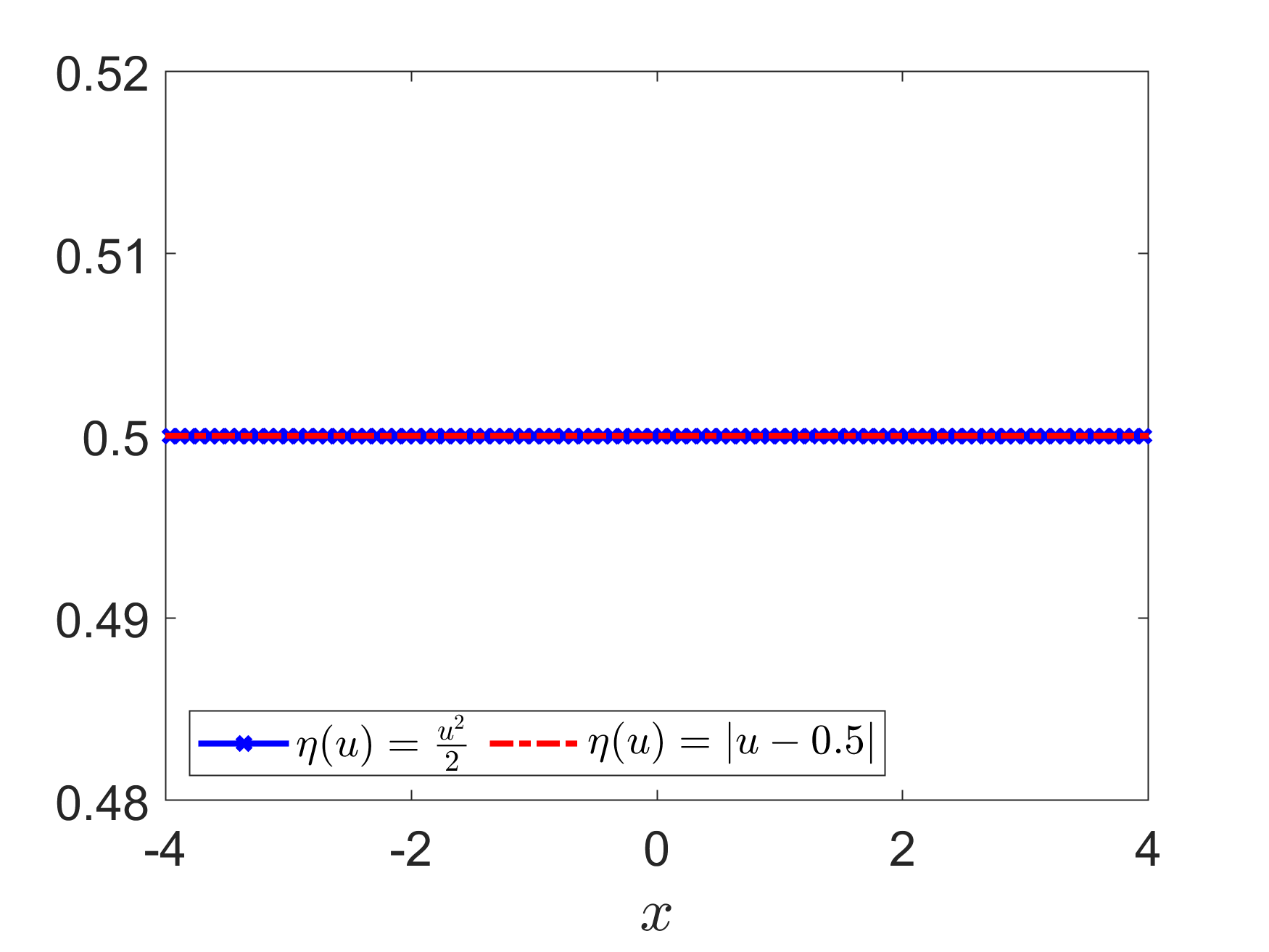}\hspace{0.5cm}
            \includegraphics[trim=0.2cm 0.7cm 0.8cm 0.1cm, clip, width=5.cm]{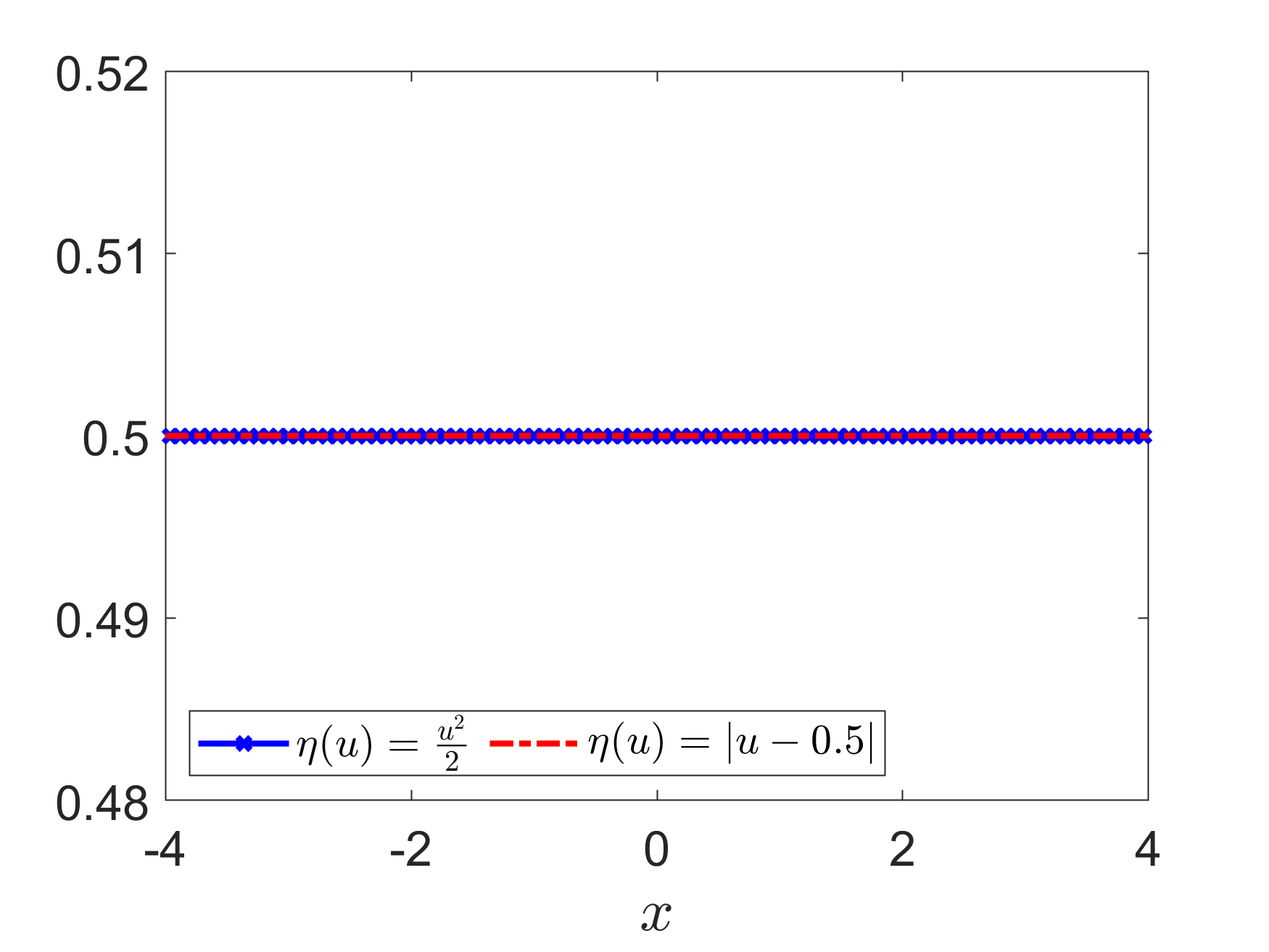}\hspace{0.5cm}
            \includegraphics[trim=0.2cm 0.7cm 0.8cm 0.1cm, clip, width=5.cm]{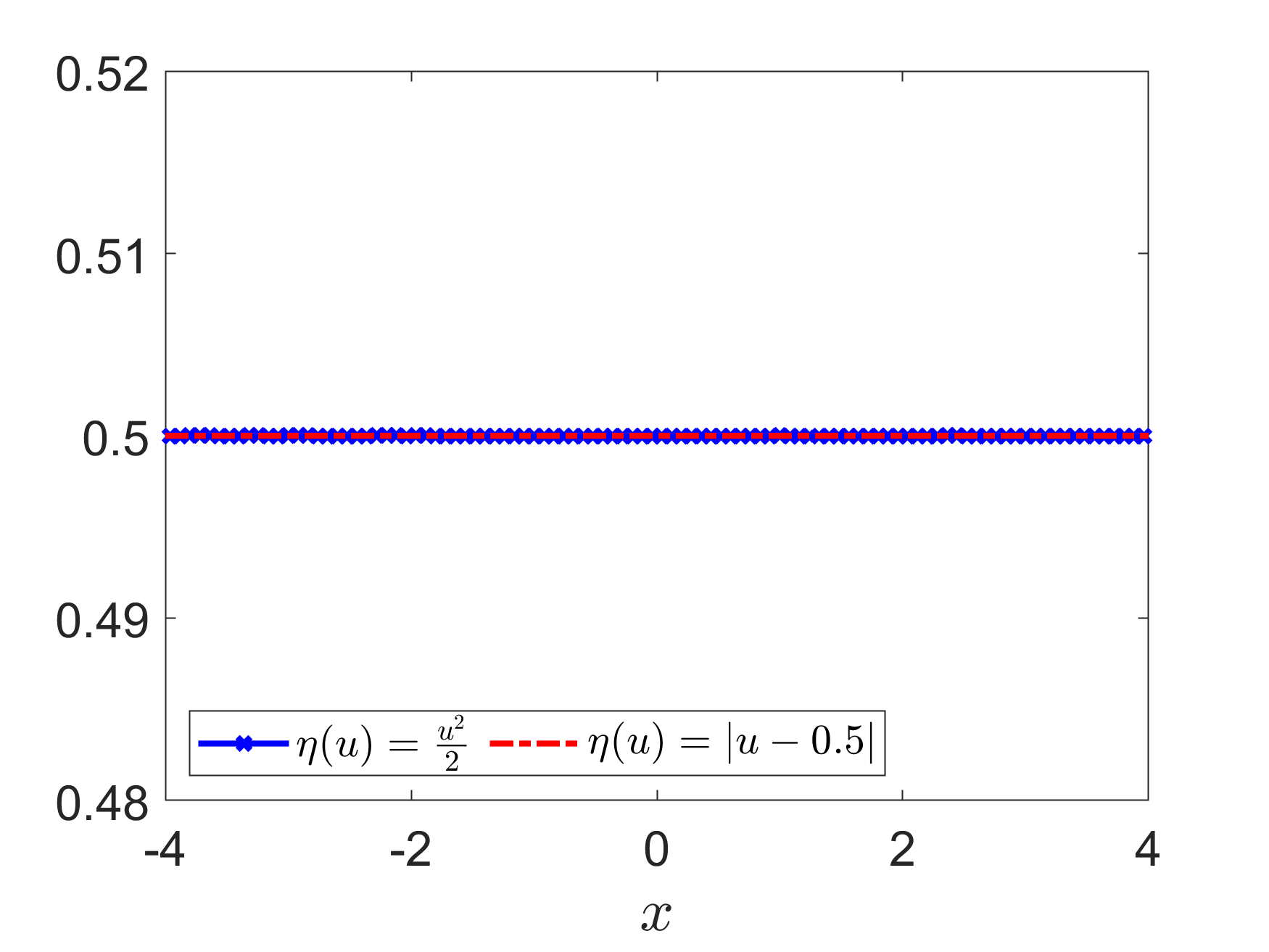}}
\caption{\sf Example 5: Moments computed by the 1-Order (left), 2-Order (middle), and 5-Order (right) Young-measure methods at the final time $t=3$ with $\eta(u)=\frac{1}{2}u^2$ and $\eta(u)=|u-0.5|$. \label{fig3.11b}}
\end{figure}

\begin{figure}[ht!]
\centerline{\includegraphics[trim=0.2cm 0.2cm 0.8cm 0.1cm, clip, width=5.cm]{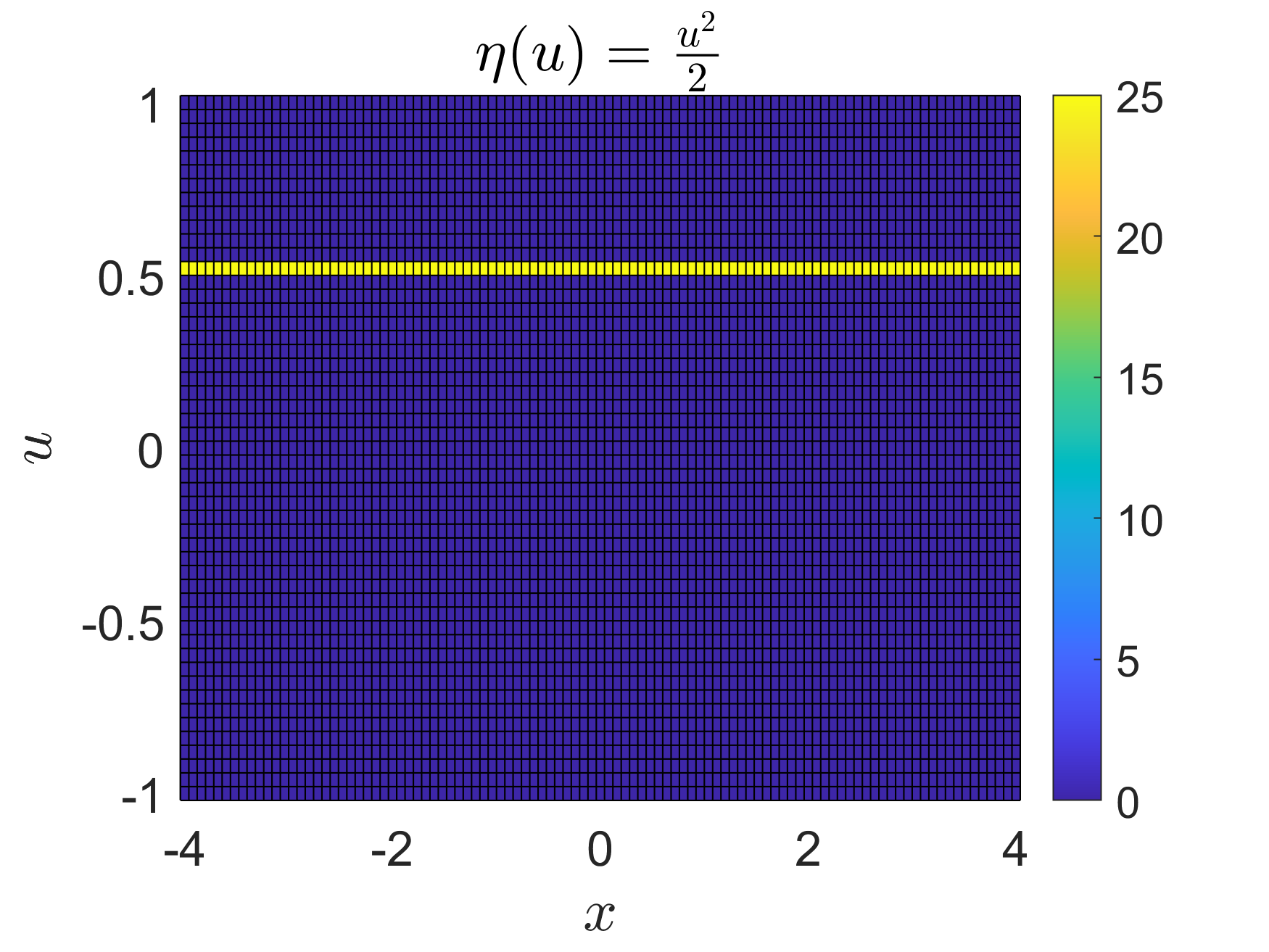}\hspace{0.5cm}
            \includegraphics[trim=0.2cm 0.2cm 0.8cm 0.1cm, clip, width=5.cm]{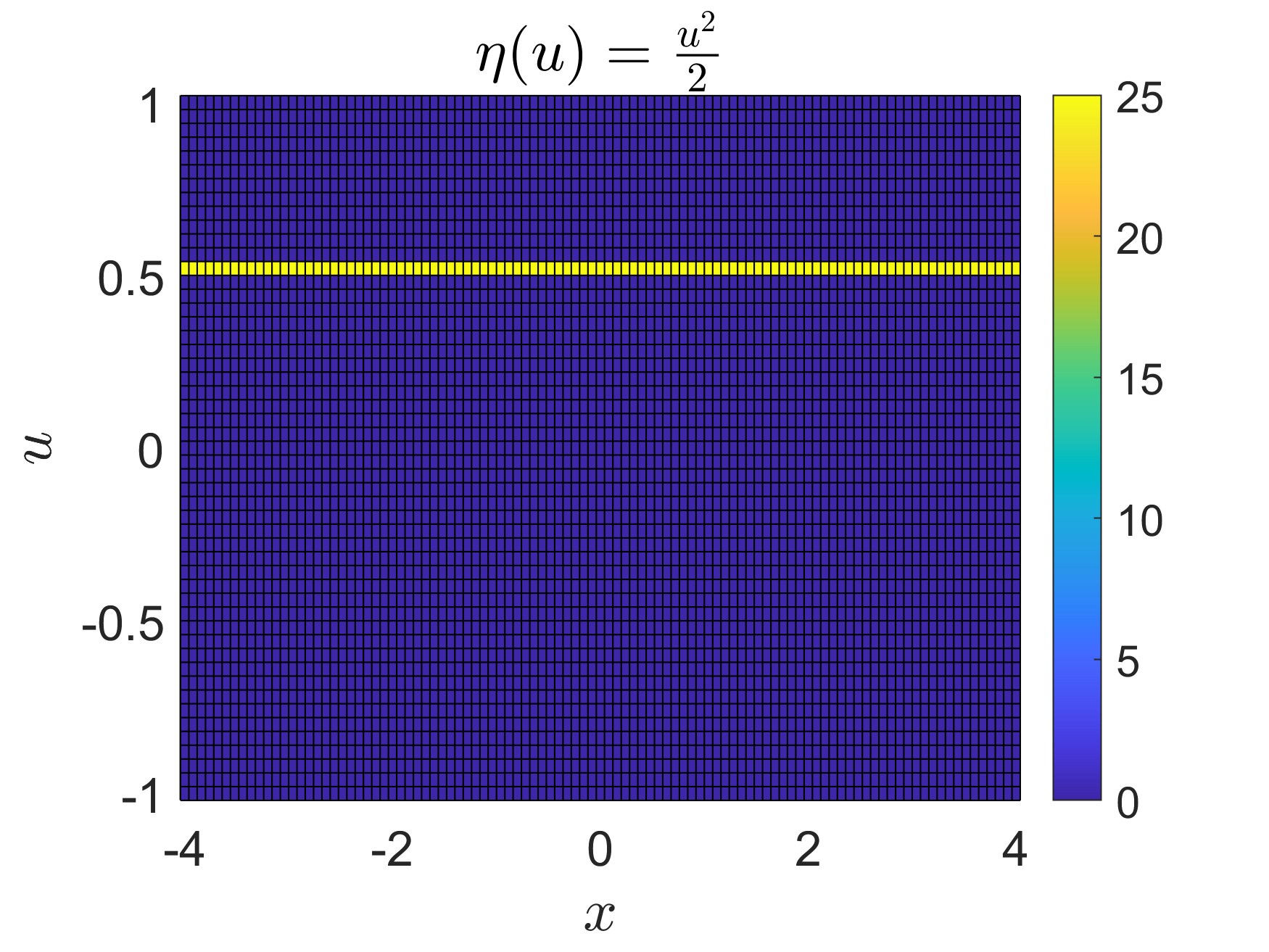}\hspace{0.5cm}
            \includegraphics[trim=0.2cm 0.2cm 0.8cm 0.1cm, clip, width=5.cm]{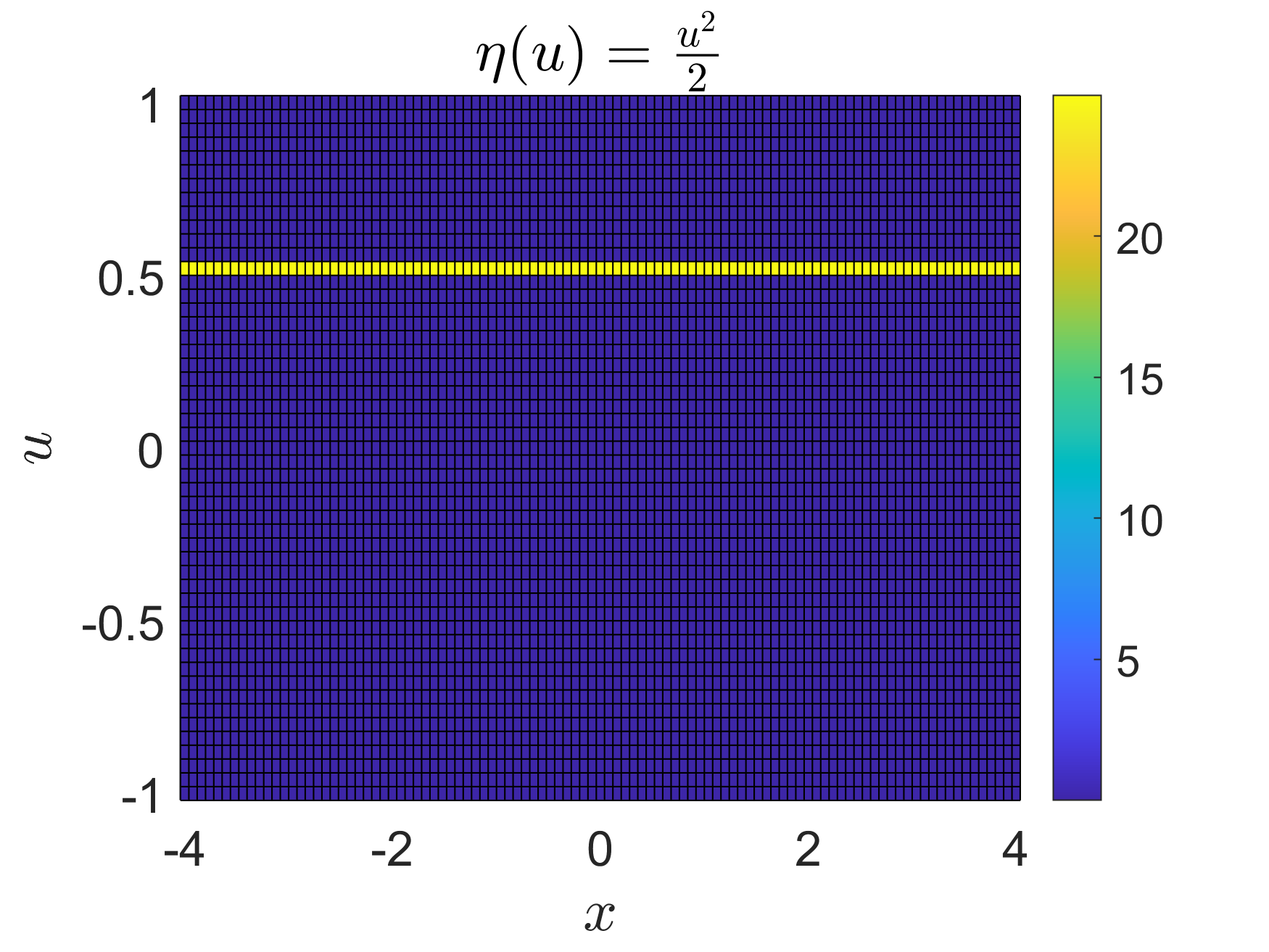}}
\vskip 15pt
\centerline{\includegraphics[trim=0.2cm 0.2cm 0.8cm 0.1cm, clip, width=5.cm]{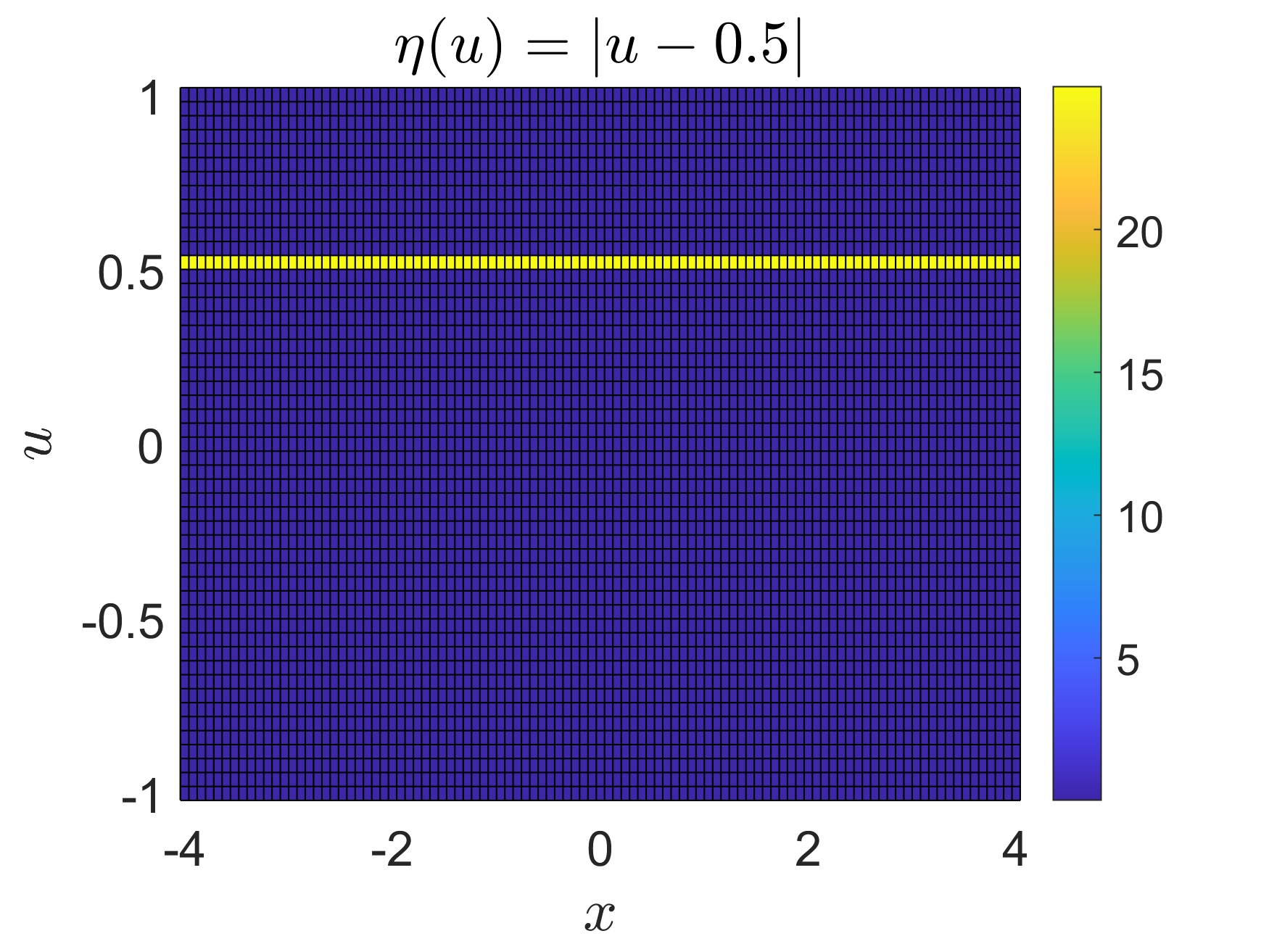}\hspace{0.5cm}
            \includegraphics[trim=0.2cm 0.2cm 0.8cm 0.1cm, clip, width=5.cm]{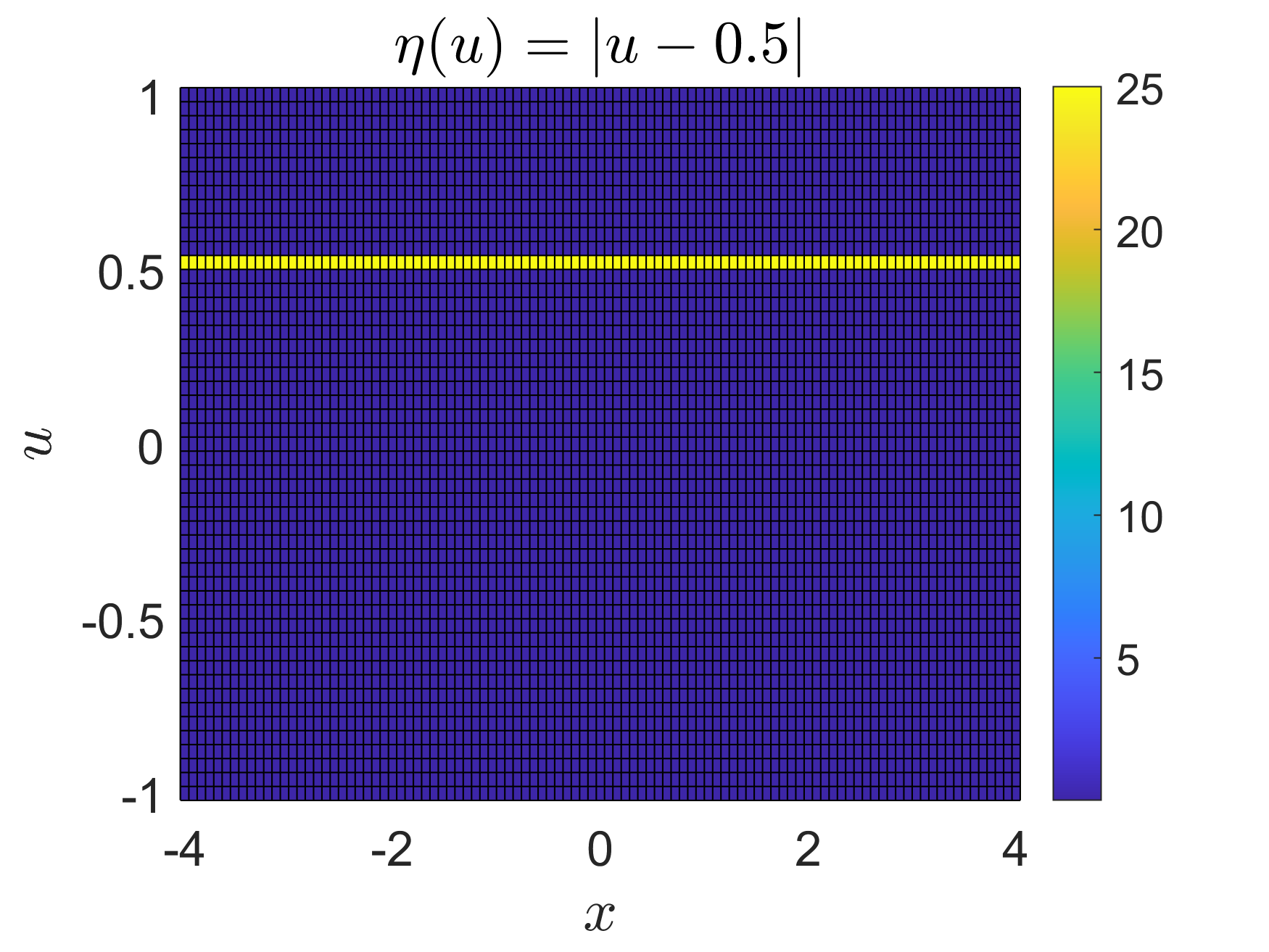}\hspace{0.5cm}
            \includegraphics[trim=0.2cm 0.2cm 0.8cm 0.1cm, clip, width=5.cm]{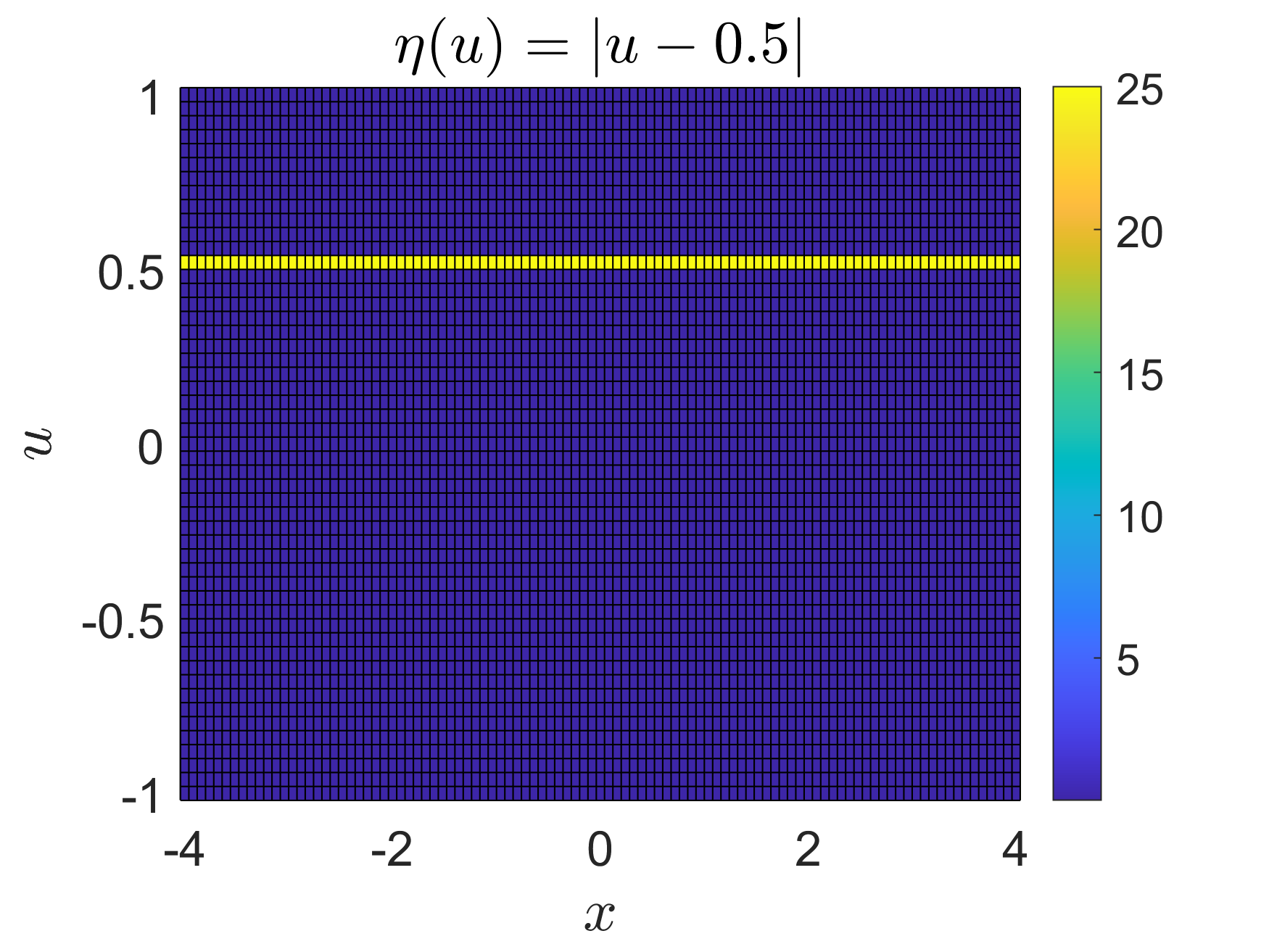}}
\caption{\sf Example 5: Supports obtained by the 1-Order (left column), 2-Order (middle column), and 5-Order (right column) Young-measure methods at the final time $t=3$ with $\eta(u)=\frac{1}{2}u^2$ (top row) and $\eta(u)=|u-0.5|$ (bottom row). \label{fig3.11c}}
\end{figure}

\begin{rmk}
Examples 4 and 5 show that in the case of discontinuous fluxes, the entropy used in the linear-programming objective plays the role of a selection criterion in the reconstruction step. The numerical results, especially those in Example 4, reflect the effect of different entropy-based selection mechanisms, while Example 5 illustrates that different objectives may also lead to essentially the same first moment. Moreover, by the moment constraint in the linear-programming problem, the first moment of the reconstructed Young measure satisfies
$$
\Delta V\sum_{\ell=1}^{N_u} z_\ell \mu^*_{i,j,\ell}=u^n_{i,j},
$$
up to the tolerance of the linear-programming solver and the phase-space discretization error. Hence, the moment plots in Figures \ref{fig3.10b} and \ref{fig3.11b} represent the evolved conservative variable, while the support plots in Figures \ref{fig3.10c} and \ref{fig3.11c} show how the corresponding reconstructed Young measure is distributed in phase space.

In a deterministic setting, if the atomic measure is admissible and a strictly convex entropy is used, the entropy-minimizing Young measure is expected to be concentrated near the deterministic state, up to the phase-space discretization. For non-strictly convex objectives such as \(\eta(u)=|u-c|\), the minimizer of the linear-programming problem may not be unique, and the solver may select a distributed admissible measure with the same prescribed first moment. This does not contradict the deterministic nature of the underlying problem, since the evolved conservative variable is still given by the first moment of the Young measure.
\end{rmk}

\subsubsection{Pressureless Gas Dynamics System}
In this section, we consider the 1-D pressureless gas dynamics system, which reads as 
\begin{equation}\label{5.6}
\begin{cases}
 \rho_t +q_x=0,& \\[1.ex]
 q_t +\big(\frac{q^2}{\rho}\big)_x =0, &
\end{cases}
\end{equation}
where $\rho$ and $v:=\frac{q}{\rho}$ are the density and velocity, respectively. The corresponding entropy function $\eta$ in \eref{linprog-a} is defined by $\eta=\frac{1}{2}\rho v^2$. 

It should be noted that the 1-D pressureless gas dynamics system is weakly hyperbolic and its solutions develop delta-shocks and vacuum states in finite time, which make the design of stable and accurate numerical methods very challenging.

\paragraph*{Example 6.} In this example, we consider the pressureless gas dynamics system \eref{5.6} and take the following initial data:
\begin{equation*}
 \bmu(x,\xi,0) =
\begin{cases}
  (1,2), & \mbox{if } x<-0.5 \\
  (1,0), & \mbox{otherwise},
\end{cases}
\end{equation*}
prescribed in the computational domain $[-1,1]$ subject to the free boundary conditions. We note that a sharp peak forms at later times in both $\rho$ and $q$,
corresponding to the formation of a delta-shock.

\begin{figure}[ht!]
\centerline{\includegraphics[trim=1.1cm 0.4cm 1.3cm 0.6cm, clip, width=5.cm]{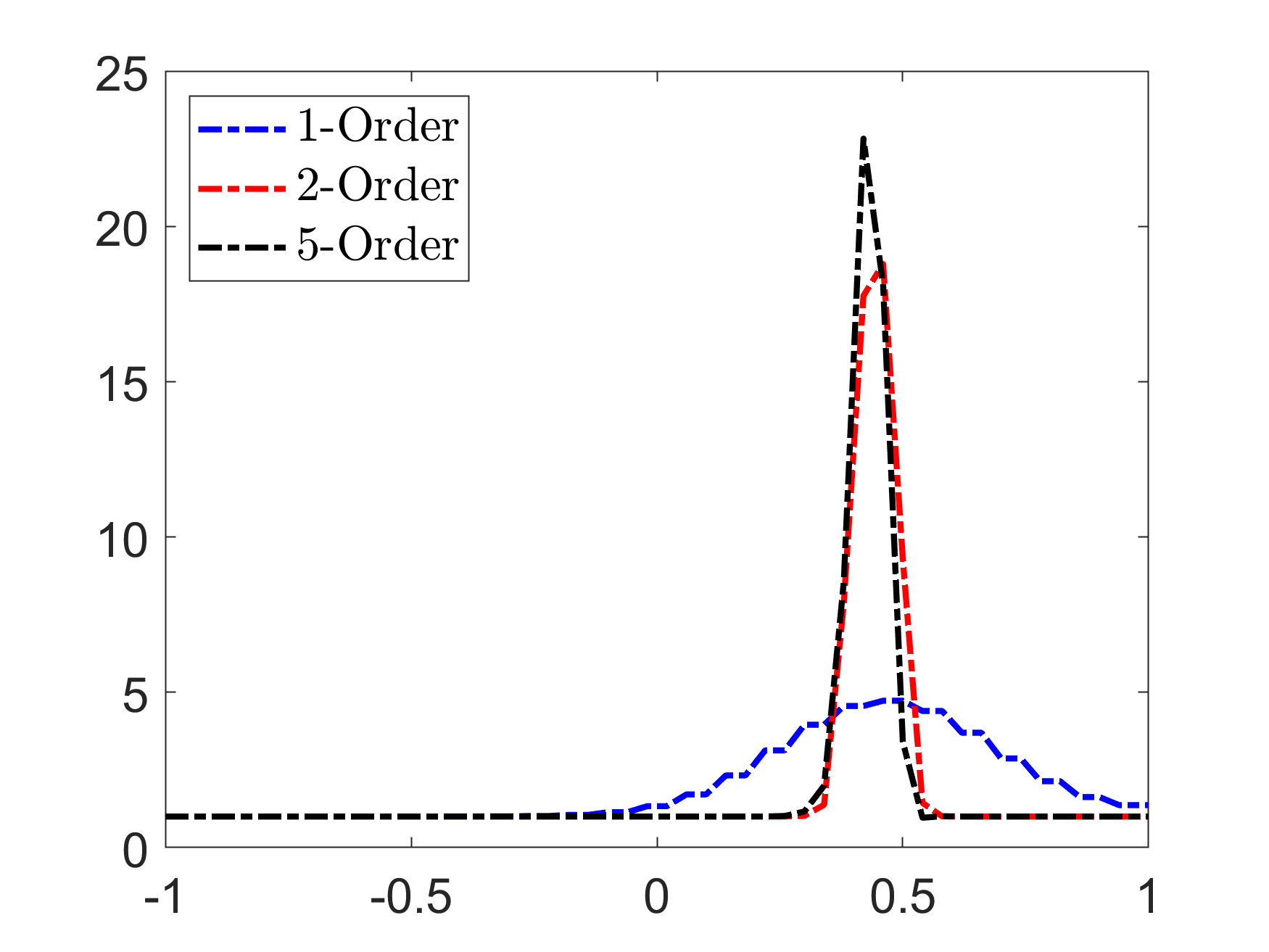} \hspace{1cm}
            \includegraphics[trim=1.1cm 0.4cm 1.3cm 0.6cm, clip, width=5.cm]{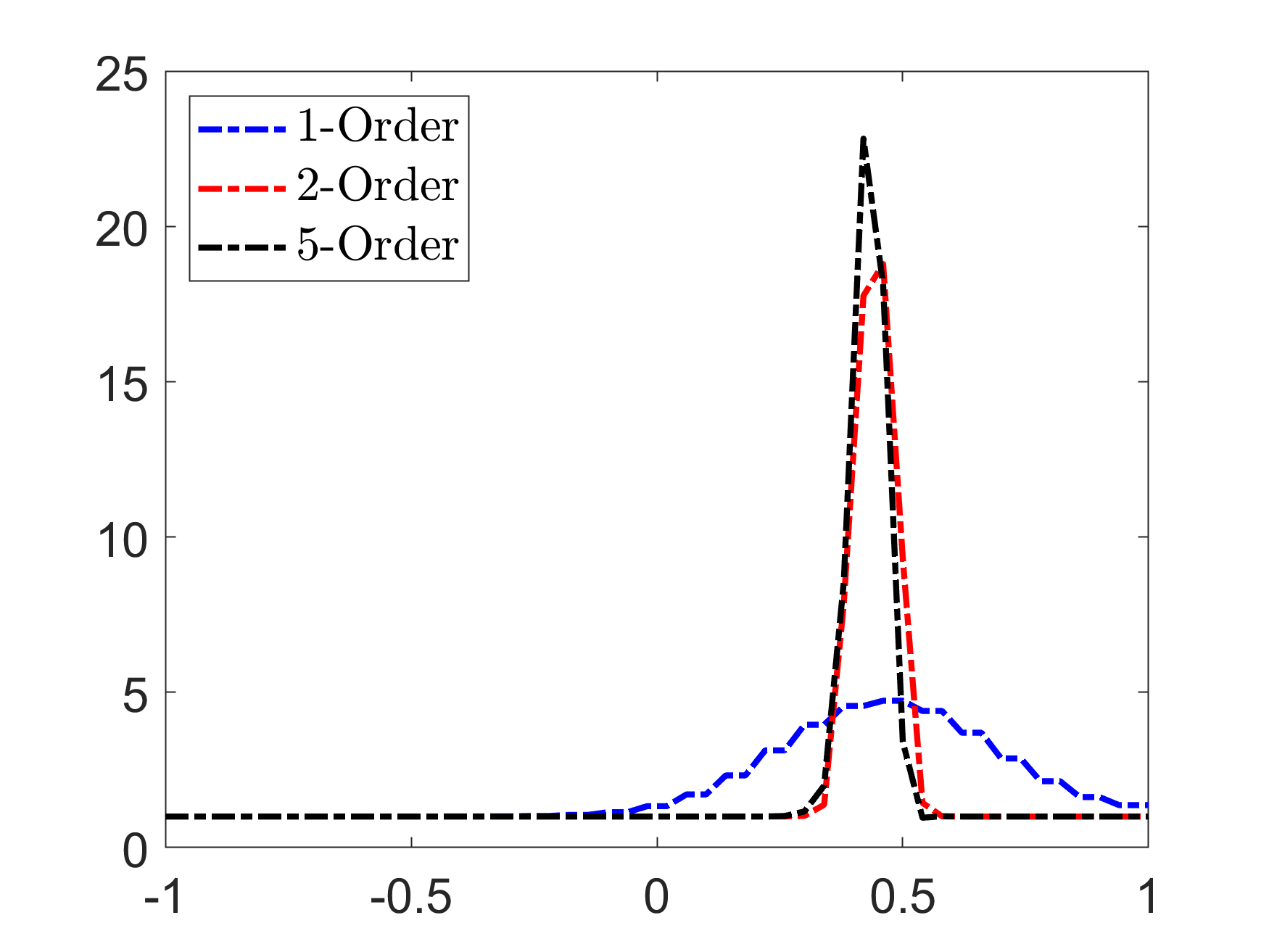} }   
\caption{\sf Example 6: Numerical results of $\rho$ (left) and $q$ (right) computed by the 1-Order, 2-Order, and 5-Order Young-measure schemes. \label{fig3.12}}
\end{figure}

We compute the numerical solutions until the final time $t=1$ by the studied 1-Order, 2-Order, and 5-Order Young-measure schemes on the uniform mesh with $N_x=50$, $N_\xi=1$. For the discretization in phase space $\bmu=(\rho, q)$, we use $[0.2,30]\times [-0.3,30]$ with $N_{\bmu}=300$ for each variable in the 1-Order, 2-Order, and 5-Order schemes, respectively. The obtained numerical results are presented in Figure \ref{fig3.12}, where one can see that the higher-order schemes produce much sharper delta-shocks than the 1-Order counterpart, demonstrating the robustness of the studied Young-measure schemes.
\medskip 

\begin{figure}[ht!]
\centerline{\includegraphics[trim=1.1cm 0.4cm 1.3cm 0.6cm, clip, width=5.cm]{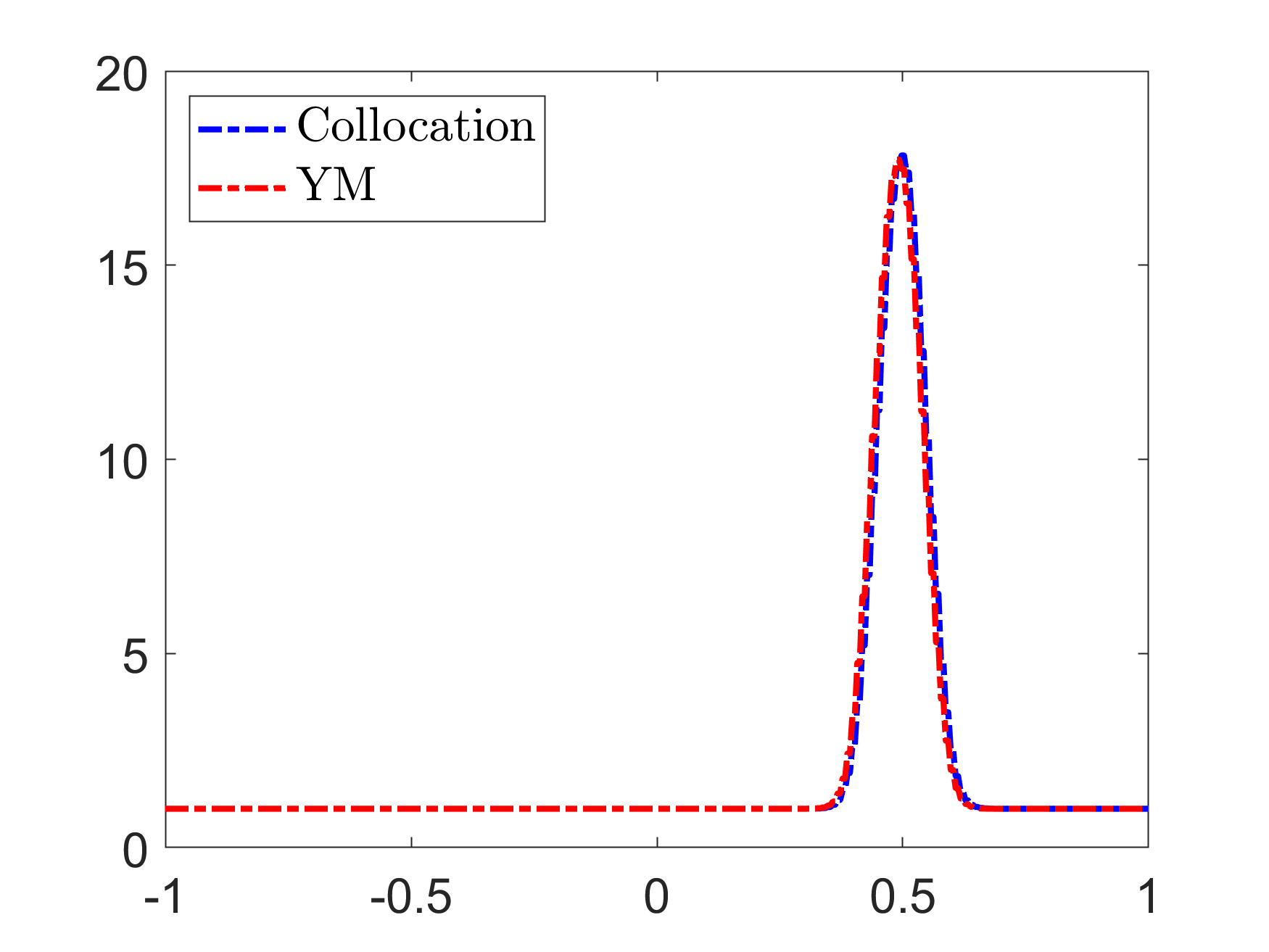} \hspace{1cm}
            \includegraphics[trim=1.1cm 0.4cm 1.3cm 0.6cm, clip, width=5.cm]{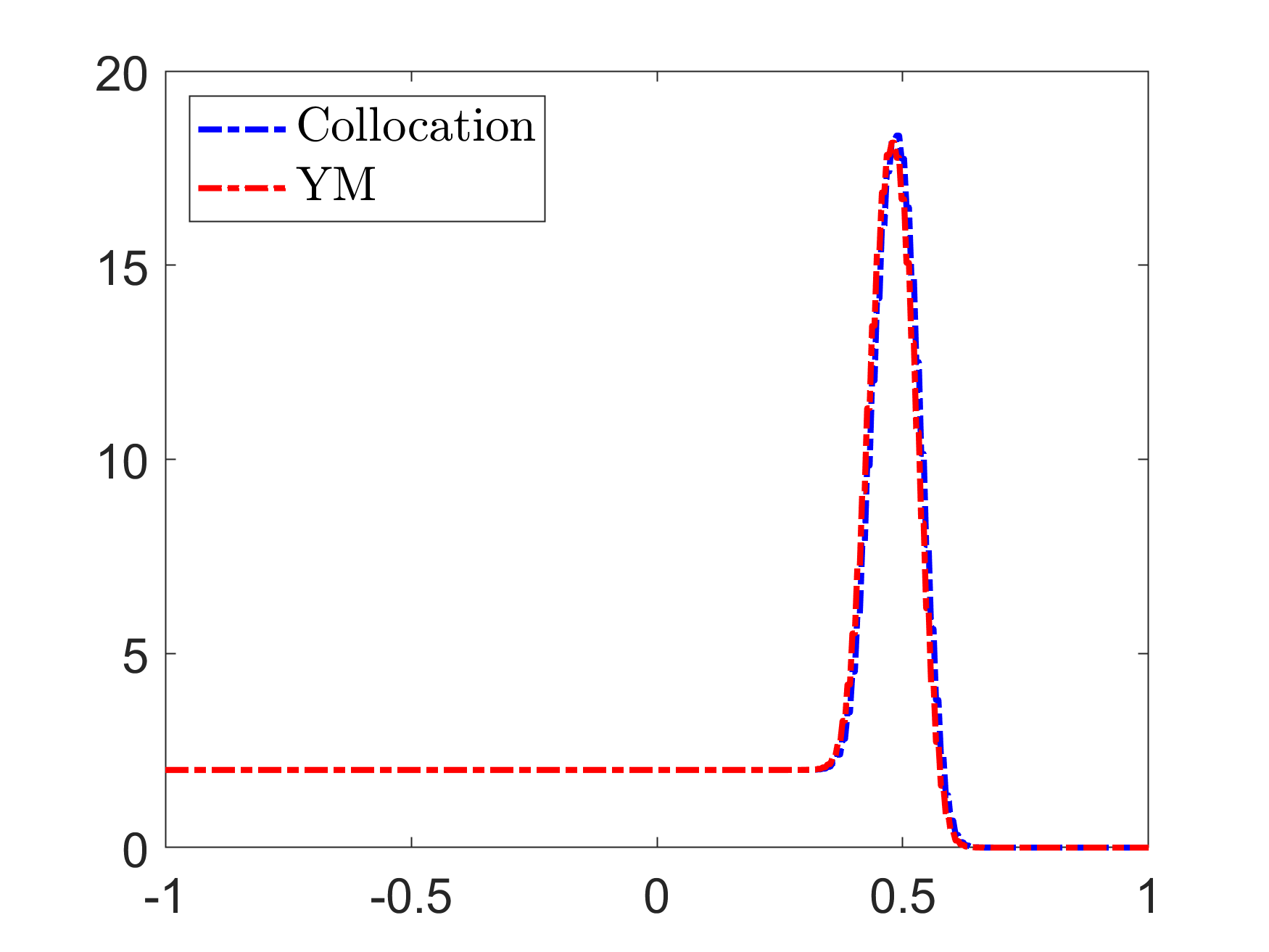} }   
\caption{\sf Example 6: Numerical results of $\rho$ (left) and $q$ (right) computed by the 1-Order collocation and Young-measure schemes. \label{fig3.12a}}
\end{figure}

\begin{figure}[ht!]
\centerline{\includegraphics[trim=1.1cm 0.4cm 1.3cm 0.6cm, clip, width=5.cm]{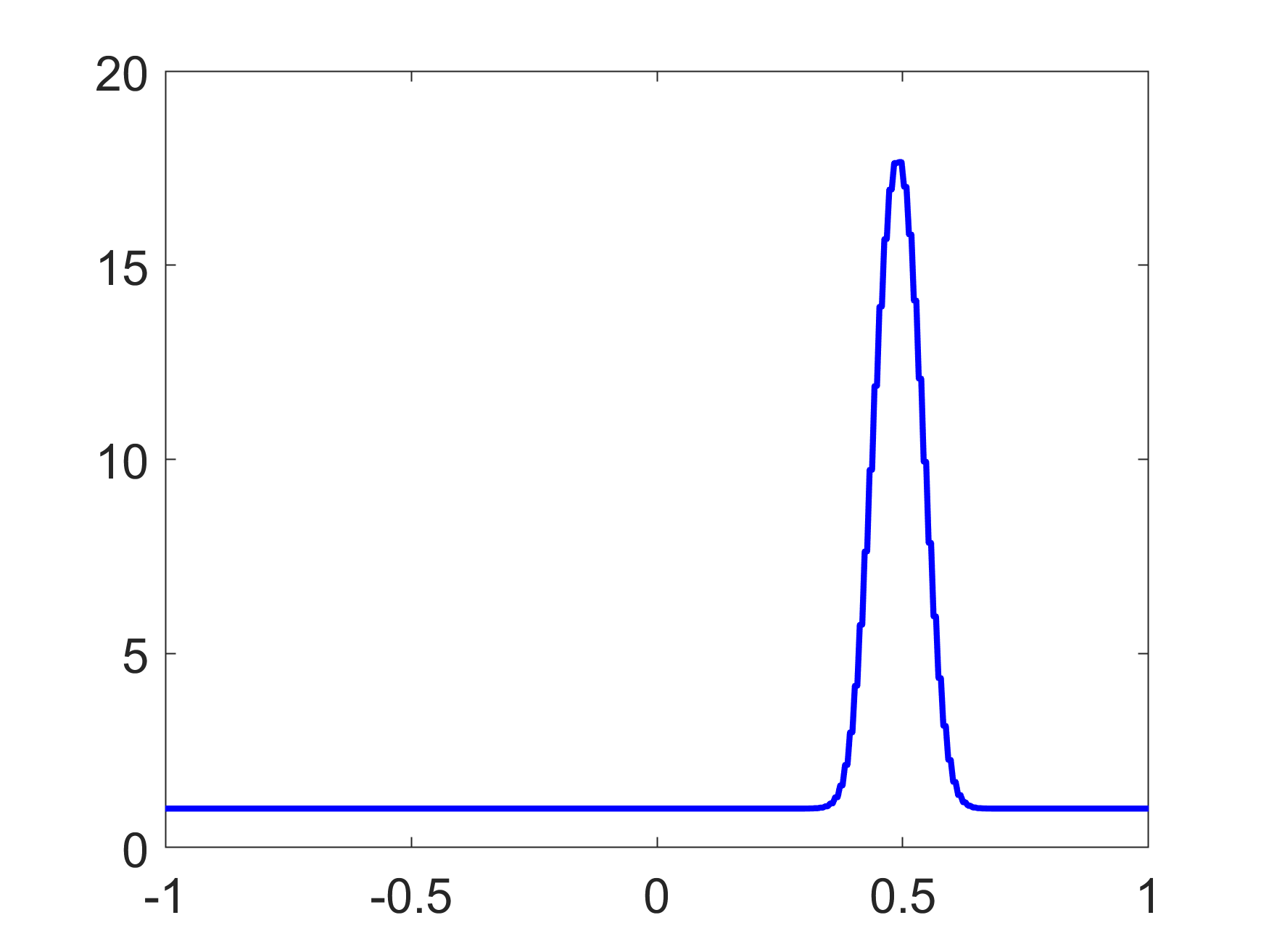} \hspace{1cm}
            \includegraphics[trim=1.1cm 0.4cm 1.3cm 0.6cm, clip, width=5.cm]{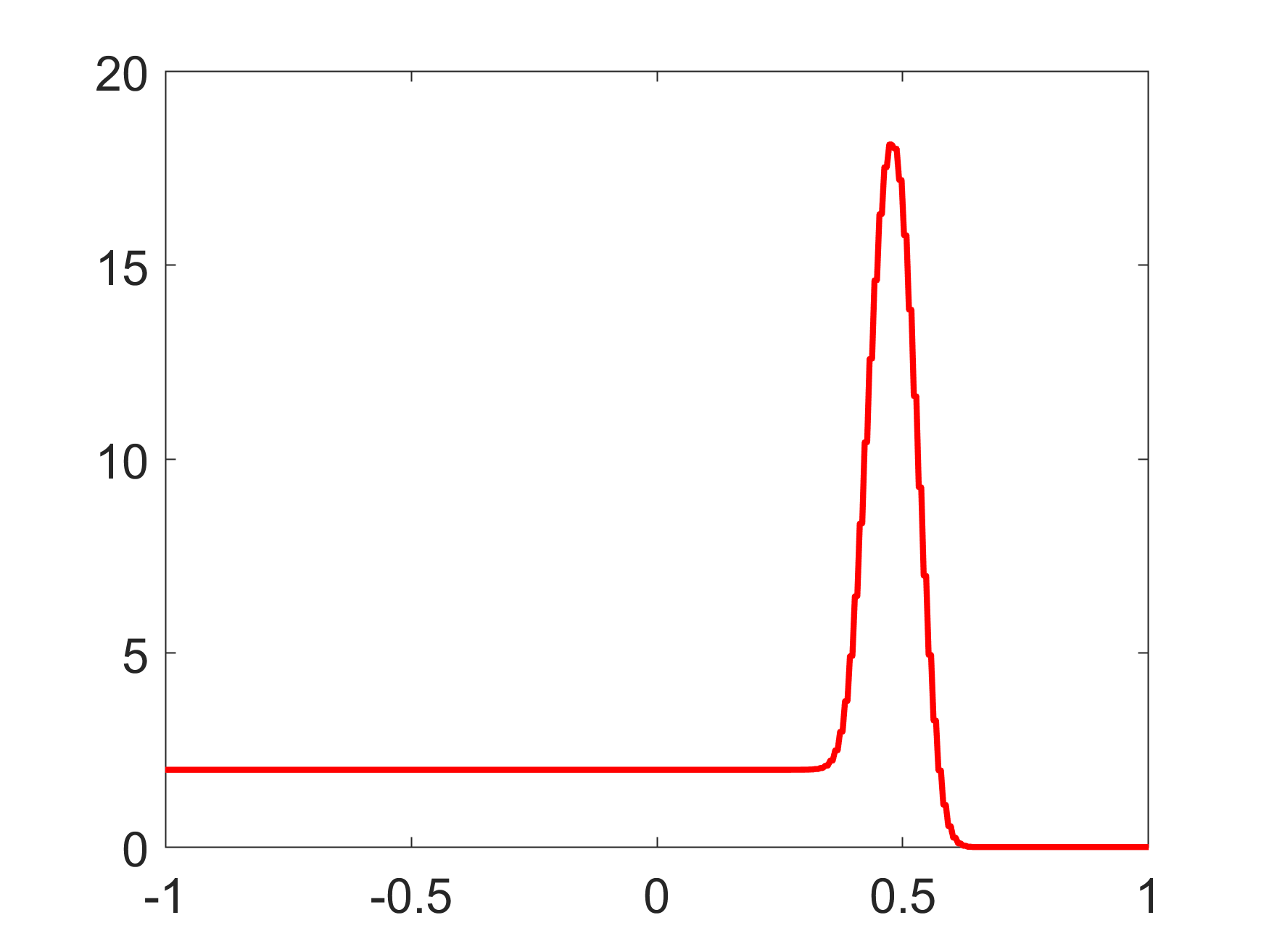} }   
\caption{\sf Example 6: Moment of $\rho$ (left) and $q$ (right) computed by the 1-Order Young-measure scheme.\label{fig3.12b}}
\end{figure}

\begin{figure}[ht!]
\centerline{\includegraphics[trim=0.2cm 0.2cm 0.9cm 0.8cm, clip, width=5.cm]{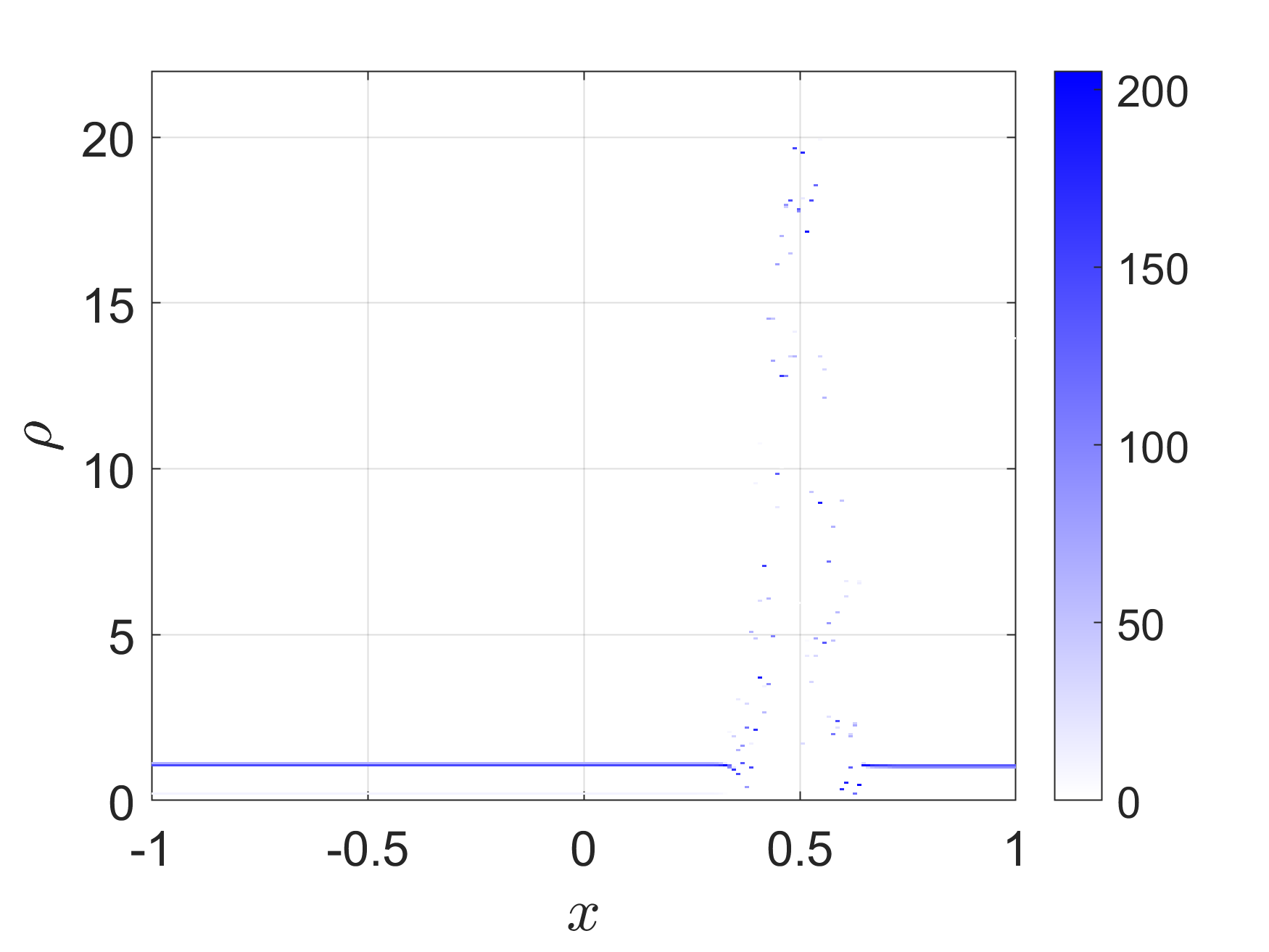} \hspace{1cm}
            \includegraphics[trim=0.2cm 0.2cm 0.9cm 0.8cm, clip, width=5.cm]{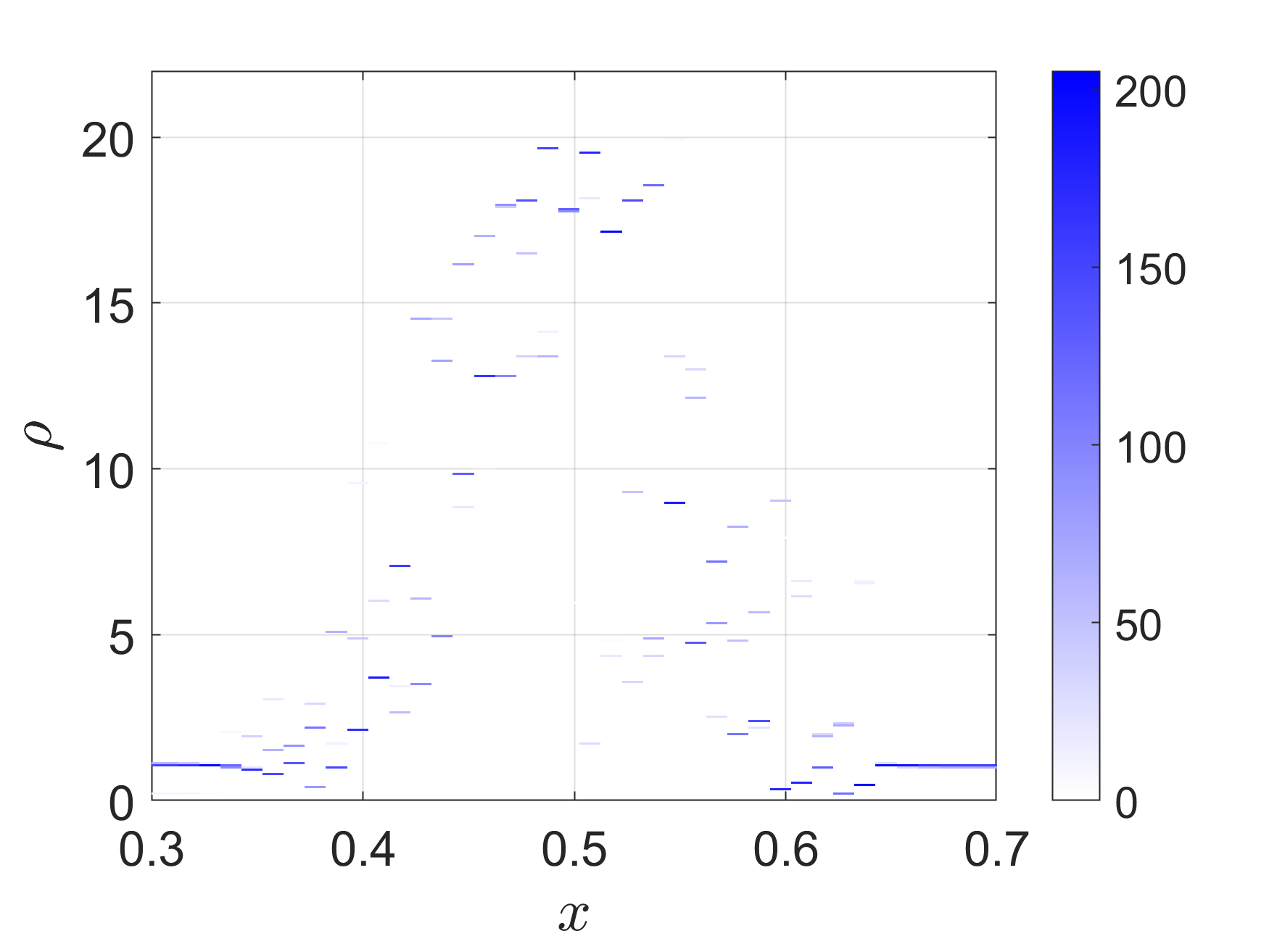} }   
\vskip 10pt
\centerline{\includegraphics[trim=0.2cm 0.2cm 0.9cm 0.8cm, clip, width=5.cm]{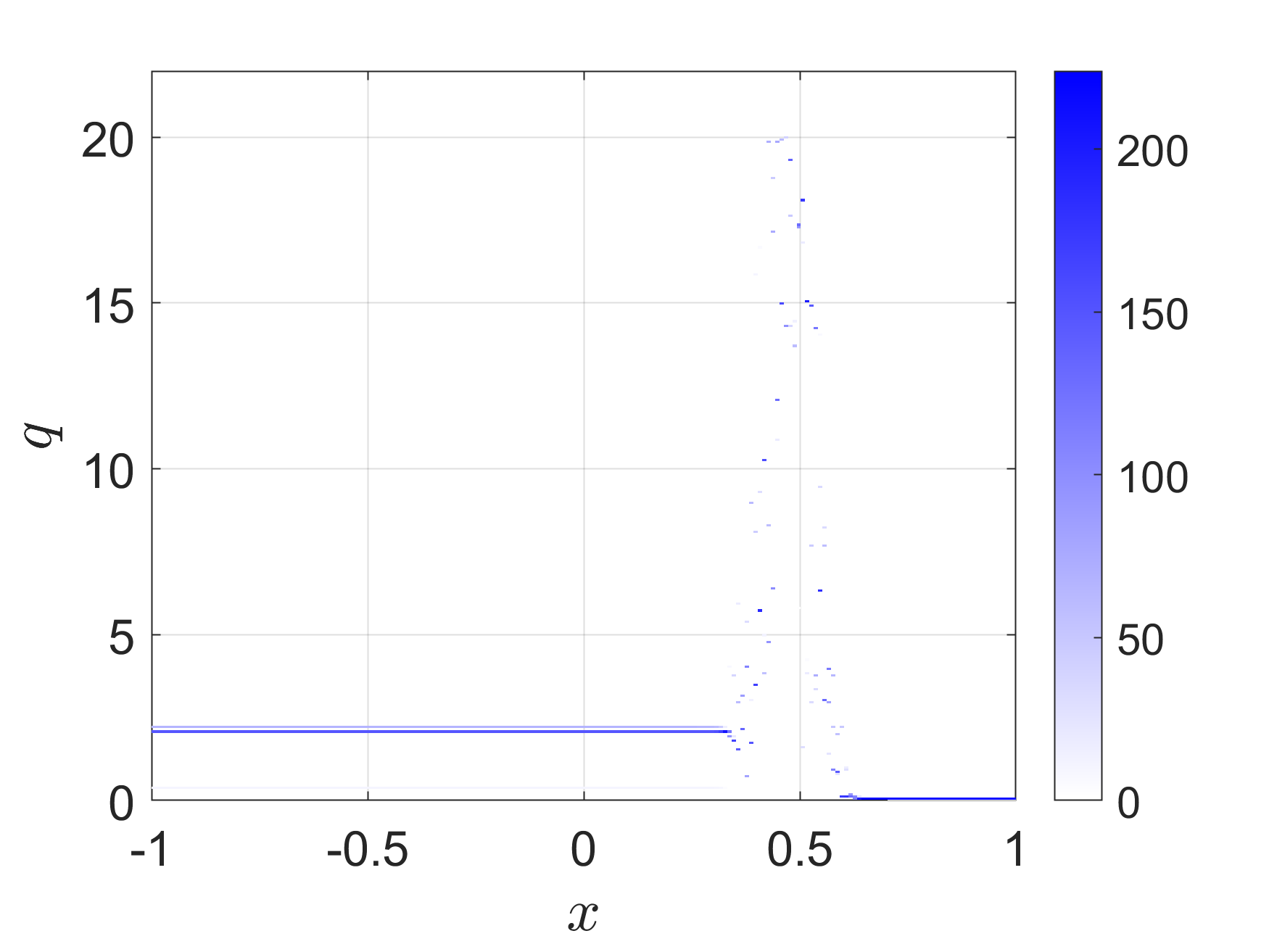} \hspace{1cm}
            \includegraphics[trim=0.2cm 0.2cm 0.9cm 0.8cm, clip, width=5.cm]{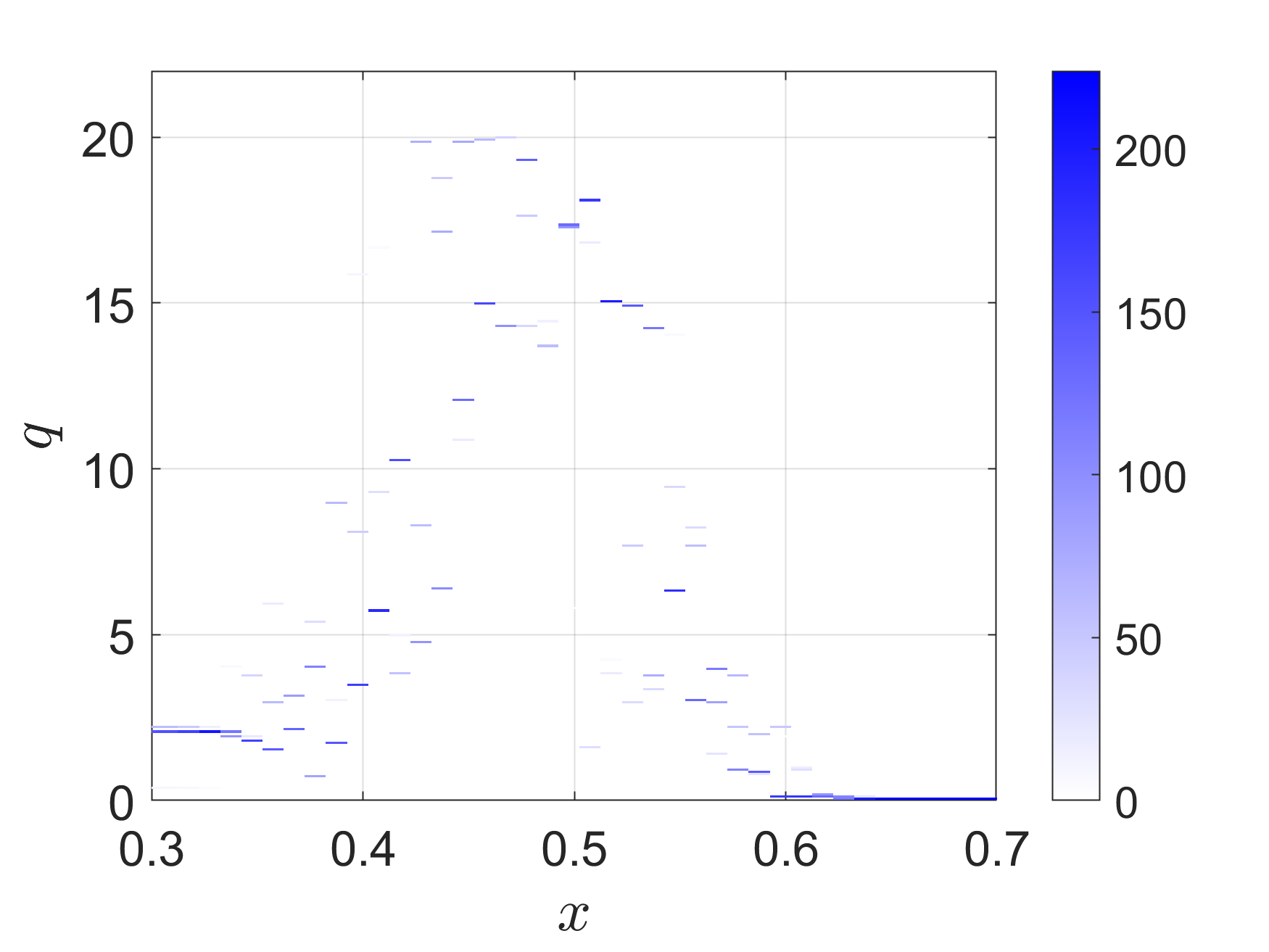} }   
\caption{\sf Example 6: $\widetilde{\rho}$ (top row) and $\widetilde{q}$ (bottom row) and zoom at $x\in [0.3, 0.7]$  computed by the 1-Order Young-measure scheme.\label{fig3.12c}}
\end{figure}

To further investigate, we also compute numerical results using the first-order scheme on a uniform mesh with $N_x=400$, $N_\xi=1$, and set $[0.2,20]\times [-0.3,20]$ with $N_{\bmu}=300$ for the discretization in phase space $\bmu=(\rho, q)$. The obtained numerical results are presented in Figure \ref{fig3.12a}. As observed, both the Young-measure and collocation schemes capture sharper delta-shocks; however, the resolution is still inferior to that obtained by the 2-Order and 5-Order schemes, even on a finer mesh. These results highlight the superior efficiency of higher-order schemes. To further illustrate the information contained in the reconstructed Young measure, we also plot
the first moments $\iint_{\mathbb{R}\times\mathbb{R}}\rho\,d\mu(\rho,q)$ and  $\iint_{\mathbb{R}\times\mathbb{R}}q\,d\mu(\rho,q)$, as well as the marginal distributions with respect to $\rho$ and $q$,
$$
\widetilde{\rho}(x,\rho)=\int_{\mathbb{R}}\mu^*(\rho,q;x)\,dq,
\qquad
\widetilde{q}(x,q)=\int_{\mathbb{R}}\mu^*(\rho,q;x)\,d\rho.
$$
The corresponding results are shown in Figures \ref{fig3.12b}--\ref{fig3.12c}. Figure \ref{fig3.12b} shows the first moments, which agree with the sharp density and momentum peaks generated by the delta-shock. Figure \ref{fig3.12c} provides complementary phase-space information: the marginal distributions show where the reconstructed Young measure places its mass in the density and momentum directions. In particular, the zoomed plots near the delta-shock indicate that the large moment values are produced by a localized concentration of the reconstructed measure in phase space. Thus, Figure \ref{fig3.12c} illustrates the phase-space representation of the delta-shock provided by the Young-measure reconstruction.

\medskip 

\medskip 

\subsubsection{Burgers Equation with Non-Atomic Support}
In the last subsection, we consider an example taken from \cite{CHLZ25}, which is inspired by \cite{Tadmor_1Dburgers}, where the Young measure may have non-atomic support. We set the factor $\lambda_F$ in \eqref{linprog-c} to 0.05. Its purpose is to demonstrate that, by imposing the upper bound with $\lambda_F<1$, the linear-programming reconstruction can produce a Young measure with distributed, non-atomic support rather than a purely atomic measure.

\begin{figure}[ht!]
\centerline{\includegraphics[trim=0.2cm 0.2cm 1.1cm 0.7cm, clip, width=7.cm]{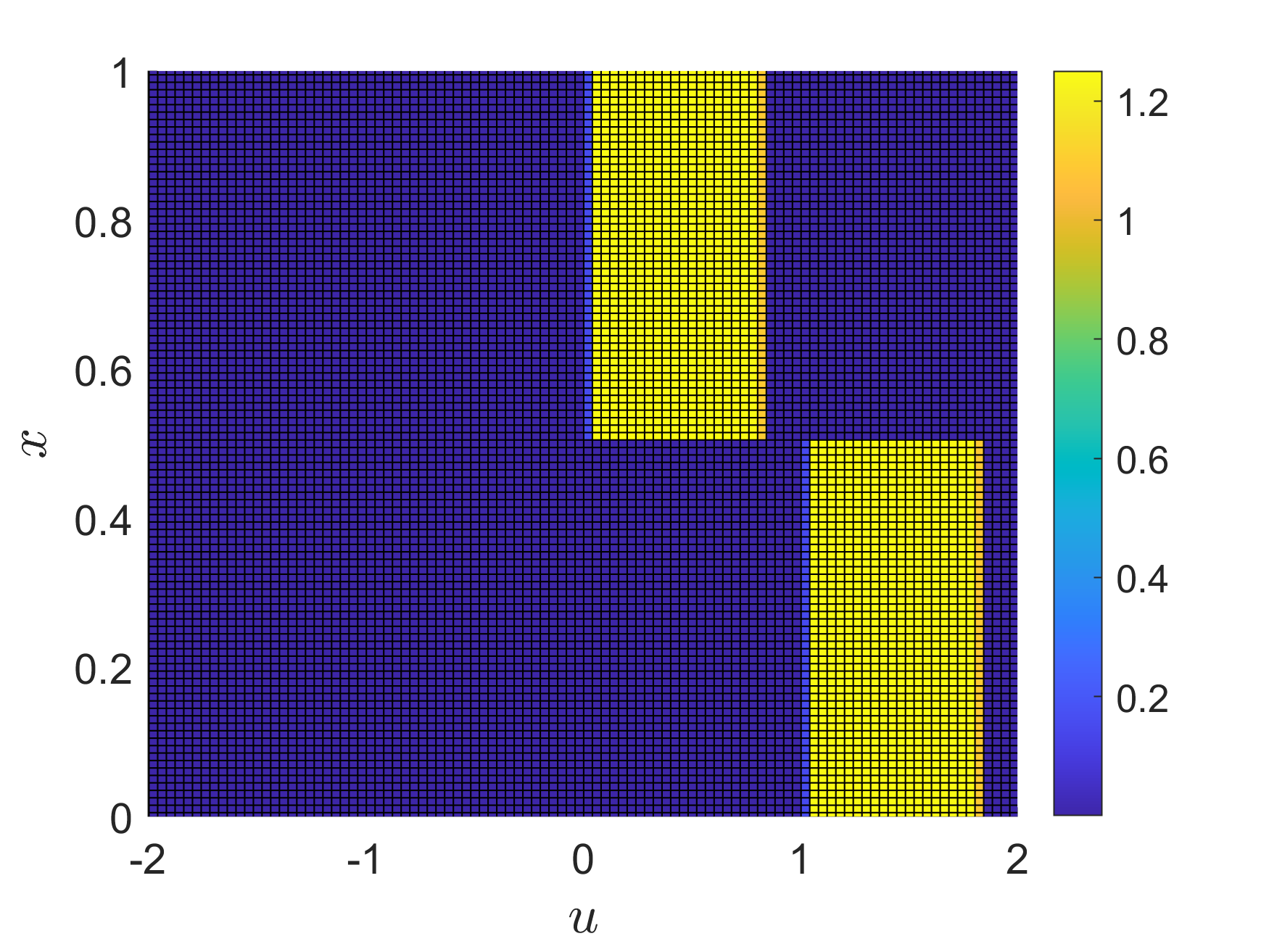}}
\caption{\sf Example 7: Support of Young measures at the initial time.\label{fig9a}}
\end{figure}

\paragraph{Example 7.} We consider the following initial data
$$
u_0(x,\xi)=\begin{cases}
             1.5 & \mbox{if } x< 0.5 \\
             0.5 & \mbox{otherwise},
           \end{cases}
$$
subject to the free boundary conditions imposed in the computational domain $[-1,1]$. In Figure \ref{fig9a}, we plot the support of the Young measures at the initial time. 

\begin{figure}[ht!]
\centerline{\includegraphics[trim=0.2cm 0.2cm 0.8cm 0.1cm, clip, width=5.cm]{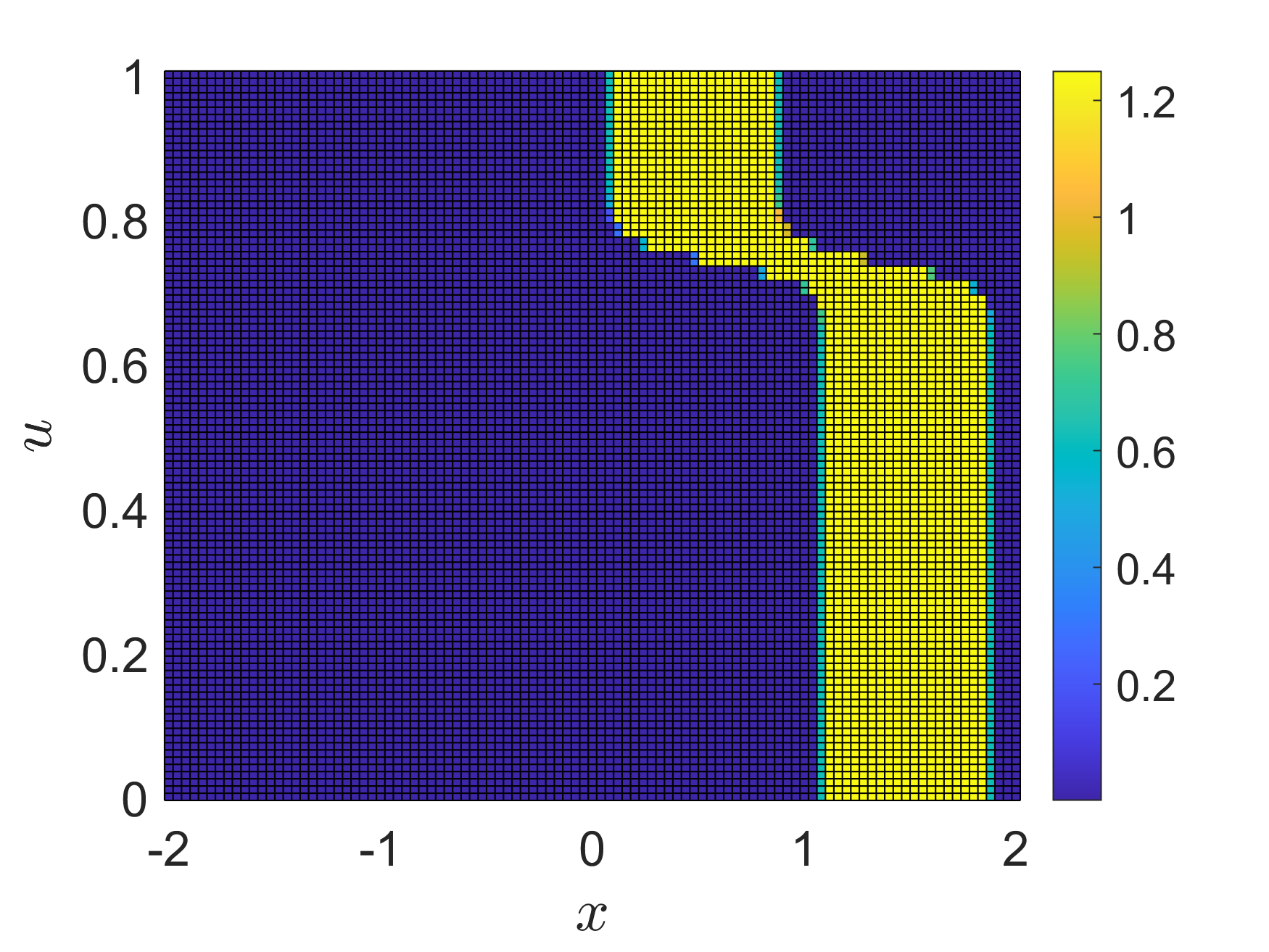}\hspace{0.5cm}
            \includegraphics[trim=0.2cm 0.2cm 0.8cm 0.1cm, clip, width=5.cm]{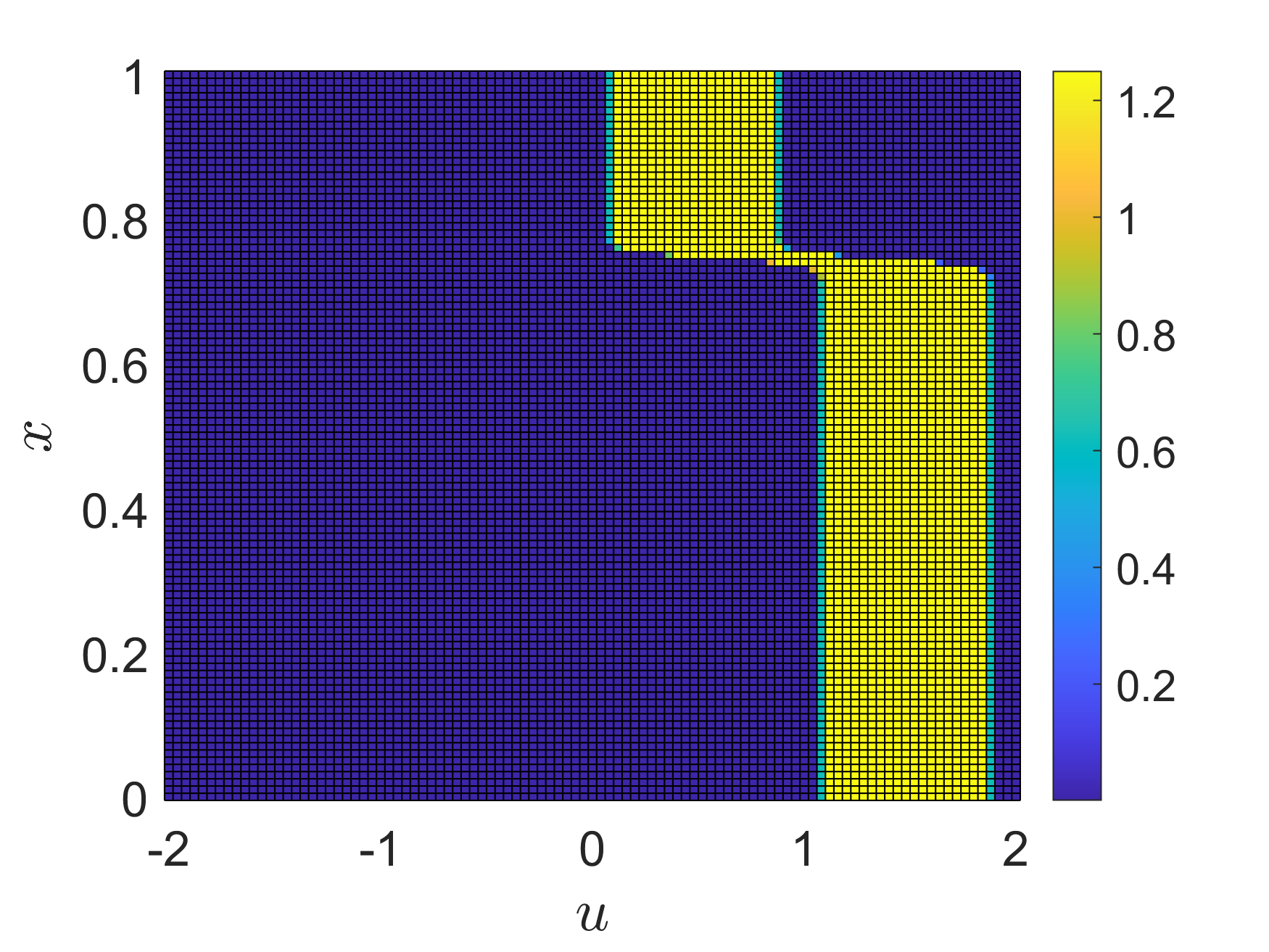}\hspace{0.5cm}
            \includegraphics[trim=0.2cm 0.2cm 0.8cm 0.1cm, clip, width=5.cm]{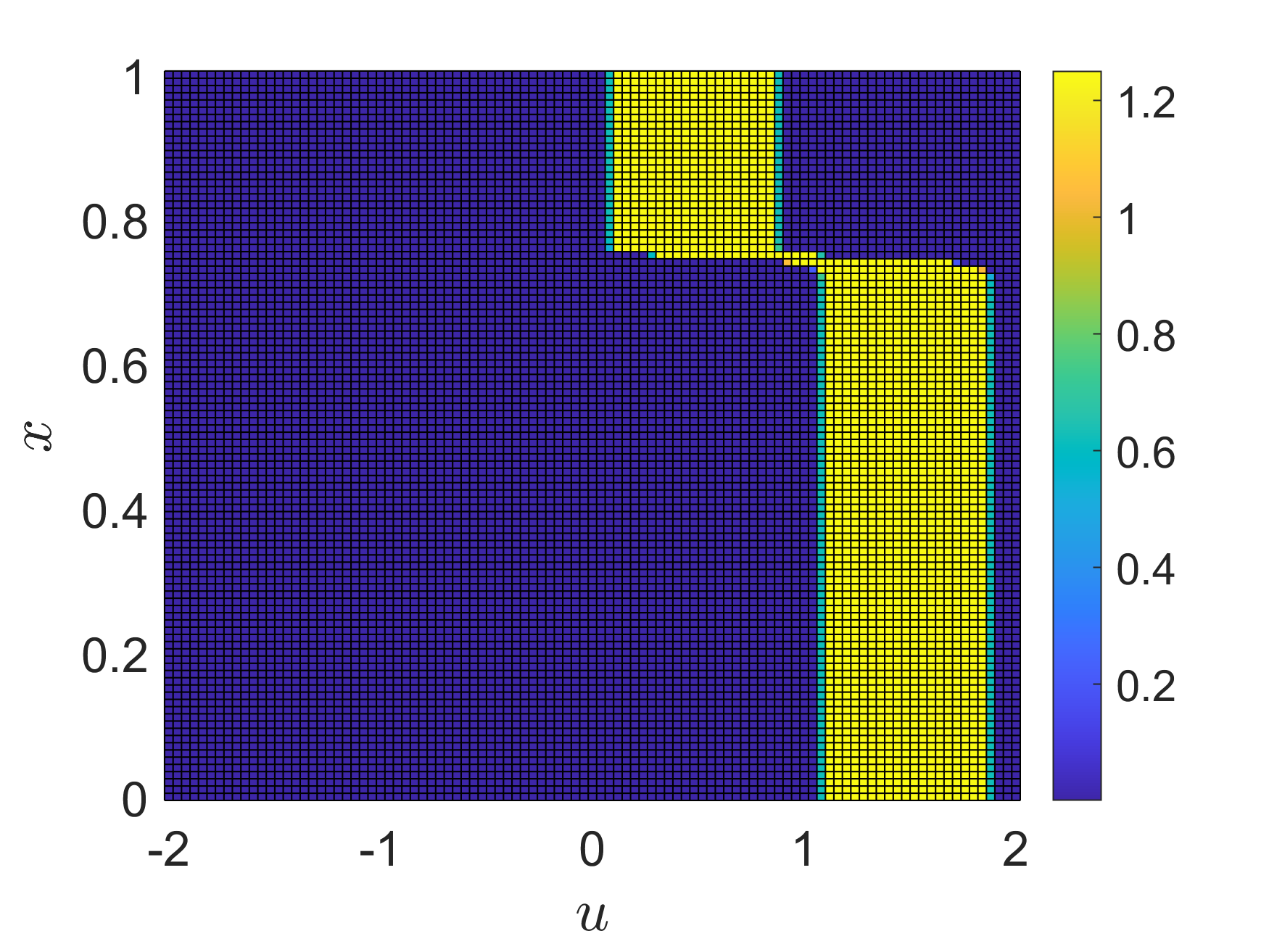}}
\caption{\sf Example 7: Support of Young measures computed by the 1- (left), 2- (middle), and 5-Order (right) schemes at the final time $t=0.25$.\label{fig9b}}
\end{figure}

We compute the numerical results until the final time $t=0.25$ by the studied 1-, 2-, and 5-Order schemes on the uniform mesh with $N_x=100$ and $N_\xi=1$. The phase space is confined to the interval $u \in [-2, 2]$ and is uniformly discretized into $N_u = 100$ equidistant cells. The obtained numerical results are presented in Figure \ref{fig9b}, where one can clearly see that the measure with non-atomic support can be captured by the linear-programming algorithm if we choose $\lambda_F<1$, and in contrast to Example 1, the support of the Young-measure solution is distributed rather than concentrated. At the same time, the use of higher-order schemes yields a significantly sharper resolution.

\section{Conclusions}
In this paper, we extended the existing first-order scheme for solving random hyperbolic conservation laws using linear programming to higher-order schemes. We first extended the first-order scheme to second order by employing piecewise linear reconstruction and then to fifth order of accuracy in the framework of finite-difference A-WENO scheme. These extensions enhance the resolution of discontinuities while preserving structure-preserving properties. We validated the proposed methods through a series of random and deterministic hyperbolic conservation laws, including the one-dimensional Burgers equation, isentropic Euler system, the scalar conservation law with a discontinuous flux function, pressureless gas dynamics system, and Burgers equation with non-atomic support. The obtained numerical results demonstrate that higher-order schemes significantly improve the resolution of complex wave structures and reduce numerical dissipation, particularly near shock and rarefaction waves. Future work includes extending these approaches to (nonconservative) multidimensional random hyperbolic problems and improving computational efficiency  via adaptive refinement strategies.
\begin{DA}
\paragraph{Acknowledgments.} 
The authors would like to thank Giuseppe Coclite (Politecnico di Bari, Italy) for providing the idea of the example on the discontinuous flux as well as interesting discussions on the topic.
\par
The work of S. Chu was supported in part by the DFG (German Research Foundation) through HE5386/19-3, 27-1. The work of M. Herty was funded by the DFG--SPP 2183: Eigenschaftsgeregelte Umformprozesse with the Project(s) HE5386/19-2,19-3 Entwicklung eines flexiblen isothermen Reckschmiedeprozesses f\"ur die eigenschaftsgeregelte Herstellung von Turbinenschaufeln aus Hochtemperaturwerkstoffen (424334423) and by the Deutsche Forschungsgemeinschaft (DFG, German Research Foundation)--SPP 2410 Hyperbolic Balance Laws in Fluid Mechanics: Complexity, Scales, Randomness (CoScaRa) within the Project(s) HE5386/26-1 (Numerische Verfahren f\"ur gekoppelte Mehrskalenprobleme,525842915) and (Zuf\"allige kompressible Euler Gleichungen: Numerik und ihre Analysis, 525853336) HE5386/27-1.

\paragraph{Conflicts of interest.} On behalf of all authors, the corresponding author states that there is no conflict of interest.

\paragraph{Data and software availability.} The data that support the findings of this study and MATLAB codes developed by the authors and
used to obtain all of the presented numerical results are available by the corresponding author upon  request.
\end{DA}

\appendix
\section{1-D Local Characteristic Decomposition Based Piecewise Linear Interpolant}\label{appa}
In this appendix, we briefly describe the 1-D local characteristic decomposition (LCD) based piecewise linear interpolant.

In order to suppress the oscillations in the piecewise linear reconstruction \eref{3.3}--\eref{3.6}, we apply the generalized minmod limiter to the local characteristic variables. To this end, we first introduce the matrix $\widehat{A}_{i,\jph}=A(\hat{\bmu}_{i,\jph})$, where $A=\frac{\partial \mf}{\partial \bmu}$ is the Jacobian and $\hat{\bmu}_{i,\jph}$ is either a simple average $(\bmu_{i,j}+\bmu_{i,j+1})/2$ or another type of average of the $\bmu_{i,j}$ and $\bmu_{i,j+1}$ states. 

We then compute the matrices $R_{i,\jph}$ and $R^{-1}_{i,\jph}$ such that $R_{i,\jph} \widehat{A}_{i,\jph}R^{-1}_{i,\jph}$ is a diagonal matrix and introduce the local characteristic variables $\bm\Gamma_k$ in the neighborhood of $x=x_\jph$:
$$
\bm\Gamma_k = R^{-1}_{i,\jph} \bmu_{i,k}, \quad k=j-1,\ldots,j+2.
$$
Equipped with the values $\bm\Gamma_{i,j-1}$, $\bm\Gamma_{i,j}$, $\bm\Gamma_{i,j+1}$, and $\bm\Gamma_{i,j+2}$, we reconstruct $\bm\Gamma$ by computing 
$$
(\bm\Gamma_x)_{i,j}={\rm minmod}\left(\theta\frac{\,\bm\Gamma_{i,j}-\,\bm\Gamma_{i,j-1}}{\dx},\,\frac{\,\bm\Gamma_{i,j+1}-\,\bm\Gamma_{i,j-1}}{2\dx},\,
\theta\frac{\,\bm\Gamma_{i,j+1}-\,\bm\Gamma_{i,j}}{\dx}\right),
$$
and 
$$
(\bm\Gamma_x)_{i,j+1}={\rm minmod}\left(\theta\frac{\,\bm\Gamma_{i,j+1}-\,\bm\Gamma_{i,j}}{\dx},\,\frac{\,\bm\Gamma_{i,j+2}-\,\bm\Gamma_{i,j}}{2\dx},\,
\theta\frac{\,\bm\Gamma_{i,j+2}-\,\bm\Gamma_{i,j+1}}{\dx}\right),
$$
where the minmod function is defined by \eref{3.6}. After that, we evaluate 
$$
\bm\Gamma^{-}_{i,\jph}=\bm\Gamma_{i,j}+\frac{\dx}{2}(\bm\Gamma_x)_{i,j} \quad {\rm and} \quad \bm\Gamma^{+}_{i,\jph}=\bm\Gamma_{i,j+1}-\frac{\dx}{2}(\bm\Gamma_x)_{i,j+1}.
$$
Finally, one can obtain the corresponding point values of $\bmu$ by 
$$
\bmu^{-}_{i,\jph}=R_{i,\jph}\bm\Gamma^{-}_{i,\jph} \quad {\rm and} \quad \bmu^{+}_{i,\jph}=R_{i,\jph}\bm\Gamma^{+}_{i,\jph}. 
$$

\begin{rmk}
The one-sided point values of $\bmu^{\pm,*}_{i,\jph}$ can be computed in the same manner by using $\bmu^*_{i,j}$, and here we omit it for the sake of brevity.
\end{rmk}

\section{1-D Local Characteristic Decomposition Based Fifth-Order WENO-Z Interpolant}\label{appb}
In this appendix, we briefly describe the 1-D LCD based fifth-order WENO-Z interpolant.

Given the point values $\bmu_{i,j}$, for the $k$-th component of $\bmu_{i,j}$, the value $u^{(k),-}_{i,\jph}$ is computed using a weighted average of the three parabolic interpolants ${{\cal P}}_0(x)$, ${{\cal P}}_1(x)$, and ${{\cal P}}_2(x)$ obtained using the stencils $[x_{j-2},x_{j-1},x_j]$, $[x_{j-1},x_j,x_{j+1}]$, and $[x_j,x_{j+1},x_{j+2}]$, respectively:
\begin{equation}
u^{(k),-}_{i,\jph}=\sum_{\ell=0}^2\omega_\ell{{\cal P}}_\ell(x_\jph),
\label{b1}
\end{equation}
where
\begin{equation}
\begin{aligned}
&{{\cal P}}_0(x_\jph)=\frac{3}{8}u^{(k)}_{i,j-2}-\frac{5}{4}u^{(k)}_{i,j-1}+\frac{15}{8}u^{(k)}_{i,j},\\
&{{\cal P}}_1(x_\jph)=-\frac{1}{8}u^{(k)}_{i,j-1}+\frac{3}{4}u^{(k)}_{i,j}+\frac{3}{8}u^{(k)}_{i,j+1},\\
&{{\cal P}}_2(x_\jph)=\frac{3}{8}u^{(k)}_{i,j}+\frac{3}{4}u^{(k)}_{i,j+1}-\frac{1}{8}u^{(k)}_{i,j+2}.
\end{aligned}
\label{b2}
\end{equation}
In order to ensure \eref{b1}--\eref{b2} is fifth-order accurate and non-oscillatory, one can take the weights $\omega_\ell$
in \eref{b1} to be
\begin{equation}
\omega_\ell:=\frac{\alpha_\ell}{\alpha_0+\alpha_1+\alpha_2},\quad
\alpha_\ell=d_\ell\left[1+\left(\frac{\tau_5}{\beta_\ell+\varepsilon}\right)^p\right],\quad\tau_5=|\beta_2-\beta_0|,
\label{b3}
\end{equation}
where $d_0=\frac{1}{16}$, $d_1=\frac{5}{8}$, and $d_2=\frac{5}{16}$, and the smoothness indicators $\beta_\ell$ in \eref{b3} are given by 
\begin{equation}
\begin{aligned}
&\beta_0=\frac{13}{12}\big(u^{(k)}_{i,j-2}-2u^{(k)}_{i,j-1}+u^{(k)}_{i,j}\big)^2+\frac{1}{4}\big(u^{(k)}_{i,j-2}-4u^{(k)}_{i,j-1}+3u^{(k)}_{i,j}\big)^2,\\
&\beta_1=\frac{13}{12}\big(u^{(k)}_{i,j-1}-2u^{(k)}_{i,j}+u^{(k)}_{i,j+1}\big)^2+\frac{1}{4}\big(u^{(k)}_{i,j-1}-u^{(k)}_{i,j+1}\big)^2,\\
&\beta_2=\frac{13}{12}\big(u^{(k)}_{i,j}-2u^{(k)}_{i,j+1}+u^{(k)}_{i,j+2}\big)^2+\frac{1}{4}\big(3u^{(k)}_{i,j}-4u^{(k)}_{i,j+1}+u^{(k)}_{i,j+2}\big)^2.
\end{aligned}
\label{b4}
\end{equation}
In all of the numerical examples reported in this paper, we have used $p=2$ and $\varepsilon=10^{-12}$. The corresponding right-sided value $u^{(k),+}_{i,\jph}$ can also be derived using a mirror-symmetric approach, and here we omit it for the sake of brevity.  

As in Appendix \ref{appa}, to ensure the nonoscillatory nature of the reconstruction \eref{b1}--\eref{b4}, we also need to adopt the reconstruction procedure within the LCD framework. To this end, with the matrices $R_{i,\jph}$ and $R^{-1}_{i,\jph}$ (see Appendix \ref{appc}), we first introduce the local characteristic variables in the neighborhood of $x=x_\jph$:
$$
\bm \Gamma_{\ell}=R^{-1}_{i,\jph}\bm u_{i,j+\ell},\quad \ell=-2,\ldots,3.
$$
Equipped with the values $\bm \Gamma_{-2}$, $\bm \Gamma_{-1}$,  $\bm \Gamma_{0}$, $\bm \Gamma_{1}$, $\bm \Gamma_{2}$, and $\bm \Gamma_{3}$, we then apply the interpolation procedure \eref{b1}--\eref{b4} to each of the components $\Gamma^{(k)}$, $k=1,\ldots,d$ of $\bm \Gamma$ to obtain $\bm {\Gamma}_{\hf}^{-}$ (the values of ${\bm \Gamma}_{\hf}^{+}$ are computed in the mirror-symmetric way). Finally, the corresponding one-sided point values of $\bm u^{\pm}_{i,\jph}$ are given by
\begin{equation*}
\bm u^{\pm}_{i,\jph}=R_{i,\jph}\bm \Gamma^{\pm}_{\hf}.
\end{equation*}

\begin{rmk}
The one-sided point values of $\bmu^{\pm,*}_{i,\jph}$ can be computed in the same manner by using $\bmu^*_{i,j}$, and here we omit it for the sake of brevity.
\end{rmk}

\section{Local Characteristic Decomposition Matrices for the 1-D Isentropic Euler System}\label{appc}
In this appendix, we introduce how to compute the LCD matrices for the 1-D isentropic Euler system \eref{5.3}--\eref{5.4}. 

We first compute the Jacobian 
\begin{equation*}
  A(\bmu)=\frac{\partial \mf}{\partial \bmu}=\begin{pmatrix}0&1\\
\kappa \gamma \rho^{\gamma-1}-v^2&2v\\\end{pmatrix},
\end{equation*}
and then introduce the matrices
$$
\begin{aligned}
\widehat A_{i,\jph}=\begin{pmatrix}0&1\\
\kappa \gamma \hat \rho^{\gamma-1}-\hat v^2&2\hat v\\
\end{pmatrix},
\end{aligned}
$$
where $\hat{(\cdot)}$ stands for the following averages:
\begin{equation*}
\hat\rho=\frac{\rho_{i,j}+ \rho_{i,j+1}}{2},~\hat v=\frac{v_{i,j}+v_{i,j+1}}{2},
\end{equation*}
with
$$
v_{i,j}=\frac{({\rho v})_{i,j}}{\rho_{i,j}}.
$$

Notice that all of the $\hat{(\cdot)}$ quantities have to have a subscript index, that is, $\hat{(\cdot)}=\hat{(\cdot)}_{i,\jph}$, but we omit
it for the sake of brevity for all of the quantities except for $\widehat A$. We then compute the matrix $R_{i,\jph}$ composed of the right
eigenvectors of $\widehat A_{i,\jph}$ and obtain
\begin{equation*}
R_{i,\jph}=\begin{pmatrix} \dfrac{1}{\hat v-\sqrt{\gamma}\hat \rho ^{(\gamma-1)/2}}&\dfrac{1}{\hat v+\sqrt{\gamma}\hat \rho ^{(\gamma-1)/2}}\\ 1&1\end{pmatrix}
\quad\mbox{and}\quad 
R^{-1}_{i,\jph}=\begin{pmatrix}
\dfrac{\hat \rho \hat v^2-\gamma \hat \rho^\gamma}{2\sqrt{\gamma} \hat \rho^{(\gamma+1)/2}}&\dfrac{1}{2}-\dfrac{\hat \rho^{(1-\gamma)/2 }\hat v}{2 \sqrt{\gamma}}\\
\dfrac{\gamma \hat \rho^\gamma-\hat \rho \hat v^2}{2\sqrt{\gamma} \hat \rho^{(\gamma+1)/2}}&\dfrac{1}{2}+\dfrac{\hat \rho^{(1-\gamma)/2}\hat v}{2 \sqrt{\gamma}}
\end{pmatrix}.
\end{equation*}

\section{A KKT-Based Solver for the Linear Programming Problem}\label{appd}
In this appendix, we introduce a KKT-based solver for the closure problem \eref{linprog-a}--\eref{linprog-e} in the case $\lambda_F=1$. In this case, the upper-bound constraint \eref{linprog-c} is redundant, since it follows from the nonnegativity and normalization constraints.

For each fixed spatial cell $j$ and time level $n$, the linear programming problem decouples with respect to the random index $i$. Therefore, in what follows, we describe the solver for one fixed $i$. We use the exponential parametrization
$$
\mu_{i,\ell}=e^{\xi_{i,\ell}},\qquad \ell=1,\ldots,N_{\bmu},
$$
which enforces the positivity of the discrete density values. We also introduce the corresponding discrete masses
$$
q_{i,\ell}:=\Delta V e^{\xi_{i,\ell}},\qquad \ell=1,\ldots,N_{\bmu}.
$$
Then the normalization and moment constraints in \eref{linprog-d}--\eref{linprog-e} become
$$
\sum_{\ell=1}^{N_{\bmu}}q_{i,\ell}=1,\qquad\sum_{\ell=1}^{N_{\bmu}}\bm z_\ell q_{i,\ell}=\bmu_{i,j}^n.
$$
For fixed $i$, the positive constant $\Delta\xi p_0$ in the objective functional does not affect the minimizer. Hence, in terms of the mass variables, the closure problem is equivalently written as
\begin{equation}
\min_{\bm q_i\geq \bm 0}\bm c^\top\bm q_i\qquad\mbox{\rm subject to}\qquad B\bm q_i=\bm b_{i,j}^n,
\label{D.1}
\end{equation}
where
$$
\bm q_i=(q_{i,1},\ldots,q_{i,N_{\bmu}})^\top,\qquad\bm c=(\eta(\bm z_1),\ldots,\eta(\bm z_{N_{\bmu}}))^\top,
$$
and
$$
B=\begin{pmatrix}
  1 & 1 & \cdots & 1\\
  \bm z_1 & \bm z_2 & \cdots & \bm z_{N_{\bmu}}
  \end{pmatrix}
  \in\mathbb{R}^{(d+1)\times N_{\bmu}},
  \qquad
  \bm b_{i,j}^n=
  \begin{pmatrix}
  1\\
  \bmu_{i,j}^n
  \end{pmatrix}
  \in\mathbb{R}^{d+1}.
$$
Here, the last $d$ rows of $B$ are understood componentwise. More precisely, if $\bm z_\ell=(z_\ell^{(1)},\ldots,z_\ell^{(d)})^\top$, then the $\ell$-th column
of $B$ is
$$
\begin{pmatrix}
1\\
  z_\ell^{(1)}\\
  \vdots\\
  z_\ell^{(d)}
  \end{pmatrix}.
$$

We now write the KKT conditions for \eref{D.1}. Let $\bm y_i\in\mathbb{R}^{d+1}$ be the Lagrange multiplier corresponding to the
equality constraint, and let
$$
\bm s_i=(s_{i,1},\ldots,s_{i,N_{\bmu}})^\top\in\mathbb{R}^{N_{\bmu}}
$$
be the nonnegative dual slack variable associated with the inequality constraint $\bm q_i\geq\bm 0$. The KKT conditions are
\begin{equation}
\begin{aligned}
  &B\bm q_i=\bm b_{i,j}^n,\\
  &\bm c-B^\top\bm y_i-\bm s_i=\bm 0,\\
  &\bm q_i\geq\bm 0,\qquad \bm s_i\geq\bm 0,\\
  &q_{i,\ell}s_{i,\ell}=0,
  \qquad \ell=1,\ldots,N_{\bmu}.
\end{aligned}
\label{D.2}
\end{equation}
The complementarity condition in \eref{D.2} means that, at the solution, a positive mass $q_{i,\ell}$ corresponds to a zero slack variable $s_{i,\ell}$, while a positive slack variable indicates that the corresponding mass is zero.

To preserve positivity during the iteration, we use the exponential parametrization
$q_{i,\ell}=\Delta V e^{\xi_{i,\ell}}$ and solve a perturbed KKT system by replacing
the exact complementarity condition in \eref{D.2} with
$$
q_{i,\ell}s_{i,\ell}=\tau,\qquad \ell=1,\ldots,N_{\bmu},
$$
where $\tau>0$ is a barrier parameter, which is gradually reduced to zero. This gives the nonlinear system
\begin{equation}
\begin{aligned}
  &B\bm q_i-\bm b_{i,j}^n=\bm 0,\\
  &\bm c-B^\top\bm y_i-\bm s_i=\bm 0,\\
  &Q_iS_i\bm 1-\tau\bm 1=\bm 0,
\end{aligned}
\label{D.3}
\end{equation}
where
$$
Q_i=\operatorname{diag}(q_{i,1},\ldots,q_{i,N_{\bmu}}),\qquad S_i=\operatorname{diag}(s_{i,1},\ldots,s_{i,N_{\bmu}}).
$$
The nonlinear system \eref{D.3} is solved by a damped Newton iteration in the variables $(\bm\xi_i,\bm y_i,\bm s_i)$, where
$$
\bm\xi_i=(\xi_{i,1},\ldots,\xi_{i,N_{\bmu}})^\top,\qquad q_{i,\ell}=\Delta V e^{\xi_{i,\ell}}.
$$
The damping parameter is chosen so that $\bm s_i$ remains positive during the iteration. The positivity of $\bm q_i$, and hence that of the discrete density values $\mu_{i,\ell}$, is automatically preserved by the exponential parametrization.

The KKT-based solver is an interior-point-type approximation of the LP solution for $\lambda_F=1$. As the barrier parameter \(\tau\to 0\), the perturbed KKT solution approaches a KKT point of the original linear programming problem. After convergence for a sufficiently small $\tau$, we set
$$
\mu_{i,j,\ell}^{*,n}=e^{\xi_{i,j,\ell}^{*,n}},\qquad \ell=1,\ldots,N_{\bmu}.
$$
This procedure is applied independently for $i=1,\ldots,N_\xi$ and for each spatial cell $j$ and time level $n$.


\end{document}